\documentclass[10pt]{amsart}
\usepackage{amssymb}
\usepackage{amsmath}
\usepackage{amsthm}
\usepackage[all]{xy}
\usepackage{color}

\newtheorem{theorem}{Theorem}
\newtheorem{proposition}[theorem]{Proposition}
\newtheorem{corollary}[theorem]{Corollary}
\newtheorem{lemma}[theorem]{Lemma}
\theoremstyle{definition}

\theoremstyle{remark}
\newtheorem{remark}[theorem]{Remark}

\newcommand{\Ksch}[1]{{{#1}}}
\newcommand{\Rsch}[1]{{\underline{#1}}}
\newcommand{\ksch}[1]{{{#1}}}
\newcommand{\quo}[1]{{\bar{\ksch{#1}}}}
\newcommand{\ac}[1]{{\operatorname{ac}(#1)}}
\newcommand{\g}{{\mathfrak{g}}}
\newcommand{\fsl}{{\mathfrak{sl}}}
\newcommand{\SL}{{\operatorname{SL}}}
\newcommand{\GL}{{\operatorname{GL}}}
\newcommand{\cInd}{{\operatorname{cInd}}}
\newcommand{\cRes}{{\operatorname{cRes}}}

\newcommand{\sgn}{{\operatorname{sgn}}}
\newcommand{\trace}{{\operatorname{trace}\,}}

\newcommand{\ceq}{{\, := \,}}
\newcommand{\tq}{{\, \vert \, }}
\newcommand{\ie}{{\it i.e.}}
\newcommand{\mes}{{\operatorname{m}}}
\newcommand{\Ad}{{\operatorname{Ad}}}

\newcommand{\ind}{{\operatorname{ind}}}
\newcommand{\EE}{{\bar\Q_\ell}}
\newcommand{\Kq}{{\mathbb{K}}}
\newcommand{\Rq}{{\mathcal{O}_{\mathbb{K}}}}
\newcommand{\Pq}[1]{{\mathfrak{p}_{\mathbb{K}}^{#1}}}
\newcommand{\Uq}{{\mathcal{O}_{\mathbb{K}}^*}}
\newcommand{\uq}{{\F_q^\times}}
\newcommand{\kq}{{\F_q}}
\newcommand{\Kalg}{{\bar{\mathbb{K}}}}
\newcommand{\Knr}{{\mathbb{K}_{nr}}}

\newcommand{\knr}{{\bar{\F}_q}}
\newcommand{\foi}{{{\hat\mu}}}
\newcommand{\mutmut}{{\textit{mut.mut.}}}
\newcommand{\Gal}[2]{{\operatorname{Gal}(#1 / #2)}}
\newcommand{\hecke}[1]{{\mathcal{H}(#1)}}
\newcommand{\vol}{{\operatorname{vol}}}

\newcommand{\AND}{{\ \wedge\ }}

\newcommand{\A}{{\mathbb{A}}}

\newcommand{\oi}{{{\mu}}}
\newcommand{\eps}{\varepsilon}
\newcommand{\ord}{\text{ord}}
\newcommand{\ring}[1]{\mathbb{#1}}
 
\newcommand{\Q}{\ring{Q}}
\newcommand{\Z}{\ring{Z}} 

\newcommand{\F}{\ring{F}} 

\newcommand{\N}{\ring{N}}

\newcommand{\frob}{{{\operatorname{TrFrob}}}}
\newcommand{\lef}{\ring{L}}
\newcommand{\eg}{{\it{e.g.},}}
\newcommand{\mot}{{{\text {Mot}}}}

\newcommand{\GKq}{{G}}
\newcommand{\gKq}{{\g}}
\newcommand{\qGkq}{{\quo{G}}}
\newcommand{\qgkq}{{\quo{\g}}}
\newcommand{\lb}{\left[}
\newcommand{\rb}{\right]}
\newcommand{\h}{{\mathfrak{h}}}
\newcommand{\gauss}{{\gamma}}
\newcommand{\st}{{\operatorname{st}}}
\newcommand{\cay}{{\operatorname{cay}}}
\newcommand{\Lie}{{\operatorname{Lie}}}
\newcommand{\hKq}{{\mathfrak{h}}}
\newcommand{\U}{{\operatorname{U}}}

\begin{document}

\title{Motivic proof of a character formula for SL(2)}

\date{\today}

\author{Clifton~Cunningham}
\address{Department of Mathematics, University of Calgary}
\email{cunning@math.ucalgary.ca}

\author{Julia~Gordon}
\address{Department of Mathematics, University of British Columbia}
\email{gor@math.ubc.ca}

\subjclass{Primary: 22E50 (Representations of Lie and linear algebraic groups over local fields); Secondary: 03C10 (Quantifier elimination, model completeness and related topics)}

\keywords{Motivic integration; supercuspidal representations; characters; orbital integrals}

\thanks{The first author thanks the Institute des Hautes \'Etudes Scientifique and the second thanks the Institute for Advanced Study for hospitality while parts of this article were being written. Both authors thank Loren Spice for helpful conversations. We also thank the first referee of this paper for a close reading and an important correction.}

\begin{abstract}
This paper provides a proof of a $p$-adic character formula by means of motivic integration.  
We use motivic integration to produce {virtual Chow} motives {that control} 
the values of the characters of all depth zero supercuspidal representations on all topologically unipotent elements of $p$-adic $\textrm{SL}(2)$; 
likewise, we find motives for the values of the Fourier transform of all regular elliptic orbital integrals having minimal non-negative depth in their own Cartan subalgebra, 
on all topologically nilpotent elements of $p$-adic $\mathfrak{sl}(2)$. 
We then find identities in the ring of virtual Chow motives over $\mathbb{Q}$ which relate these two classes of motives. 
These identities provide explicit expressions for the values of characters of all depth zero supercuspidal representations of $p$-adic $\textrm{SL}(2)$ as linear combinations of Fourier transforms of semi-simple orbital integrals, thus providing a proof of a $p$-adic character formula.
\end{abstract}

\maketitle

\section*{Introduction}
The representation theory of groups over the finite fields is stated entirely in terms of algebraic geometry by the work of Deligne-Lusztig \cite{DL} and Lusztig \cite{La}. Naturally, it is desirable to ``lift'' the geometric constructions to the $p$-adic groups. Motivic integration offers an approach to this problem which is very different from all previous ones; although it is not aiming at a complete understanding of the geometry underlying various objects of representation theory, it could in principle provide an algorithm for computing them in each individual case. This paper is essentially an experiment illustrating this point. 

In this paper we work with the group ${\SL}(2)$ over a $p$-adic field for the simple reason that the geometric objects that arise as a result of computing the ``motivic'' volumes turn out to be easily computable by hand; in fact, all the ones that are non-trivial turn out to be conics. For bigger groups both the algorithms and the results  are of course much more complicated, but we believe that still there are computable geometric objects responsible for the character values in general, 
and it will be the subject of a future article. 
{ There is already substantial evidence that the objects arising in harmonic analysis on reductive $p$-adic groups are ``motivic'' (in particular, computable). 
The most spectacular result in 
this direction appears in a very recent article \cite{CHL} that shows that the 
orbital integrals arising in the Fundamental Lemma are motivic. Earlier papers on this 
topic  include \cite{CH} that provides an explicit description of a certain class of 
semisimple orbital integrals, and \cite{G}, where it is proved that the 
values of characters of 
depth zero representations are motivic near the identity.}

The virtual Chow motives we calculate appear in the context of the following expansion for the character of any depth zero supercuspidal representation of $p$-adic ${ \SL}(2)$. (The precise statement is Theorem~\ref{thm:C}.)
Let $p$ be an odd prime, let $\Kq$ be a $p$-adic field with residue field $\kq$ (so $q$ is a power of $p$) and let $\GKq = { \SL}(2,\Kq)$.
{We use the modified Cayley transform 
$\cay(Y)=(1+Y/2)(1-Y/2)^{-1}$ to pass between the topologically nilpotent elements in the Lie algebra and the topologically unipotent elements in the group.}
We show that there is a finite set of regular elliptic orbits represented by the elements $X_z$ in the Lie algebra $\gKq \ceq { \fsl(2,}\Kq)$, each having minimal non-negative depth in its Cartan subalgebra, with the following property.
Let $\pi$ be an arbitrary depth zero supercuspidal representation of $\GKq$ and let $\Theta_\pi$ be the distribution character of $\pi$ in the sense of Harish-Chandra.
Then there are { rational numbers} $c_z(\pi)$ such that
\begin{equation}\label{eqn:C0}
 	\Theta_\pi\left(\cay (Y)\right)
 	= 
	\sum_{z} c_z(\pi)\ \foi_{X_z}(Y),
\end{equation}
for all regular, topologically nilpotent elements $Y$ of $\gKq$.  (Here, and $\foi_{X_z}$ denotes the generalized function on $\gKq$ corresponding to the Fourier transform of the orbital integral at $X_z$.) The coefficients $c_z(\pi)$ are given in Table~\ref{table:C}, from which one sees that {, up to sign (given by a quadratic character at $-1$),} $c_z(\pi)$ is in fact a rational function in $q$ (the cardinality of the residue field of $\Kq$) with integer coefficients { which are independent of $q$}; of course, this is a meaningful observation only if one clarifies in what sense the left- and right-hand sides of equation (\ref{eqn:C0}) may be viewed as functions in $q$, as we shall do later in the paper. We shall refer to equation (\ref{eqn:C0}) as a \emph{semi-simple character expansion}.
 The question of uniqueness of the coefficients is discussed in Section~\ref{subsection: comments}.

It is fair to say that the existence of a semi-simple character expansion is well-known
in one form or another, as is the fact that it extends,
\mutmut, to a much larger class of representations and groups. The
novelty of this paper lies in our proof.
We use motivic integration to ``separate'' the character and each of  the
orbital integrals into a linear combination of the values of the
corresponding function defined over the residue field, with
coefficients that are virtual Chow motives. Then we can see directly
that on each side of the semi-simple character expansion 
we have the combination of the same
values with  the same coefficients.
This approach clearly shows two ways in which algebraic geometry
appears in the values of the characters and orbital integrals: the
finite field function is in fact the characteristic function of a
character sheaf; and the ``motivic'' coefficients come from the
process of inflation and induction that connects our representation
with the representation of the group over the finite field. 

As a consequence of this perspective, we find much more than the
coefficients $c_z(\pi)$ promised above. In fact, we produce expressions for the values of the characters of all depth zero supercuspidal representations on all topologically unipotent elements of $\GKq$. 
Likewise, we produce expressions for the values of the Fourier transform of all regular elliptic orbital integrals having minimal non-negative depth 
in their own Cartan subalgebra, on all topologically nilpotent elements of $\gKq$. Comparing these leads to Table~\ref{table:C} and the proof of equation~(\ref{eqn:C0}).
{ The character tables for $\text{SL}(2, \Kq)$ were computed by Sally and Shalika in \cite{SSh}. Since they use a different construction of the supercuspidal representations  
of $\text{SL(2)}$ than induction from a compact subgroup used here, it is difficult to 
match our calculations with their classical calculation before we see the result. 
However, once we have the character values, it becomes easy to find them in the
character tables computed by Sally and Shalika.}
{This can, in fact, be used to match the representations of depth
zero obtained by the modern construction with the equivalent representations appearing in those classical tables.}
 
\tableofcontents

\section{Basic Notions}\label{section: basic notions}

Throughout this paper, $\Kq$ denotes a $p$-adic field; by this we mean { that}
$\Kq$ is a field equipped with a non-archimedean valuation such that { it} 
is complete and locally compact with respect to the topology
determined by the norm, and such that the residue field 
(which is necessarily finite)
has characteristic $p$. 
Notice that we do not put any
condition on the characteristic of $\Kq$. As is well-known, any such
$\Kq$ is a finite extension of $\Q_p$ or of $\F_p((t))$, { where $t$ is 
transcendental over $\F_p$}.  
The ring of integers of $\Kq$ will be denoted $\Rq$ and the maximal
ideal in $\Rq$ will be denoted $\Pq{}$.  The residue field will be
denoted $\kq$. 
We reserve $q$ for the cardinality of $\kq$, and $p$ for the
characteristic of $\kq$, so that $q$ is a power of $p$.

Next, we fix a prime $\ell$ different from $p$ and an algebraic closure $\EE$ of $\Q_\ell$. We will shortly assume $p\ne 2$.
Henceforth, by `representation' we mean a representation in a vector space over $\EE$.

In Sections~\ref{subsection:representations} and \ref{subsection:good} we make brief reference to the Bruhat-Tits building and Moy-Prasad filtrations for $\GKq$ (see \cite{MP}). In those sections of the paper, we will use the notation of Moy-Prasad.
In particular, for each pair
$(x,r)$, where $x$ is a point and $r$ is a non-negative real number as
above, $\GKq_{x,r}$ is a compact open subgroup of $\GKq$; when $r=0$, this is often abbreviated to $\GKq_{x}$. For each pair $(x,r)$ as above, let $\Rsch{G}_{x,r}$ be the smooth integral model for $G_\Kq$ introduced in \cite{Y}; 
the set of integral points in $\Rsch{G}_{x,r}$ coincides with $\GKq_{x,r}$. We write $\quo{G}_{x,r}$ for the reductive quotient of the special fibre of $\Rsch{G}_{x,r}$. Analogous notions apply to the Lie algebra $\gKq$: for any point $x$ and any real number $r$, $\gKq_{x,r}$ denotes the Moy-Prasad lattice in $\gKq$ and $\Rsch{\g}_{x,r}$ denotes the corresponding integral model for $\g_\Kq$ while $\quo{\g}_{x,r}$ denotes the corresponding Lie algebra scheme over $\kq$. When we use $\gKq_{x,r}$ and $\GKq_{x,r}$, we will often write these down explicitly, with the hope that readers unfamiliar with Moy-Prasad filtrations should have little difficulty with the essential features of this paper.

Note that we have not yet chosen a uniformizer for $\Kq$, and that none of the constructions above require such a choice.

\subsection{Depth-zero representations}\label{subsection:representations}

In this section we remind the reader how to construct all depth zero supercuspidal irreducible representations { of $\GKq$}, up to equivalence. 

Let $(\pi,V)$ be a representation of $\GKq$. For each point $x$ in the Bruhat-Tits building for $\GKq$, let $V_{x}$ denote the subspace consisting of those $v\in V$ such that $\pi(k)v = v$ for each $k\in \GKq_{x,0^+}$. Then $\GKq_{x,0}$ acts on $V_x$, and the resulting representation is denoted $\pi_x$. By definition, the representation $\pi$ has \emph{depth zero} if $V_x$ is non-trivial for some point $x$ in the Bruhat-Tits building for $\GKq$. 

{ The canonical map  $\Rq\to \kq$ will be denoted by $x\mapsto\bar{x}$ and, with $x$ and $r$ as above, we write $\rho_{x,r} : \gKq_{x,r} \to \quo{\g}_{x,r}$ for the quotient map with kernel equal to the pro-nilpotent radical $\gKq_{x,r^+}$ of $\gKq_{x,r}$.}

For each depth zero $\pi$ and each point $x$ in the Bruhat-Tits building for $\GKq$, the representation $\pi_x$ factors through $\rho_{x,0} : \GKq_{x,0} \to \qGkq_{x,0}$. Let $\quo{\pi}_x$ denote the representation of $\qGkq_{x,0}$ on $V_x$ such that $\pi_x(k) = (\quo{\pi}_x \circ \rho_{x,0})(k)$ for all $k\in \GKq_x$; then $\quo{\pi}_x$ is called the representation of $\qGkq_{x,0}$ obtained by {compact restriction} and is also denoted $\cRes^{\GKq}_{\GKq_{x,0}}(\pi,V)$.

From the other direction, let $x$ be a point in the Bruhat-Tits building for
$\GKq$ and let $(\sigma,W)$ be a representation of $\qGkq_x$,
and let $(\rho_{x,0}^*\sigma,W)$ be the representation of $\GKq_x$ defined by $(\rho_{x,0}^*\sigma)(k) = \sigma(\rho_{x,0}(k))$, for each $k\in \GKq_x$. The right-regular representation of $\GKq$ on the space of compactly supported functions $f : \GKq \to W$ such that $f(kg) = (\rho_{x,0}^*\sigma)(k) f(g)$ for all $g\in \GKq$ and all $k\in \GKq_x$ is called the representation of $\GKq$ obtained from $(\sigma,W)$ (or $(\rho_{x,0}^*\sigma,W)$) by {compact induction} and denoted $\cInd^{\GKq}_{\GKq_x} (\sigma,W)$.

Although $(\cRes^{\GKq}_{\GKq_x},\cInd^{\GKq}_{\GKq_x})$ is an adjoint pair of functors (see \cite{V}), one must be careful: even if $\sigma$ is a cuspidal irreducible representation of $\qGkq_x$, it does not follow in general that $\cInd^{\GKq}_{\GKq_x}\sigma$ is an admissible representation of $\GKq$, let alone supercuspidal. To clarify matters somewhat, we have the following result, which we paraphrase from independent results by Lawrence Morris and Moy-Prasad. For each point $x$ in the Bruhat-Tits building for $\GKq$, the compact restriction functor $\cRes^{\GKq}_{\GKq_x}$ restricts to a surjection from supercuspidal irreducible representations of $\GKq$ to cuspidal irreducible representations of $\qGkq_x$. Moreover, if $\pi$ is a depth zero supercuspidal representation of $\GKq$ then there is a {vertex} $x$ in the Bruhat-Tits building for $\GKq$ such that $\cInd^{\GKq}_{\GKq_x}(\cRes^{\GKq}_{\GKq_x} \pi)$ is equivalent to  $\pi$.

There are exactly two $\GKq$-orbits of vertices in the Bruhat-Tits building for $\GKq$; let $(0)$ and $(1)$ be the adjacent vertices corresponding to the maximal compact subgroups below.
\begin{equation}\label{eqn:vertices}
\begin{aligned}
\GKq_{(0)} &\ceq \left\{ \begin{bmatrix} a & b \\ c & d \end{bmatrix} \Big\vert a, b, c, d\in \Rq;\ ad-bc =1 \right\}\\
\GKq_{(1)} &\ceq \left\{ \begin{bmatrix} a & b \\ c & d \end{bmatrix} \Big\vert a, d\in \Rq, b\in\Pq{-1}, c\in \Pq{1};\ ad-bc =1 \right\}.
\end{aligned}
\end{equation}
The special fibres of the integral models $\Rsch{G}_{(0)}$ and $\Rsch{G}_{(1)}$ are both ${ \SL}(2)_\kq$; it follows that $\quo{G}_{(0)}$ and $\quo{G}_{(1)}$ are both ${ \SL}(2,\kq)$. Therefore, up to equivalence, all depth zero supercuspidal irreducible representations of $\GKq$ are produced by compact induction from cuspidal irreducible representations of ${ \SL}(2,\kq)$. Accordingly, to enumerate all depth zero supercuspidal irreducible representations of $\GKq$ it is necessary to recall the construction of cuspidal irreducible representations of ${ \SL}(2,\kq)$. This we will do in Section~\ref{subsection: cuspidal}. In the meantime, we record the result: each depth zero supercuspidal irreducible representation of $\GKq$ is equivalent to one in Table~\ref{table:rep}, where all terms are defined in Section~\ref{subsection: cuspidal}.

\begin{table}[htbp]  \caption{The representations $\pi$ appearing in Theorem~\ref{thm:C}}
  \centering
  \begin{tabular}{@{} |c|c| @{}}
    \hline
&\\
${ \pi(0,\theta)} \ceq \cInd^{\GKq}_{\GKq_{(0)}}(\sigma_\theta)$
	&	${ \pi(1,\theta)} \ceq \cInd^{\GKq}_{\GKq_{(1)}}(\sigma_\theta)$\\
&\\
\hline
&\\
${ \pi(0,+)} \ceq \cInd^{\GKq}_{\GKq_{(0)}}(\sigma_+)$
	&	${ \pi(1,+)} \ceq \cInd^{\GKq}_{\GKq_{(1)}}(\sigma_+)$\\
&\\
\hline
&\\
${ \pi(0,-)} \ceq \cInd^{\GKq}_{\GKq_{(0)}}(\sigma_-)$
	&	${ \pi(1,-)} \ceq \cInd^{\GKq}_{\GKq_{(1)}}(\sigma_-)$\\
&\\
\hline  
\end{tabular}
  \label{table:rep}
\end{table}

\subsection{Gauss sums}\label{subsection: gauss}

Before moving on to a review of the cuspidal irreducible representations of ${ \SL}(2,\kq)$, we say a few words about Gauss sums, which are the prototype for much that follows. Let $\sgn:\kq \to \EE$ be the quadratic character of $\uq$ extended by zero to $\kq$.  \emph{Arbitrarily but irrevocably}, we fix an additive character $\bar\psi : \kq \to \EE$. Let $\widehat\sgn$ denote the Fourier transform of $\sgn$ with respect to $\bar\psi$ in the following sense:
\begin{equation*}
\widehat{\sgn}(a) = \sum_{x\in \A^1(\kq)} {\sgn}(x) \bar\psi(ax).
\end{equation*}
(Note that this Fourier transform is not unitary.)
Consider the Gauss sums $\gauss_\pm : \kq\to \EE$ defined by
\begin{equation}\label{eqn: gauss}
\begin{aligned}
\gauss_+ (a) &\ceq \sum_{\{ x\in \A^1(\kq) \vert {\sgn}(x) = +1\}} \bar\psi(xa)\\
\gauss_- (a) &\ceq \sum_{\{ x\in \A^1(\kq) \vert {\sgn}(x) =-1\}} \bar\psi(xa).
\end{aligned}
\end{equation}
Then $\gauss_+ - \gauss_- = \widehat\sgn$.
Elementary arguments show that $\gauss_+ + \gauss_-  = -1$ and $(\gauss_+ - \gauss_-)^2 = q\, {\sgn}(-1)$. Fix a square-root $\sqrt{q}$ of $q$ in $\EE$. Then there is a {unique} square-root $\zeta\in \EE$ of ${\sgn}(-1)$ (determined by $\bar\psi$) such that
\begin{equation}\label{eqn:root-of-unity}
\widehat\sgn =  \sqrt{q} \zeta \sgn.
\end{equation}
Observe that $\sgn$ is an eigenvector for the Fourier transform with eigenvalue $\sqrt{q} \zeta$.

\begin{remark}\label{remark: zeta}
We wish to emphasize here that $\zeta\in \EE$ is determined by two choices: the square-root $\sqrt{q}$ of $q$ in $\EE$ {and} the additive character $\bar\psi: \kq \to \EE$. Also, although $\zeta$ is a fourth root-of-unity, it is not necessarily a {primitive} fourth root-of-unity; indeed, $\zeta^2 = {\sgn}(-1)$.
\end{remark}

\subsection{Cuspidal representations of ${ \SL}(2,\kq)$}\label{subsection: cuspidal}

Let $\ksch{T}$ denote the anisotropic torus of ${ \SL}(2)_\kq$ with $\kq$-rational points
\begin{equation}
\ksch{T}(\kq)
	= \left\{ \begin{bmatrix} x & y \\ \epsilon y & x\end{bmatrix}  \in { \SL}(2,\kq)
		\tq x^2- y^2 \epsilon =1\right\},
\end{equation}
where $\epsilon$ is a non-square in $\uq$. 

Let $R^{\ksch{G}}_{\ksch{T}}$ denote the Deligne-Lusztig induction functor of \cite{DL}; this takes virtual representations of $\ksch{T}(\kq)$ to virtual representations of ${ \SL}(2,\kq)$. If $\theta$ is non-trivial and the order of $\theta$ is not $2$ (so $\theta$ is in `general position') then $(-1)R^{\ksch{G}}_{\ksch{T}}(\theta)$ {is} an irreducible cuspidal representation of ${ \SL}(2,\kq)$, which we henceforth denote $\sigma_\theta$; in other words, as virtual representations,
\begin{equation}
	R^{\ksch{G}}_{\ksch{T}}(\theta) 
	= 
	- \sigma_\theta.
\end{equation}
Let $Q_\ksch{T}$ denote the restriction of $\trace R^{\ksch{G}}_{\ksch{T}}(\theta)$ to the set of unipotent elements of ${ \SL}(2,\kq)$, where $\theta$ is in general position; as the notation suggests, $Q_\ksch{T}$ is independent of $\theta$. This is the (classical) Green's polynomial for ${ \SL}(2,\kq)$. For future reference:
\begin{equation}\label{eqn:ff7.T}
\begin{aligned}
Q_{\ksch{T}}\left(g\right)
&= 
		\begin{cases} 
			1, 	& g\ne [\begin{smallmatrix} 1 & 0 \\ 0 & 1 \end{smallmatrix}] \AND \trace g =2 \\
			1-q , & g=[\begin{smallmatrix} 1 & 0 \\ 0 & 1 \end{smallmatrix}] \\
			0,	& \text{otherwise}.
		\end{cases}
\end{aligned}
\end{equation}
{ (See, for example, the appendix to \cite{DM}.)}

If $\theta_0$ is the quadratic character of $\ksch{T}(\kq)$, then the Deligne-Lusztig virtual representation $R_T^G(\theta_0)$ contains two irreducible cuspidal representations, comprising the Lusztig series for $(\ksch{T},\theta_0)$. It is well-known that the difference between the characters of these two representations is supported by the set of regular unipotent elements of ${ \SL}(2,\kq)$ (see Remark~\ref{remark: character sheaves}).
We label the representations in the Lusztig series for $(\ksch{T},\theta_0)$ by $\sigma_+$ and $\sigma_-$ and define
\begin{equation}\label{eqn:QGdefined}
Q_\ksch{G} \ceq \trace\sigma_+ - \trace\sigma_-
\end{equation}
in such a way that
\begin{equation}\label{eqn:zetarule}
Q_\ksch{G}\left(\begin{bmatrix} 1 & 1 \\ 0 & 1 \end{bmatrix}\right) = \sqrt{q}\zeta^3.
\end{equation}
Then,
\begin{equation}\label{eqn:ff7.G}
\begin{aligned}
Q_{\ksch{G}}\left( \begin{bmatrix} a & b \\ c & d \end{bmatrix}\right)
&=		\begin{cases} 
			\sqrt{q} \zeta^3 {\sgn}(b) , 	& a+d=2 \AND b\ne 0\\
			\sqrt{q} \zeta {\sgn}(c), 	& a+d=2 \AND c\ne 0\\
			0,				& \text{otherwise}.
		\end{cases}
\end{aligned}
\end{equation}
In particular $Q_\ksch{G}$ is supported by regular unipotent elements.
{ (For this calculation we refer readers to the lovely table from the appendix to \cite{DM} with the caveat that what they denote $\sqrt{q\, {\sgn}(-1)}$ is here denoted $\sqrt{q} \zeta^3$.) }

For a variety of reasons, it is the characters $\sigma_+$ and $\sigma_-$ which are the most interesting. 

Since $Q_{\ksch{G}}$ and $Q_{\ksch{T}}$ are supported by unipotent elements, 
for all primes $p$ except $2$, we can and shall often abuse notation by considering these as functions on the nilpotent elements of ${ \fsl(2,}\F_q)$ { by composing them with the modified Cayley transform $\cay(X) = (1+(X/2))(1-(X/2))^{-1}$. We then have:
\begin{equation*}
\begin{aligned}
Q_\ksch{G}\left([\begin{smallmatrix} 0 & 1 \\ 0 & 0 \end{smallmatrix}]\right)&= \sqrt{q} \zeta^3 \\
Q_\ksch{G}\left([\begin{smallmatrix} 0 & -1 \\ 0 & 0 \end{smallmatrix}]\right)&= \sqrt{q} \zeta \\
Q_\ksch{G}\left([\begin{smallmatrix} 0 & 0 \\ 0 & 0 \end{smallmatrix}]\right)&= 0;
\end{aligned}
\end{equation*}
and
\begin{equation*}
\begin{aligned}
Q_\ksch{T}\left([\begin{smallmatrix} 0 & 1 \\ 0 & 0 \end{smallmatrix}]\right)&= 1 \\
Q_\ksch{T}\left([\begin{smallmatrix} 0 & -1 \\ 0 & 0 \end{smallmatrix}]\right)&= 1 \\
Q_\ksch{T}\left([\begin{smallmatrix} 0 & 0 \\ 0 & 0 \end{smallmatrix}]\right)&= 1-q.
\end{aligned}
\end{equation*}
}

\begin{remark}\label{remark: signs}
We wish to emphasize that the choice of $\sqrt{q}$ and $\bar\psi$ determined $\zeta$, and that we labelled the representations in the Lusztig series $\{\sigma_+,\sigma_-\}$ for $(\ksch{T},\theta_0)$ precisely so that equation~(\ref{eqn:zetarule}) would be true.
\end{remark}

We close this section by recalling { one small consequence of Lusztig's celebrated work on character sheaves and representations of finite groups of Lie type} (see \cite{La}).
Let $\sigma$ be any cuspidal representation of ${ \SL}(2,\kq)$ and let $\chi_\sigma$ denote the restriction of $\trace\sigma$ to unipotent elements. Then there are unique {$a_Q(\sigma)\in \EE$} 
such that
\[
\chi_\sigma =  \sum_Q a_Q(\sigma) Q,
\]
where the sum is taken over the set $\{ Q_{\ksch{T}}, Q_{\ksch{G}}\}$ of (generalized) Green's polynomials. The values of $a_Q(\sigma)$ are given in Table~\ref{table:a}; { they will play a role in the calculations for the coefficients $c_z(\pi)$ appearing in Table~\ref{table:C}}.

\begin{table}[htbp]  \caption{$\chi_\sigma =  \sum_Q a_Q(\sigma) Q$.}  
  \centering
  \begin{tabular}{@{} | c | cc | @{}}
    \hline
    $a_Q(\sigma)$ & $Q= Q_{\ksch{T}}$ & $Q= Q_{\ksch{G}}$ \\ 
    \hline
	&&\\
    $\sigma =\sigma_\theta$ & $-1$ & $0$ \\ 
	&&\\
    $\sigma = \sigma_\pm$ & $-\frac{1}{2}$ & $\mp \frac{1}{2} $ \\ 
	&&\\
     \hline
  \end{tabular}
  \label{table:a}
\end{table}

\begin{remark}\label{remark: character sheaves}
The characters of the representations $\sigma_\theta$ and $\sigma_+$ and $\sigma_-$ introduced above are perhaps best understood in terms of characteristic functions of character sheaves.  Let $\theta$ be any character of $\ksch{T}(\kq)$ and let $\mathcal{L}_\theta$ be the corresponding Frobenius-stable Kummer local system on the \'etale site of $\ksch{T}_\knr$; in this case, the characteristic function $\chi_{\mathcal{L}_\theta} :\ksch{T}(\kq) \to \EE$ of $\mathcal{L}_\theta$ coincides with $\theta$ (see \cite[8.4]{La} for the definition of `characteristic function'). Let $\ind_{\ksch{T}}^{\ksch{G}}$ denote the cohomological induction functor of \cite{La}; this takes Frobenius-stable character sheaves for $\ksch{T}_\knr$ to Frobenius-stable perverse sheaves for $\ksch{G}_\knr$. If $\theta$ is non-trivial and not quadratic (\ie, in general position for ${ \SL}(2)$) then $\ind_{\ksch{T}}^{\ksch{G}}\mathcal{L}_\theta[1]$ is a Frobenius-stable character sheaf and the characteristic function of this perverse sheaf is the character of $R^{\ksch{G}}_{\ksch{T}}(\theta)$; in other words,
\[
\trace R^{\ksch{G}}_{\ksch{T}}(\theta) = \chi_{\ind_{\ksch{T}}^{\ksch{G}}\mathcal{L}_\theta[1]},
\]
when $\theta$ is in general position. (Note that
$\ind_{\ksch{T}}^{\ksch{G}}\mathcal{L}_\theta[1]$ is not a
\emph{cuspidal} character sheaf, even though $\trace
R^{\ksch{G}}_{\ksch{T}}(\theta)$ is a cuspidal function.) On the other
hand, $\ind_{\ksch{T}}^{\ksch{G}}\mathcal{L}_{\theta_0}[1]$ is not a
character sheaf; rather, it is a direct sum of character sheaves. The
algebraic group $\ksch{G}_\knr$ admits two cuspidal character sheaves:
one, denoted $C_0$, is unipotent (\ie, supported by the unipotent cone in
${ \SL}(2,\knr)$) while the other, denoted $C_1$, is supported by $-1$
times the unipotent cone. Comparing with the notation above, we have $Q_\ksch{G} = \chi_{C_0}$. In this paper we are chiefly concerned with the characters of cuspidal representations of ${ \SL}(2,\kq)$ when restricted to unipotent elements in ${ \SL}(2,\kq)$; the vector space spanned by characters of cuspidal representations restricted to unipotent elements is two-dimensional and a basis for this space is given by $Q_{\ksch{G}}$ (the restriction of $\chi_{C_0}$ to the unipotent cone) and $Q_{\ksch{T}}$ (the restriction of $\chi_{C_\theta}$ to the unipotent cone) where $\theta$ is any fixed character in general position. 
\end{remark}

\subsection{The elements $X_z$}\label{subsection:Xz}

Before specifying the orbits appearing in equation (\ref{eqn:C0}) we say a few words about Cartan subalgebras of $\gKq$, or equivalently, about conjugacy classes of forms of ${\GL}(1)_\Kq$ in ${ \SL}(2)_\Kq$. (These are the proper twisted-Levi subgroups of $\GKq$.) This is precisely the kernel ${\mathfrak{C}(\Kq)}$ of the map in Galois cohomology $\alpha_\Kq: H^1(\Kq, N) \to H^1(\Kq, G)$ which is induced by the inclusion $N \to G$, where $N$ denotes the normalizer of the split torus of diagonal matrices in $G$. To study this kernel, it is useful to first pass to the unramified closure $\Knr$ of $\Kq$ in $\Kalg$ and observe that ${\mathfrak{C}(\Knr)} = \ker\alpha_\Knr$ contains exactly one non-trivial element and this element corresponds to the conjugacy class of forms which split over a {(totally ramified)} 
quadratic extension of $\Knr$. 
Let $\tau_0$ be the trivial element in 
${\mathfrak{C}(\Knr)}$
and let $\tau_1$ be the non-trivial element in 
${\mathfrak{C}(\Knr)}$
. The inclusion $\Gal{\Kalg}{\Knr} \to \Gal{\Kalg}{\Kq}$ induces a surjection in Galois cohomology which in turn restricts to a surjection 
${\mathfrak{C}(\Kq) \to \mathfrak{C}(\Knr)}$
, as pictured below. 
\[
\xymatrix{
H^1(\Kq, G) \ar@{->>}[r] & H^1(\Knr, G) \\
\ar[u]^{\alpha_\Kq} H^1(\Kq, N) \ar@{->>}[r] & \ar[u]_{\alpha_\Knr} H^1(\Knr, N) \\
\ar[u] {\mathfrak{C}(\Kq)} \ar[r] & \ar[u] {\mathfrak{C}(\Knr)}\\
\ar[u] { 1} \ar[r] & \ar[u] { 1}
}
\]
The fibre above $\tau_0$ consists of three classes in $H^1(\Kq, N)$ and the fibre above $\tau_1$ consists of four classes in $H^1(\Kq, N)$ if ${\sgn}(-1)=1$ and two classes in $H^1(\Kq, N)$ if ${\sgn}(-1)=-1$, where $\sgn: \Kq\to \EE$ is the quadratic character of $\Uq$ extended by zero to all $\Kq$.

We now return to the orbits appearing in equation (\ref{eqn:C0}). In
order to define these orbits, we now pick { a uniformizer $\varpi$ for $\Kq$ and} a non-square unit $\eps$ in
$\Rq$. For reasons which will be  apparent later, we also now introduce a new parameter, $v$, taking values in $\Uq$ and define $X_z(v)$ in Table~\ref{table:orbits}. Representatives for the orbits appearing in equation (\ref{eqn:C0}) are given by $X_z\ceq X_z(1)$. Notice that each $X_z(v)\in \gKq$ above is \emph{good}, in the sense of \cite[\S 2.2]{A}.

\begin{table}[htbp]\caption{Elements $X_z(v)$. Orbits $X_z\ceq X_z(1)$ appear in Theorem~\ref{thm:C}.}
  \centering
  \begin{tabular}{@{} | c || c  c | c  c  c  c | @{}}
\hline
$z$ & ${s_1}$ & ${s_2}$ &  $t_0$ & $t_1$ & $t_2$ & $t_3$\\
\hline
&&&&&&\\
 $X_z(v)$ 
& {$\lb\begin{smallmatrix} 0 & v \\ \eps v & 0\end{smallmatrix}\rb$}
& { $\lb\begin{smallmatrix} 0 & \varpi^{-1} v \\ \eps\varpi v& 0\end{smallmatrix}\rb$}
& {$\lb\begin{smallmatrix} 0 & v \\ \varpi v & 0\end{smallmatrix}\rb$}
& {$\lb\begin{smallmatrix} 0 & \eps v\\ \varpi\eps^{-1} v & 0\end{smallmatrix}\rb$ }
 & {$\lb\begin{smallmatrix} 0 & v \\ \varpi\eps v& 0\end{smallmatrix}\rb$ }
& {$\lb\begin{smallmatrix} 0 & \eps v\\ \varpi v& 0\end{smallmatrix}\rb$}
\\
&&&&&&\\
\hline
 \end{tabular}
  \label{table:orbits}
\end{table}

Each of the $X_z(v)$ determines a compact Cartan subalgebra $\gKq$ (independent of $v\in \Uq$), and therefore a \emph{cocycle} representing a class in $\ker\alpha_\Kq$. Thus, the summation set in equation (\ref{eqn:C0}) is the set of cocycles $\{ s_1,s_2,t_0,t_1,t_2,t_3\} \subset Z^1(\Kq,N)$. However, some of these cocycles represent the same cohomology class in $H^1(\Kq,N)$, depending on the sign of $-1$ in $\Kq$; specifically, if ${\sgn}(-1)=-1$ then $t_0$ and $t_1$ represent to same cohomology class, and $t_2$ and $t_3$ also represent the same class in $H^1(\Kq, N)$. 

\begin{remark}\label{remark: eps}
Observe that picking { $\varpi$ and} $\eps$ is exactly equivalent to picking representative cocycles for the cohomology classes in $\ker\alpha_\Kq$. In fact, our choice for $s_1$ corresponds to our choice of $\eps$ { and our choice for $t_0$ corresponds to our choice of $\varpi$.}
\end{remark}

\subsection{Normalization of the measures}\label{subsection: measures}

We choose the Haar measure on ${ \SL}(2, \Kq)$ that coincides with the Serre-Oesterl\'e measure on ${ \SL}(2, \Rq)$; that is, our Haar measure for $\GKq$ is normalized so that the
volume of the maximal compact subgroup ${ \SL}(2, \Rq)$ is the cardinality
of ${ \SL}(2, \F_q)$, which is $q(q^2-1)$. In this setting, all the fibres
of the projection ${ \SL}(2,\Rq)\to { \SL}(2,\kq)$ have volume $1$. We denote
this measure on $\GKq$ by $\mes$. This choice determines a Haar measure on $\gKq$ such that volume of
of the kernel of the map ${ \fsl(2,}\Rq)\to { \fsl(2,}\kq)$ is also $1$; \ie, the volume
of ${ \fsl(2,}\Rq)$ equals the cardinality of ${ \fsl(2,} \kq)$.  We denote
this measure on  $\gKq$ by $\vol$. Observe that the volume
of ${ \fsl(2,}\Rq)$ is $\frac{q^2}{q^2-1}$ times that of ${ \SL}(2,\Rq)$.

{ With these choices, the formal degree of a representation $\pi$ produced by compact induction from a cuspidal irreducible representation $\sigma$ on the finite reductive quotient of a maximal parahoric subgroup of $\GKq$, coincides with the dimension $\trace\sigma(1)$ of $\sigma$.}

\subsection{Statement of the character formula}\label{subsection: coefficients}

We now state the character formula appearing in the introduction
more carefully. 
Fix measures for $\GKq$ and $\gKq$ as specified in Section~\ref{subsection: measures}. 

For any smooth representation $\pi$ of $\GKq$, let $\Theta_{\pi}$
denote the generalized function corresponding to the distribution
character $\pi$ with respect to the Haar measure on $\GKq$
above. For any regular semi-simple $X$ in $\gKq$ and any locally
constant compactly supported $f : \gKq \to \EE$, let $\oi_{X}(f)$
denote the orbital integral of $f$ at $X$. Fix an additive character $\psi$
of $\Kq$ with conductor $\Rq$ such that the induced additive character of the residue field $\kq$ is $\bar\psi$ (see Section~\ref{subsection: gauss}). For any $f$ as above, let $\hat{f}$ denote the Fourier transform of $f$ taken with respect to the Killing form $\langle Y,Z\rangle \ceq \trace(YZ)$, character $\psi$ and Haar measure on $\gKq$ specified in Section~\ref{subsection: measures}; thus,
\begin{equation*}
\hat{f}(Y) = \int_{\gKq} f(Z)\, \psi(\langle Y,Z\rangle)\, dZ.
\end{equation*}
Let $\foi_{X}(f)$ denote the orbital integral of $\hat{f}$ at $X$.

\begin{theorem}\label{thm:C}
Let $\Kq$ be a $p$-adic field with $p\ne 2$. For each
depth zero supercuspidal representation $\pi$ of $\GKq$ and for
each cocycle $z\in \{ s_1,s_2,t_0,t_1,t_2,t_3\}$ (see Section~\ref{subsection:Xz}) there is a
regular elliptic $X_z\in \gKq$ and $c_z(\pi)\in \Q$ such that
\begin{equation}\label{eqn:C}
 	\Theta_\pi(\cay^*f)
 	= 
	\sum_{z\in \{ s_1,s_2,t_0,t_1,t_2,t_3\}} c_z(\pi) \foi_{X_z}(f)
\end{equation}
for all { Schwartz functions $f$ supported by topologically nilpotent elements in $\gKq$}. 
{ The coefficients $c_z(\pi)$ are all rational numbers.} 
The coefficients $c_z(\pi)$ are given in Table~\ref{table:C}. 
\end{theorem}

This paper is devoted to understanding Theorem~\ref{thm:C} from the motivic perspective. 
To that end, we present a { brief} review of motivic integration in Section~\ref{subsection:motref}.

Our proof of Theorem~\ref{thm:C} is sketched in Section~\ref{section:outline} and executed in Section~\ref{section: motivic_proof}; we will calculate the { coefficients} $c_z(\pi)$ for each depth zero supercuspidal $\pi$, up to equivalence.  

{
\begin{remark}
More properties of the coefficients $c_z(\pi)$ are given in Section~\ref{subsection: comments}.
We show there that each coefficient may be interpreted as  a virtual Chow motive.
The issue of uniqueness of the coefficients is also addressed in Section~\ref{subsection: comments}.
\end{remark}
}

{
\begin{table}[htbp]  \caption{The coefficients $c_z(\pi)$ appearing in Theorem~\ref{thm:C}. (Here we write $\zeta^2$ for ${\sgn}(-1)$ in order to save space.) }
  \centering
  \begin{tabular}{@{} |c || c  c |  cc | c c | @{}}
    \hline
  &&&&&&\\
$c_z(\pi)$	& $z=s_1$ &  $z=s_2$ & $z=t_0$ & $z=t_1$ & $z=t_2$ & $z=t_3$ \\ 
 &&&&&&\\
    \hline\hline
 &&&&&&\\
${ \pi=\pi(0,\theta)}$ & $q-1$ &  $0$ & $0$ & $0$ & $0$ & $0$  \\
 &&&&&&\\
${ \pi= \pi(1,\theta)}$ & $0$  & $q-1$ &  $0$ & $0$ & $0$ & $0$ \\ 
 &&&&&&\\
 \hline\hline
 &&&&&&\\
${ \pi= \pi(0,+)}$ & $\frac{q-1}{2}$ &  $0$ &  $-\frac{q^2-1}{2^3q}\zeta^2$ & $+\frac{q^2-1}{2^3q}\zeta^2$ & $-\frac{q^2-1}{2^3q}\zeta^2$ & $+\frac{q^2-1}{2^3q}\zeta^2$ \\
 &&&&&&\\
${ \pi=\pi(0,-)}$ & $\frac{q-1}{2}$  & $0$ & $+\frac{q^2-1}{2^3q}\zeta^2$ & $-\frac{q^2-1}{2^3q}\zeta^2$ & $+\frac{q^2-1}{2^3q}\zeta^2$ & $-\frac{q^2-1}{2^3q}\zeta^2$ \\
 &&&&&&\\
${ \pi=\pi(1,+)}$ & $0$  &  $\frac{q-1}{2}$ &  $-\frac{q^2-1}{2^3q}\zeta^2$ & $+\frac{q^2-1}{2^3q}\zeta^2$ & $+\frac{q^2-1}{2^3q}\zeta^2$ & $-\frac{q^2-1}{2^3q}\zeta^2$ \\
 &&&&&&\\
${\pi=\pi(1,-)}$ & $0$  &  $\frac{q-1}{2}$ & $+\frac{q^2-1}{2^3q}\zeta^2$ & $-\frac{q^2-1}{2^3q}\zeta^2$ & $-\frac{q^2-1}{2^3q}\zeta^2$ & $+\frac{q^2-1}{2^3q}\zeta^2$ \\
 &&&&&&\\
\hline
  \end{tabular}
  \label{table:C}
\end{table}
}

\subsection{Motivic integration}\label{subsection:motref}

The goal of arithmetic motivic integration is essentially to reduce the
calculation of $p$-adic volumes to the computation of the
number of points on varieties over finite fields. In particular, 
suppose we are talking about subsets of an affine space defined over a
$p$-adic field. Consider the measurable subsets that can be described
in a language (of logic) that depends neither on $p$ nor on the choice of the
uniformizer of the valuation (one such language is the language of
Denef-Pas that is described in the next subsection).  Once a set is
described by a logical formula in such a language, it is possible to
associate a geometric object defined over $\Q$ with it, in such a way
that for almost all primes $p$, the volume of our set can be
recovered  from the number of points of the reduction of this object
over $\F_q$, where $q=p^r$ is the cardinality of the residue field of
the given local field.  

Motivic integration is based on the algorithm of elimination of
quantifiers (there is also a more recent version \cite{Cluckers.Loeser} 
where in some parts quantifier elimination is replaced with cell
decomposition, also an algorithmic procedure). 
Since at present motivic integration as an algorithm is not implemented, it has not previously been used for computing examples; nevertheless, it can be used to prove general existence results for the required
geometric objects (which already has implications, as explored in \cite{TH}, \cite{CH} and \cite{G}).

We refer to the original papers \cite{DL.ar}, \cite{Cluckers.Loeser}
and to the beautiful exposition \cite{Tom.whatis} for the detailed
description the concept of motivic integration, and the statements of
``comparison theorems'' that relate the classical $p$-adic volumes
with the motivic volumes. In the next few subsections we just give a list of the techniques
we use.

\subsubsection{The language of rings and the language of Denef-Pas}\label{subsub:language}
The language of rings is the language of logic such that its formulas
can be interpreted in any ring with identity; thus, any such ring is a
structure for this language. The language of rings has symbols for '$0$' and '$1$', symbols
for countably many variables, and the
symbols for the operations of addition '$+$' and multiplication
'$\times$'. A formula is a syntactically correct expression built
from finitely many of these symbols, and also parentheses '$()$', quantifiers 
'$\exists$', '$\forall$', and symbols of conjunction '$\wedge$',
disjunction '$\vee$', and negation '$\neg$'. We will usually use this
language to work with the residue field of our valued field. 
 
There is also a first-order language that is perfectly suited for
defining subsets of non-archimedean valued fields: it is the language
of Denef-Pas. Formulas in this language have variables of three {\it
  sorts}: the valued field sort, the residue field sort, and the value
sort. There also are symbols for the function '$\ord(\cdot)$' which takes variables of the valued field
sort to the variables of the value sort, and the function '$\ac{\cdot}$' that
takes the variables of the valued field sort to the variables of the
residue field sort.  The formulas are built using algebraic operations
on variables of the same sort, the symbols '$\ord(\cdot)$' and
'$\ac{\cdot}$' (which can be applied  to variables of the valued field
sort), quantifiers, and the symbols for the logical operations of
conjunction, disjunction and negation.

Every valued field with a choice of a uniformizer  $(\Kq,\varpi)$ is a
structure for the language of Denef-Pas; then the functions '$\ord$' and
'$\text{ac}$' match the usual valuation and ``angular component'' maps,
where the ``angular component'' $\ac{x}$ equals the first non-zero
coefficient of the $p$-adic expansion of $x$, \ie, if $x$ is a unit,
$\ac{x}=x\mod (\varpi)$; if $x$ is not a unit,
$\ac{x}=x\varpi^{-\ord(x)} \mod (\varpi)$.

A subset of an affine space over a local field $\Kq$ is called {\it
  definable} if it can be defined by a formula in the language of
Denef-Pas.
For a formula $\Phi$ with $m$ free variables of the residue field sort
  and no other free variables, given a local field with the
  ring of integers $\Rq$ and a 
choice of the uniformizer, we will denote by $Z(\Phi, \Rq)$ the
  subset of $\Rq^{\oplus m}$ defined by the formula $\Phi$.   

\subsubsection{The ring of virtual Chow motives}\label{subsubsection: ch.m.}
In all versions of motivic integration, the motivic volume takes
values in some ring built from the Grothendieck ring of the category
of algebraic varieties over the base field, which is $\Q$ in our case.
 
Strictly speaking, the so-called ring of virtual Chow motives is a
natural choice of the ring of values for the arithmetic motivic volume,
for reasons described in \cite{Tom.whatis}. 
We will not define Chow motives here (see \cite{Scholl} for a good 
introduction). As a first approximation, it is possible to think just
of formal linear combinations of isomorphism classes of varieties
with rational coefficients.

The motivic volume takes values in the ring $\mot$, which we will now define.
Let ${\mathcal Mot}_{\Q}$ be the category of Chow motives over $\Q$. 
We take its Grothendieck ring $K_0({\mathcal Mot}_{\Q})$, \ie, the ring
of formal linear combinations of the isomorphism classes of the Chow
motives, with natural relations (see \cite{Scholl}), and the   
product coming from tensor product. This ring has a unit $1$.
For every smooth projective variety $V$ there is a Chow motive 
that corresponds to $V$ in a natural way. 
It is a deep theorem that  this map from varieties to
motives extends to all (not just smooth projective) varieties, and 
induces  a homomorphism from 
the Grothendieck ring of the category of varieties 
$K_0(Var_{\Q})$ to $K_0({\mathcal Mot}_{\Q})$. Let us denote the image of
this map by $K_0({\mathcal Mot})^v$.   
The image of the class of the affine line under this map is usually
denoted by $\lef$ (called the Lefschetz motive); the point maps to $1$. 
Note that naturally, in $K_0({\mathcal Mot}_{\Q})$ we have
$\lef+1=[{\mathbb P}^1]$, where $[{\mathbb P}^1]$ is the class of the
Chow motive corresponding to the projective line. 

The ring $\mot$ is obtained from  $K_0({\mathcal Mot})^v\otimes_{\Z}\Q$ by
localizing at
$\lef$ followed by localizing further at the
set $\{\lef^{-i}-1\mid {i>0}\}$ (\ie, all these elements can be
formally inverted). 
Note that the localization at 
$\{\lef^{-i}-1\mid {i>0}\}$ (which is equivalent to adding in the sums
of geometric series with quotient $\lef^i$) replaces { the} completion that
was done in the original version \cite{DL.ar}. 
 
Let us also adopt the following convention: for a polynomial 
$F(q)\in \Z[q]$, we will denote by $[F(\lef)]$ the class of $F(\lef)$
in the  ring ${\mot}$. 

For every prime power $q$, there is an action of $\frob_q$ on the 
elements of $\mot$ that comes from the Frobenius action on the 
Chow motives. It can be thought of as a generalization of counting
points on a variety over $\kq$. The trace of the
Frobenius operator on a Chow motive is the alternating sum 
of the traces of Frobenius acting on its $\ell$-adic cohomology 
groups, and it is, {\it a priori}, 
an element of $\EE$ (see \cite[Section 3.3]{DL.ar}), but in fact this number 
lies in  $\Q$ for all elements of $\mot$ that arise from the motivic 
volumes    
(in particular, for us the choice of $\ell$ doesn't matter). 
It is the trace of 
Frobenius action that allows us to relate the values of the motivic 
measure (elements of $\mot$) with the usual $p$-adic volumes (rational numbers).

\subsubsection{The motivic volume}\label{subsection:comparison theorem}

 We cannot describe the construction of the motivic volume here, but
 we only need very simple examples in this paper, so we list just a
 few main principles and introduce the notation.
Following the pattern of \cite{Tom.whatis}, we start by declaring that
there is a map from formulas in the language of rings to the ring
 $\mot$ defined above. We denote the image of a formula $\phi$ under
 this map by $[\phi]$. We will use two properties of this map: 
\begin{itemize}
\item The motivic volume of a formula in the language of rings with no
  quantifiers that defines a smooth algebraic variety is the class of
  that variety.
\item If $\phi_1$ and  $\phi_2$ are equivalent ring formulas, then they
  have equal motivic volumes. If $\phi_1$ is an $n$-fold cover of
  $\phi_2$ (both the equivalence, and the term ``cover'' are understood
  in the sense of \cite{Tom.whatis}), then
  \[[\phi_2]=\frac1n[\phi_1].\]
\end{itemize}

Now let us say a few words about the construction of the motivic
volume for definable sets that is compatible with the
classical volume of $p$-adic sets. First, we need to talk about sets
in a way that does not  depend  on the field -- and that is done by
considering {\it definable} sets, as described above. 
In order to pass back and forth between sets and
formulas defining them, it is very convenient to use the notion of a 
{\it definable subassignment} introduced in \cite{DL.ar}. 

Let $h_{\A^m}$ be the functor of points of the affine space that takes
fields to sets ($K\mapsto \A^m(K)$). A {\it subassignment} of
$h_{\A^m}$ is simply a collection of subsets of $\A^m(K)$, one for
each field $K$. A formula $\Phi$ in the language of Denef-Pas with $m$ free
variables naturally defines a subassignment of $h_{\A^m}$ on {
a suitable category of local fields ({\it e.g.}, the category of all non-archimedean completions of a given global field):} 
for each local field $K$, take the subset of $\A^m(K)$ on which the formula
takes the value 'true'. Such subassignments are called {\it definable}.

The {\it motivic volume} is defined on definable subassignments and
 takes values in the ring of virtual Chow motives $\mot$. We denote
 the motivic volume by $\mu$, so that $\mu(\Phi)$ denotes the motivic
 volume of the subassignment defined by the formula $\Phi$. Sometimes,
 if we have a definable set $H$ (defined by some formula $\Phi$), we
 will write $\mu(H)$, meaning $\mu(\Phi)$.  

 The strategy of motivic integration is to replace the formula that defines our subassignment
 with an equivalent formula that has no quantifiers ranging over the
 valued field (the language of Denef-Pas admits quantifier
 elimination).  The next step is to ``approximate'' the new formula
 by the formulas in the language of rings whose variables range only
 over the finite field. Finally, ring formulas can be mapped to the
 ring $\mot$ defined above, as mentioned in the beginning of this subsection. 

The main result that makes motivic integration relevant for us is the 
comparison theorem \cite[Theorem 8.3.1, 8.3.2]{DL.ar}.   
Let $\Phi$ be a formula in the language of
Denef-Pas. The comparison theorems (one of which deals with 
the fields of characteristic zero, and the other one -- with function
fields), state that    
for all but finitely many primes $p$, for any local field $\Kq$ with
the ring of integers $\Rq$ and the residue field $\kq$ of cardinality 
$q=p^r$, the volume of the set $Z(\Phi,\Rq)$ equals the trace of
Frobenius action on $\mu(\Phi)$ (which is an element of $\mot$). 

In fact, it is possible to some extent to keep track of the ``bad
primes'' where the comparison theorem fails, and we will do this as we
compute the motivic volumes in our examples.
 
\subsubsection{Motivic integration with parameters}\label{subsub:parameters}

Here we quote one more technical notation from
\cite{Cluckers.Loeser} that will be used only in the proof of Lemma \ref{lemma:unram.elements}.
Cluckers and Loeser \cite{Cluckers.Loeser} have a more general
construction of the motivic volume: they define motivic
volume on definable subassignments of the functor of points of 
$\A_{K}^m\times \A_k^n\times \Z^r$, which is denoted by $h[m,n,r]$
(here $K$ stands for the valued field, and $k$ for its residue field),
which corresponds to considering families of definable subassignments
of $h_{\A_K^m}$ with parameters in $k$ and in $\Z^r$.

\subsection{Outline of the proof of Theorem~\ref{thm:C}}\label{section:outline}

Since the rest of the paper is devoted to the proof of Theorem~\ref{thm:C}, we describe the strategy before we get into the technicalities.

We need to prove an equality of two distributions represented by locally constant locally integrable functions on the set of regular topologically unipotent elements in the group $\GKq$.  
Each of these functions takes countably many values. We prove the equality by partitioning the domain into sets on which the both sides are constant, and showing equality on each of these sets. Moreover, the problem is in fact not just countable, it is in some sense finite. Each one of the two functions is obtained by a Frobenius-like formula from a function which takes finitely many values (explicitly, three distinct values for the function on the left-hand side and five distinct values for the function on the right-hand side). Motivic integration allows us to isolate these values and the coefficients at each one of them that are acquired in the Frobenius formula. Those latter coefficients also turn out to be computable. It is in this sense that the proof is better suited for a computer than for a human -- and that is the point of this paper, really.

In order to be more precise we must first say a few words about regular, topologically
unipotent elements of $\GKq$. To do this, we start by studying
regular, topologically nilpotent elements in the Lie algebra
$\gKq$; we will then invoke our assumption that $p$ is odd,
in which case these correspond exactly, 
via the { modified Cayley transform $\cay(X) = (1+(X/2))(1-(X/2))^{-1}$,
to regular, topologically unipotent elements of $\GKq$.
(In fact, $\cay$ establishes a bijection from topologically nilpotent elements of $\gKq$ to topologically unipotent elements of $\GKq$.)}

What we need, ideally, is a partition of the set of regular topologically
unipotent elements into {\it definable} sets such that all the functions
involved (the characters and Fourier transforms of the 
orbital integrals) would be constant on
them. We will not explicitly use any local constancy results (except in Section~\ref{subsection: comments}){;} instead,
we will first make the partition, and then see that it is the right
one.

Recall that at the moment, $\eps$ (a non-square unit) 
and $\varpi$ (the uniformizer of the valuation) are fixed.  One of the
points of motivic integration approach is to do the calculation
independently of these choices, but that requires some tricks
(discussed below), since they are both very deeply ingrained in all the
definitions: $\eps$ is linked with the cocycles $z$, and 
$\varpi$ plays a role in the construction of our representations.
For now, we write down explicitly seven Cartan subalgebras of $\gKq$ 
(one non-compact and six compact).
\begin{equation*}
\begin{aligned}
\mathfrak{h}_{{ s_0}} &= \left\{ \begin{bmatrix} x & 0 \\ 0 & -x \end{bmatrix} \ \Big\vert\  x\in \Kq \right\}
\end{aligned}
\end{equation*}

\begin{equation*}
\begin{aligned}
\mathfrak{h}_{{s_1}} &= \left\{ \begin{bmatrix} 0 & x \\ x \eps  & 0 \end{bmatrix} \ \Big\vert\  x\in \Kq \right\}\\
\mathfrak{h}_{{s_2}} &= \left\{ \begin{bmatrix} 0 & x \varpi^{-1} \\  x \eps \varpi  & 0 \end{bmatrix} \ \Big\vert\  x\in \Kq \right\}
\end{aligned}
\end{equation*}

\begin{equation*}
\begin{aligned}
\mathfrak{h}_{{t_0}} &= \left\{ \begin{bmatrix} 0 & x \\ x \varpi  & 0 \end{bmatrix} \ \Big\vert\  x\in \Kq \right\}\\
\mathfrak{h}_{{t_1}} &= \left\{ \begin{bmatrix} 0 & x \eps \\  x \varpi \eps^{-1} & 0 \end{bmatrix} \ \Big\vert\  x\in \Kq \right\}.
\end{aligned}
\end{equation*}

\begin{equation*}
\begin{aligned}
\mathfrak{h}_{{t_2}} &= \left\{ \begin{bmatrix} 0 & x \\ x \varpi  \eps & 0 \end{bmatrix} \ \Big\vert\  x\in \Kq \right\}\\
\mathfrak{h}_{{t_3}} &= \left\{ \begin{bmatrix} 0 & x \eps \\ x \varpi  & 0 \end{bmatrix} \ \Big\vert\  x\in \Kq \right\}
\end{aligned}
\end{equation*}

We have labelled our Cartan subalgebras by cocycles 
$z\in \{ { s_0, s_1, s_2, t_0, t_1 , t_2 , t_3}\}$, one from each 
cohomology class { in $\mathfrak{C}(\Kq)$ (see Section~\ref{subsection:Xz})} in the case ${\sgn}(-1)=1$. 
In the case ${\sgn}(-1)=-1$, 
the Cartan subalgebras labelled by ${t_0}$ and ${t_1}$ are conjugate;
and the ones labelled by ${t_2}$ and ${t_3}$ are conjugate (which means
that the corresponding cocycles represent the same cohomology
class). We will always consider each one of these subalgebras
separately, which means that we are doing one third more work than
necessary half the time. 

Our Cartan subalgebras are filtred by an index $n$ which is closely related to depth.
To write the filtration lattices explicitly, in Table~\ref{table:Yzn} we define, for each 
$z\in \{{ s_0,s_1,s_2,t_0,t_1,t_2,t_3}\}$, integer $n$ 
 and unit $u$ in $\Rq$, an element $Y_{z,n}(u)\in \h_z$.
Let  $\h_{z,n} \ceq \left\{ Y_{z,n}(u) \tq u\in \Uq  \right\}$. 
{ If $z$ is $s_0$, $s_1$ or $s_2$ then $\h_{z,n}$ is the set of elements of $\h_z$ with depth $n$;
however, if $z$ is $t_0$, $t_1$, $t_2$ or $t_3$ and $\h_{z,n}$ is the set of elements of $\h_z$ with depth $\frac{1}{2} + n$.}

\begin{table}[htbp]\caption{Elements $Y_{z,n}(u)$.}
  \centering
  \begin{tabular}{@{} | c || c  c  c | @{}}
\hline
$z$ & ${s_0}$ & ${s_1}$ & ${s_2}$ \\
\hline
&&&\\
$Y_{z,n}(u)$ 
& $\lb\begin{smallmatrix} \varpi^n u& 0 \\ 0 & -\varpi^n u \end{smallmatrix}\rb$ 
& $\lb\begin{smallmatrix} 0 & \varpi^n u\\ \eps \varpi^n u& 0\end{smallmatrix}\rb$ 
& $\lb\begin{smallmatrix} 0 & \varpi^{n-1} u\\ \eps\varpi^{n+1} u & 0\end{smallmatrix}\rb$\\
&&&\\
\hline
 \end{tabular}
   \begin{tabular}{@{} | c || c  c  c  c | @{}}
\hline
$z$ &  $t_0$ & $t_1$ & $t_2$ & $t_3$\\
\hline
&&&&\\
$Y_{z,n}(u)$ 
& $\lb\begin{smallmatrix} 0 & \varpi^n u\\ \varpi^{n+1} u& 0\end{smallmatrix}\rb$
& $\lb\begin{smallmatrix} 0 & \eps \varpi^n u \\ \eps^{-1}\varpi^{n+1} u  & 0\end{smallmatrix}\rb$ 
 & $\lb\begin{smallmatrix} 0 & \varpi^n u \\ \eps \varpi^{n+1}u& 0\end{smallmatrix}\rb$ 
& $\lb\begin{smallmatrix} 0 & \eps \varpi^n u\\ \varpi^{n+1} u& 0\end{smallmatrix}\rb$\\
&&&&\\
\hline
 \end{tabular}
  \label{table:Yzn}
\end{table}

In this paper we are particularly interested in understanding our distributions on 
{ topologically nilpotent elements of $\gKq$}.
In anticipation of the correct partition of regular, topologically nilpotent elements of
$\gKq$ we { now} define 
\[
\h_{z,n,+} \ceq \left\{ Y_{z,n}(u) \tq u\in \Uq
\AND {\sgn}(u)=+1 \right\},
\]
 and 
\[
\h_{z,n,-} \ceq \left\{ Y_{z,n}(u) \tq
u\in \Uq \AND {\sgn}(u)=-1  \right\},
\]
and denote by ${}^G\h_{z,n,\pm}$ the corresponding $\GKq$-invariant sets. 
{ If $z=s_0$, $s_1$ or $s_2$ then ${}^G\h_{z,n,\pm}$ contains topologically nilpotent elements if and only if $n\geq 1$;
if $z=t_0$, $t_1$, $t_2$ or $t_3$ then ${}^G\h_{z,n,\pm}$ contains topologically nilpotent elements if and only if $n\geq 0$.}

 Further, let  ${}^G\h_{z,n}$ ({ resp. ${}^G\h_{z,n,+}$,
 ${}^G\h_{z,n,-}$}) denote the smallest $\GKq$-invariant set containing
 $\h_{z,n}$, where $\GKq$ acts on $\gKq$ by adjoint action.  
The disjoint union of all the sets ${}^{G}\h_{z,n}$ 
(as $z$ ranges over 
{ the fixed set of representatives for $\mathfrak{C}(\Kq)$} and $n$ ranges over all
non-negative integers, with the exception of $n=0$ in ${}^G\h_{s_0,n}$, ${}^G\h_{s_1,n}$ and ${}^G\h_{s_1,n}$) 
coincides with the set of regular topologically nilpotent
elements in $\gKq$.

In fact, conjugation by $\GKq$ removes the ambiguity caused by the choice
 of specific cocycles (and therefore the choice of $\eps$): different
 representatives of the same cohomology class correspond to different,
 but conjugate, Cartan subalgebras. Hence, it would be more
 appropriate to label the sets ${}^G\h_{z,n}$ by the cohomology
 classes, and call them  ${}^G\h_{\tau,n}$, where $\tau$ is the cohomology class represented by $z$. We would like to stress
 that this produces a partition of the topologically nilpotent set
 that depends only on ${\sgn}(-1)$ (because the number of the distinct
 cohomology classes depends on this sign), and not on any choices.
However, we do need the specific representatives $Y_{z,n}(u)$ of each
 one of the sets ${}^G \h_{\tau,n}$ in order to carry out our calculations. 
Our choice of these representatives depends on $\eps$ and $\varpi$.
The strategy is to carry out all the explicit calculations first, get
 to the level of logical formulas, and at that stage $\varpi$ would
 disappear entirely. To remove the dependence of $\eps$, in
Section~\ref{W_Clambda} we turn it
 into a parameter ranging over the residue field, and then average
 over it, also showing that in fact we are averaging a constant function. 

Note that with the way we defined our subalgebras, in the case 
${\sgn}(-1)=-1$ when two pairs of them become conjugate, this conjugation
switches the $\pm$ sign: ${}^G\h_{t_0,n,+}$ corresponds to 
${}^G\h_{t_1,n,-}$, and ${}^G\h_{t_3,n,+}$ corresponds to ${}^G\h_{t_2,n,-}$
(we will see that this is reflected in our character table). 

On the group side, we write $\Gamma_{z,n}$ (resp. $\Gamma_{z,n,+}$,
$\Gamma_{z,n,-}$) for the image of $\h_{z,n}$ (resp. $\h_{z,n,+}$,
$\h_{z,n,-}$) under the { modified Cayley transform $\cay$,} when defined; we denote the corresponding $\GKq$-invariant sets by 
${}^G\Gamma_{z,n}$ (resp. ${}^G\Gamma_{z,n,+}$, ${}^G\Gamma_{z,n,-}$).
The characteristic functions of the sets ${}^G\Gamma_{z,n,\pm}$ will
be denoted by $f_{z,n,\pm}$ and the characteristic functions of the
sets ${}^G\h_{z,n,\pm}$ will be denoted by $\tilde f_{z,n,\pm}$
throughout the paper. 

Given a depth zero representation $\pi$ induced from a maximal compact
subgroup $G_x$, in
Section~\ref{section:frobenius characters} we
associate three virtual Chow motives (denoted $M_{z,n,\pm}^{x,0}$, $M_{z,n,\pm}^{x,1}$ and $M_{z,n,\pm}^{x,\eps}$) 
with each triple $(z,n,\pm)$ as above, so that 
the value of the distribution character of $\pi$ at
$f_{z,n,\pm}$
can be recovered
from this triple of virtual Chow motives  for almost all residual
characteristics $p$ (in fact, for all $p\neq 2$). Moreover, we will see that the only
way in which these motives depend on $\pi$ is through the compact
subgroup on which $\pi$ has non-trivial compact restriction (\ie, through the choice of the
vertex in the building $x=(0)$ or $x=(1)$, see Section~\ref{subsection:representations}). 

In Section~\ref{section:frobenius orbital}, we also associate with
each
$(z,n, \pm)$ five  
virtual Chow motives, so that all orbital integrals that appear in the
right-hand side of the semi-simple character expansion can be recovered
from these motives.

In Section~\ref{subsection: finite}, we put together all the results
about groups over the finite fields that are crucial for our
understanding of the $p$-adic ``lifts''.

Since both sides of semi-simple character expansion are invariant under
conjugation by $\GKq$, all we have to do is check the equality on each of the sets
$\Gamma_{z,n}$. This is done in Section~\ref{section: motivic_proof}. Due to a very mechanical nature of this proof, we do
not include all the details. In all proofs, we
include all the details  in one ``unramified'' case
$z=s_1$ and one ``ramified'' case $z=t_2$; in the other cases we just
indicate the differences and summarize the results in the tables.

\section{Motives corresponding to our characters}\label{section:frobenius characters}

Throughout Section~\ref{section:frobenius characters}, $x$ denotes either the standard vertex $(0)$ or the vertex $(1)$ in the Bruhat-Tits building for $\GKq$ (see equation~\ref{eqn:vertices}); we also reserve the symbol $y$ for the baricentre of the facet with boundary $\{ (0),(1)\}$.

\subsection{Consequences of the Frobenius formula for the character}\label{sub:tricks}

This section follows the method of expressing the character as a sum
over conjugacy classes in the reductive quotient that was used in
\cite{G}. 

Let $\pi$ be a supercuspidal depth zero representation of $\GKq$ and
let $\Theta_{\pi}$ be its distribution character in the sense of
Harish-Chandra. Let $f$ be a test function supported on some compact
subset $H$ of the set of regular topologically unipotent elements in
$\GKq$. We can assume that $f=f_H$ is the characteristic function of such
a set without loss of generality. Let $\chi$ be the character of the
representation of the finite group  $\qGkq_x$ that gave rise to $\pi$
(see Section~\ref{subsection:representations}). 
Since $\GKq_x$ is a maximal compact subgroup of $\GKq$, the Frobenius formula gives the following expression for the character:
\begin{equation}\label{eq:frobenius_char}
\begin{aligned}
\Theta_{\pi}(f_H)
&= \int_{\GKq} \int_{\GKq} f_H(ghg^{-1})\chi_{x,0}(h)\,dh\,dg\\
&= \int_{\GKq} \int_H \chi_{x,0}(g^{-1}hg)\,dh\,dg\\
&= \int_{{\GKq}/\GKq_x} \int_{\GKq_x}\int_H\chi_{x,0}(k^{-1}g^{-1}hgk)\,dh\,dk\,dg\\
&= \mes(\GKq_x)\int_{{\GKq}/\GKq_x}\int_H\chi_{x,0}(g^{-1}hg)\,dh\,dg.
\end{aligned}
\end{equation}
(Recall from Section~\ref{subsection: measures} that the formal degree $d(\pi)$ of $\pi$ is exactly $\chi_{x,0}(1)$.)
The real work is to prove that these integrals converge.
Of course, Harish-Chandra did this long ago when the characteristic of $\Kq$ is $0$.
This is also done in the proof of \cite[Theorem A.14(ii)]{Bushnell-Henniart}. 
We will see the convergence when we do the calculations by hand. 

The next step is to rewrite the outside integral in equation~(\ref{eq:frobenius_char}) as a sum
using the Cartan decomposition $\GKq = \GKq_x A \Ksch{G}_x$, where $A$ is the set of elements of the form  $a_\lambda \ceq \text{diag}(\varpi^\lambda,\varpi^{-\lambda})$ with $\lambda$ a non-negative integer.
This is also done in \cite[Theorem A.14(ii)]{Bushnell-Henniart}, see
equation (\ref{char}) below.
We will also take this formula one step further by collecting the terms
corresponding to each value of the character of ${ \SL}_2(\kq)$ at
unipotent elements. 

Before we can carry out this plan, we need to introduce some more notation.
Recall that there are three unipotent conjugacy classes in ${ \SL}(2,\kq)$: the class $U_0$ of the identity (one element), the class $U_1$ of the element $[\begin{smallmatrix} 1 &1\\0&1\end{smallmatrix}]$, and the  class $U_{\epsilon}$ of the element 
$[\begin{smallmatrix} 1&\epsilon\\0&1\end{smallmatrix}]$, where $\epsilon$ is a non-square in $\kq$.
Suppose $U$ is a unipotent conjugacy class in ${ \SL}(2,\kq)$,  $\lambda$ is a non-negative integer and $x$ is a vertex in the Bruhat-Tits building for $\GKq$. For any regular topologically unipotent element $h$ of $\GKq$,
let $N_{U,\lambda}^{x}(h)$ denote the number of right $\GKq_x$-cosets
$g\GKq_x$ inside the double coset $\GKq_x a_{\lambda}\GKq_x$ that satisfy the following condition:
\[
g^{-1}hg\in \GKq_x
\quad \AND \quad 
\rho_{x,0}(g^{-1}hg)\in U.
\]

\begin{proposition}\label{proposition: N_U}
Let $\pi$ be a supercuspidal depth zero representation of $\GKq$ and let $\Theta_{\pi}$ be its distribution character in the sense of Harish-Chandra. Let $H$ be a compact set of regular topologically unipotent elements in $\GKq$ and let $f_H$ be the characteristic function of $H$. Let $\chi$ be the character of the representation of the finite group $\qGkq_x$ that gave rise to $\pi$ (see Section~\ref{subsection:representations}). Then 
\begin{equation}\label{eqn: N_U}
\Theta_{\pi}(f_H)=\sum_{\lambda}\sum_U\chi(U)\int_H N_{U,\lambda}^{x}(h)\,dh,
\end{equation}
where $U$ runs over unipotent conjugacy classes in ${ \SL}(2,\kq)$ and
$\lambda$ runs over non-negative integers. The sum over $\lambda$ has
only finitely many nonzero terms. 
\end{proposition}

\begin{proof} 
Consider a fixed double coset $\GKq_x a_{\lambda}\GKq_x$.
We observe that if $g_1, g_2 \in \GKq_x$ and   
$g_1 a_{\lambda}\GKq_x=g_2 a_{\lambda}\GKq_x$ then 
$\chi_{x,0}(g_1^{-1}hg_1)= \chi_{x,0}(g_2^{-1}hg_2)$. Hence, using
Cartan decomposition, the double integral from the formula~(\ref{eq:frobenius_char})
can be rewritten as a  double sum, and we obtain the formula  
(cf. \cite[Theorem A.14(ii)]{Bushnell-Henniart}):    
\begin{equation}\label{char}
\begin{aligned}
\Theta_{\pi}(f_H)=\mes(\GKq_x) \sum_{a\in \GKq_x \backslash \GKq /\GKq_x} \sum_{g\in
  {\GKq_x a \GKq_x / \GKq_x }} \int_H \chi_{x,0}(g^{-1}hg)\, dh,
\end{aligned}
\end{equation}
where the summation index in the outside sum runs over the set of
representatives of the double cosets, that is, over $A$, and the
summation index in the inner sum runs over the set of representatives
of the left $\GKq_x$-cosets inside the given double coset. We refer to
\cite[Theorem A.14(ii)]{Bushnell-Henniart}) for the proof that there
are only finitely many  values of $\lambda$ that give nonzero summands. 

Now, consider the contribution of the double coset 
\[
\GKq_x \left[\begin{smallmatrix}\varpi^{\lambda}&0\\0&\varpi^{-\lambda}\end{smallmatrix}\right] \GKq_x
\]
to the sum in the right hand side of equation~(\ref{char}).  The
following two observations make the calculation of the character
fairly simple. First, the function $\chi_{x,0}$ vanishes outside the maximal
compact subgroup $\GKq_x$, while on $\GKq_x$, the value $\chi_{x,0}(g)$
depends only on $\rho_{x,0}(g)$. 
Second, the set $H$ is contained in the set of
topologically unipotent elements, so  for every $h$ in $H$, the element $\rho_{x,0}(g^{-1}hg)$ is a unipotent element in $\qGkq_x$, provided that $g^{-1}hg$ lies in $\GKq_x$. Using these observations it is possible to rewrite the formula for the character (equation (\ref{char})) in the following form:
\begin{equation}\label{ch2}
\Theta_{\pi}(f_H)=\sum_{\lambda}\sum_{U}\chi(U)\int_H N_{U,\lambda}^{x}(h)\,dh.
\end{equation}
\end{proof}

In the next section, we will see that the numbers $N_{U,\lambda}^x(h)$
give rise to geometric objects. Note that these numbers essentially 
depend only on
the group; they depend on the representation only through the choice
of the vertex $x$ that indexes the maximal compact subgroup from which
our representation is induced.

\subsection{Almost definable sets $W_{U,\lambda}^x(h)$}
In this section we begin the process of finding Chow motives related
to the numbers $N_{U,\lambda}^x(h)$ introduced in
Section~\ref{sub:tricks} by expressing these numbers 
as $p$-adic volumes of some sets $W_{U,\lambda}^{x}(h)$; 
we will then use the comparison theorem of Denef and Loeser (see
Section~\ref{subsection:comparison theorem}) to recover these
$p$-adic volumes from the motivic volumes. But first we need to
introduce yet more notation.

Suppose $U$ is a unipotent conjugacy class in $\qGkq_x$,  $\lambda$ is a non-negative integer, and $x$ is a vertex in the Bruhat-Tits building for $\GKq$ as before. For any regular topologically unipotent element $h$ of $\GKq$,
define
\begin{equation}\label{eqn:W_U}
W_{U,\lambda}^{x}(h) \ceq \{y\in \GKq_x\mid
a_{\lambda}^{-1}y^{-1}h ya_{\lambda}\in \GKq_x
 \wedge \rho_{x,0}({(ya_{ \lambda} )^{-1}h (ya{ \lambda})})\in U\}.
\end{equation}

\begin{lemma}\label{lemma: W_U}
With notation as above,
\begin{equation}\label{eqn: W_U}
N_{U,\lambda}^{x}(h)=(q+1)q^{2\lambda-1}\frac{\mes(W_{U,\lambda}^{x}(h))}{\mes(\GKq_x)},
\end{equation}
for all regular topologically unipotent $h\in \GKq$.
\end{lemma}

\begin{proof}
Let $a_{\lambda}$ denote the element
$\text{diag}(\varpi^{\lambda},\varpi^{-\lambda})$, as before. For each $\lambda\in \N\cup \{0\}$, we say that two elements $y_1$ and $y_2$ in $\GKq_x$ are $\lambda$-equivalent, if $y_1a_{\lambda}$ and $y_2a_{\lambda}$ are in the same left $\GKq_x$-coset, that is, if $a_{\lambda}^{-1}y_1^{-1}y_2a_{\lambda}\in \GKq_x$. The $\lambda$-equivalence class of an element $y$ is denoted by $[y]_{\lambda}$.
Each set $W_{U,\lambda}^{x}(h)$ is a finite disjoint union of $\lambda$-equivalence classes (see \cite{G}), and if it contains an element $y$, it contains the whole class $[y]_{\lambda}$. With this notation, the number $N_{U,\lambda}^{x}(h)$ equals the number of  $\lambda$-equivalence classes in the set $W_{U,\lambda}^{x}(h)$.
Now, for $\lambda>0$, the cardinality of $\GKq_x a_{\lambda} \GKq_x/\GKq_x$ is $(q+1)q^{2\lambda-1}$, as can be shown using the affine Bruhat decomposition for $\GKq$ (see \cite{BT}, for example). Since all equivalence classes have equal volumes (see \cite[Lemma 4]{G}), it follows that
\begin{equation}\label{equation: W_U.1}
N_{U,\lambda}^{x}(h)=(q+1)q^{2\lambda-1}\frac{\mes(W_{U,\lambda}^{x}(h))}{\mes(\GKq_x)},
\end{equation}
as claimed.
\end{proof}

The sets $W_{U,\lambda}^x(h)$ are almost 
definable, but they depend on the parameter $h$
that cannot be specified within the language. Still, it is possible to  
use motivic integration to calculate their volumes, by using the
version of motivic integration that allows parameters.

\subsection{Motivic volumes of the sets $W_{U,\lambda}^{x}(h)$}\label{W_Clambda}

In this section, we write down the formulas defining the sets 
$W_{U,\lambda}^{x}(h)$ of the previous section, and then calculate
their motivic volumes. This is done case-by-case in $z$. We start with the
most interesting case ${z=s_1}$. We include the details only for
$x=(0)$. 
Before we start the calculation, two remarks are due.

First, recall that we had to fix  a non-square unit $\eps$ and 
a uniformizer $\varpi$, as discussed  in Section~\ref{section:outline}.
In this section, we introduce a variable $\delta$ that would be
allowed to range over non-square units. 
As discussed in Section~\ref{subsub:language}, a $p$-adic field
together with the choice of the uniformizer is a structure for the
language of Denef-Pas (\ie, such a choice provides an interpretation
of all the formulas). If we choose a non-square unit $\eps$ in the
given field, carry out all the constructions that appeared so
far, and then plug in the value of $\eps$ for $\delta$ in our formulas, 
we will get the corresponding sets $W_{U,\lambda}^{x}(h)$. We will
see, however, that their motivic volumes (and therefore, also $p$-adic
volumes) are independent of the choice of $\eps$. 

Second, the comparison theorem relates the motivic volumes to
$p$-adic volumes  for all but finitely many primes $p$. There are two
sources of ``bad'' primes in this statement: singularities of the
varieties that come up, and quantifier elimination. Since we are doing
all the calculations by hand, and our cases are very simple, we will
see that the only prime that needs to be excluded is $2$.

\subsubsection{The logical formulas for $W_{U,\lambda}^{(0)}(h)$ for $h\in{\Gamma}_{{s_1},n}$}\label{subsub: formulas for W_s1}

Recall that all the elements of the set 
$\Gamma_{s_1,n}$ have the form
\begin{equation}\label{theelement}
h=(1+\ast)
\left[
\begin{matrix}
1+\ast& u\varpi^n\\
\eps u\varpi^n &1+\ast
\end{matrix}
\right],
\end{equation}
where $u\in\Uq$,  
and '$\ast$' stands for 
 elements of order at least $2n$. 
We write $y=[\begin{smallmatrix}a&b\\c&d\end{smallmatrix}]\in \GKq_{(0)}$
and use the entries $a,b,c,d$ as free variables in our formulas. 
We begin by fixing $U$ and $h$ as above and finding explicit
conditions on $a,b,c,d$ and $\lambda$ which ensure that $y\in
W_{U,\lambda}^{(0)}(h)$. 
 A direct computation gives:
\begin{equation}
y^{-1}h y=(1+\ast)\left[
\begin{matrix}
1+\text{h.o.t} & u(d^2-b^2\eps)\varpi^n +\text{h.o.t}\\
u(a^2\eps-c^2)\varpi^n+\text{h.o.t}&1+\text{h.o.t} 
\end{matrix}
\right],
\end{equation}
where \text{h.o.t} means ``higher order terms'', which refers to the terms with valuation greater than that of the leading term. 

As we see from equation~\ref{eqn:W_U} we must now conjugate this element by $a_\lambda$.
Conjugation by $a_{\lambda}$ induces multiplication of the entry in the upper right-hand corner by $\varpi^{-2\lambda}$, and multiplication of the entry in the lower left-hand corner by $\varpi^{2\lambda}$. It follows that $a_{\lambda}^{-1}y^{-1}h y a_{\lambda}\in \GKq_{(0)}$ if and only if
\begin{equation}\label{conditions}
-2\lambda+\ord(u(d^2-b^2\eps)\varpi^n))\ge 0 \quad \AND \quad 
2\lambda+\ord(u(a^2\eps-c^2)\varpi^n))\ge 0.
\end{equation}
Consequently, we have the following list of cases:

\noindent{\bf Case 1}. $2\lambda<n$.
In this case,  for any $y\in \GKq_{(0)}$ 
the element $a_{\lambda}^{-1}y^{-1}h y a_{\lambda}$ lies in  $\GKq_{(0)}$
and projects to the identity under the reduction mod $\varpi$.

\noindent{\bf Case 2.} $2\lambda>n$.
We observe that $\ord(d^2-b^2\eps)=0$: indeed, either both $a$ and $d$ or $b$ and $c$ have to be units because $ad-bc=1$, and $\rho_{(0),0}(d^2-b^2\eps)\neq 0$ since $\eps$ is a non-square.  Similarly, $\ord(a^2\eps-c^2)=0$. This implies that if $2\lambda>n$, the set of $y$s satisfying the  conditions is empty.

\noindent{\bf Case 3.} $2\lambda =n$. (This case corresponds to the interesting situations.) In this case the image of $a_{\lambda}^{-1}y^{-1}h ya_{\lambda}$ under $\rho_{(0),0}$ is an element of the form $[\begin{smallmatrix}1&\beta\\ 0&1\end{smallmatrix}]$, where $\beta=\ac{d^2-b^2\eps}\ac{u}\in\kq$. Then for $h\in \Gamma_{{s_1},n,+}$, the reduction 
$\rho_{(0),0}(a_{\lambda}^{-1}y^{-1}h y a_{\lambda})$ falls into $U_1$ if $\ac{d^2-b^2\eps}$ is a square and into $U_{\eps}$ otherwise.

Notice that when $\lambda$ is fixed and $h$ is fixed,  the condition that $y\in W_{U,\lambda}^{(0)}(h)$ does not depend on $u$, as long as $h$ is confined to one of the sets $\Gamma_{{s_1},n,+}$ and 
$\Gamma_{{s_1},n,-}$.

We now reformulate what we have found concerning the character in
terms of formulas in Pas's language. We
start by observing that the maximal compact subgroup $\GKq_{(0)}$ is defined by the formula 
$$\text{`}
ad-bc=1\quad\wedge \ord(a)\ge 0\quad\wedge\ord(b)\ge 0\quad\wedge
\ord(c)\ge 0\quad\wedge \ord(d)\ge 0.
\text{'}$$

For the cases $2\lambda<n$ or $2\lambda>n$, no other calculations are
needed.
Indeed, if $2\lambda<n$, then we are within Case 1, and 
therefore both sets
$W_{U_{1,\eps},\lambda}^{(0)}(h)$
are empty, and $W_{U_0,\lambda}^{(0)}(h)=\GKq_{(0)}$. 
If $2\lambda>n$, then we are within Case 2, and the set of $y$s
satisfying the conditions is also empty. 

Let us now consider the case $n$ even, and  $2\lambda =n$.
Let $\psi(b,d,\delta)$ be the formula (in Pas's language)
$$
\psi(b,d,\delta)=
\text{'}\exists x (d^2-b^2\delta=x^2)\text{'}.
$$

Since we will often pass back and forth between the valued field and
the finite field, let us introduce an abbreviation for the 
``reduction $\mod \varpi$'' map: let
$$
\bar x=\begin{cases}
\ac{x}, & \ord(x)=0\\
0,& \ord(x)>0.
\end{cases}
$$

If we  know that $\delta$ is a non-square unit, then by Hensel's
Lemma, to check whether the
triple $(b,d,\delta)$ with $(\bar b,\bar d)\neq (0,0)$ 
satisfies the formula $\psi$, it is enough
to know whether its reduction $\mod (\varpi)$ satisfies the same
formula
(where now the variables are interpreted as residue-field variables).
Hensel's Lemma is applicable because $\bar d^2-\bar b^2\bar\delta$
cannot be $0$ when we assume that $\bar\delta$ is a non-square and $\bar
b,\bar d$ are not simultaneously zero.  

Let us consider the family of formulas depending on a parameter 
$\eta$ that ranges over the set of non-squares in $\kq$ (note that this
is a definable set):
\begin{equation}
\phi_{\eta}(a,b,c,d)\quad=\text{'}ad-bc=1 \quad\wedge\exists \xi
(\bar b^2-\bar d^2\eta=\xi^2)\text{'}.
\end{equation}

For every value of $\eta\in\kq^{\ast}\setminus {\kq^{\ast}}^2$, 
if we let all the variables range over $\Rq$, 
the formula $\phi_{\eta}(a,b,c,d)$
defines the set $W_{U_1,\lambda}^{(0)}(h)$ for any $h\in
\Gamma_{{s_1},n,+}$ if the unit $\eps$ that was fixed in  order to define
the sets $\Gamma_{{s_1},n,\pm}$ has the property $\bar\eps=\eta$.  
The same formula also defines the set 
$W_{U_{\eps},\lambda}^{(0)}(h)$ with  $h\in\Gamma_{{s_1},n,-}$.
In the next subsection we find the motivic volumes of these sets.

\subsubsection{The motivic volumes of $W_{U,\lambda}^{(0)}(h)$,  
$h\in\Gamma_{{s_1},n,\pm}$}\label{subsection: anything}

The calculation of these volumes is very simple. Recall that our Haar 
 measure on $\GKq$ is normalized in such a way that the fibres of the
 projection $\GKq\to G(\kq)$ have volume $1$. The formula
 $\phi_{\eta}$ imposes a condition on the variables $b,d$ that 
apparently depends only on $\bar b, \bar d$. This means that either a
 whole fibre over a point $(\bar a, \bar b, \bar c, \bar d)\in
 { \SL}(2,\kq)$ satisfies this condition, or the whole fibre does not
 satisfy it. Hence, to calculate the volume of the set 
$W_{U_1,\lambda}^{(0)}(h)$, all we need to do is
 count the number of points in 
$[\begin{smallmatrix}\bar a &\bar b\\ \bar c & \bar
 d\end{smallmatrix}]\in { \SL}(2,\kq)$ that satisfy the condition
$\exists \xi: \bar b^2-\bar d^2\eta=\xi^2$, where $\eta$ is a parameter that is
 a quadratic residue in $\kq$.  
{ This calculation is carried out
 carefully in \cite{JG: overview}. Here we state the result as a lemma}.

In fact, this is the only place in the present paper where
interesting geometric objects come up from $p$-adic volumes. 
{ Indeed, as we see from the argument above,
the $p$-adic volumes of the sets $W_{U,\lambda}^{(0)}(h)$
with the ``border-line'' value of $\lambda=n/2$ (for $n$ even) are
connected with the number of points of a conic over the finite field}.  

\begin{lemma}\label{lemma:unram.elements}
Suppose $h\in\Gamma_{{s_1},n}$.
If $n$ is even and $\lambda =n/2$ 
then the motivic volumes of the sets $W_{U_1,\lambda}^{(0)}(h)$ and $W_{U_{\eps},\lambda}^{(0)}(h)$ both equal $\frac{1}{2}\lef(\lef-1)(\lef+1)$, which is half of the motivic volume of the maximal compact subgroup $\GKq_{(0)}$; otherwise, the motivic volumes of the sets $W_{U_1,\lambda}^{(0)}(h)$ and $W_{U_{\eps},\lambda}^{(0)}(h)$ both equal $0$.
\end{lemma}
\begin{proof} 
{ The statement follows from the calculation of the motivic 
volume of the formula 
\begin{equation}
\begin{aligned}
& \Phi(a,b,c,d,\eta)\\
&=\text{`}
ad-bc=1\quad\wedge
\exists \xi \neq 0(\bar d^2-\bar b^2\eta=\xi^2)
\quad\wedge
\nexists \beta (\eta=\beta^2) 
\text{'}.
\end{aligned}
\end{equation}
that is carried out in \cite{JG: overview}. Note that the calculation 
ultimately boils down to computing the class of the ``hyperbola'' 
$x^2-y^2=1$, which is $\lef-1$. 
This is the reason that we get an answer that is polynomial in $\lef$. }
\end{proof}

\subsubsection{The motivic volumes of $W_{U,\lambda}^{(0)}(h)$ for 
$h\in\Gamma_{{s_2},n}$}

In the case ${z=s_2}$ the answer is essentially the same as in the case
$z={s_1}$, but the calculation is slightly more complicated. Here we sketch the
calculation of the motivic volume of the sets
$W_{U,\lambda}^{(0)}(h)$ for $h\in \Gamma_{{s_2},n,\pm}$, 
indicating the  differences with the case $z={s_1}$. 

Exactly as before, we consider the element 
$a_{\lambda}^{-1}y^{-1}\gamma_{{s_2},n}ya_{\lambda}$, and 
use the entries of the matrix  
$y=\left[\begin{smallmatrix}a&b\\c&d\end{smallmatrix}\right]\in
\GKq_{(0)}$ as free variables in our logical formulas.
 
The conditions~(\ref{conditions}) are now replaced with
\begin{equation}\label{conditions.s2}
\begin{aligned}
&-2\lambda+\ord(d^2u\varpi^{n-1}-b^2u\eps\varpi^{n+1}))\ge 0 \quad\\ 
&2\lambda+\ord(a^2u\eps\varpi{n+1}-c^2u\varpi^{n-1})\ge 0.
\end{aligned}
\end{equation}
The second condition is satisfied automatically if $\lambda\ge 0$, so
we only need to focus on the first one.
Similarly to the case $z={s_1}$, this condition implies that when 
$\lambda>\frac{n+1}2$, the element 
$a_{\lambda}^{-1}y^{-1}hya_{\lambda}$ is outside $\GKq_{(0)}$;
if $\lambda<\frac{n-1}2$, this element is in $\GKq_{(0)}$ and projects to $U_0$.

Suppose for now that $n$ is odd.
Unlike the case $z={s_1}$, now there are two interesting cases: $\lambda=\frac{n-1}2$ and
$\lambda=\frac{n+1}2$.
If  $\lambda=\frac{n-1}2$, and $\ord(d)= 0$, then  
$\rho_{(0),0}(a_{\lambda}^{-1}y^{-1}hya_{\lambda})\in U_1$
if ${\sgn}(u)=1$ and
$\rho_{(0),0}(a_{\lambda}^{-1}y^{-1}hya_{\lambda})\in U_{\eps}$ 
if ${\sgn}(u)=1$. If $\ord(d)> 0$, then  
$\rho_{(0),0}(a_{\lambda}^{-1}y^{-1}hya_{\lambda})\in U_0$.
If $\lambda=\frac{n+1}2$, then 
$a_{\lambda}^{-1}y^{-1}hya_{\lambda}$ is in $\GKq_{(0)}$
only if $\ord(d)>0$. In this case, the projection of the element 
$a_{\lambda}^{-1}y^{-1}hya_{\lambda}$ depends on the sign 
of $\ac{d^2-b^2\eps}$ and on ${\sgn}(u)$, which is similar to the case  
$z={s_1}$. The only difference is that here there is an additional
condition $\ord(b)=0$ (if $\ord(d)>0$ then the determinant condition
forces $\ord(b)=0$).
{ One finds that the motivic volume of the subset of $G$ defined by the formula 
'$\exists \xi\neq 0 (\ac{d^2-b^2\eps}=\xi^2)
\wedge (\ord(b)=0) \wedge (\ord(d)>0)$' equals $\frac12(\lef-1)^2$.
The motivic volume of its complement, \ie, the set defined by the formula 
'$\nexists \xi\neq 0 (\ac{d^2-b^2\eps}=\xi^2)
\wedge (\ord(b)=0) \wedge (\ord(d)>0)$' equals 
$\lef(\lef^2-1)-\frac12(\lef-1)^2=\frac12\lef^2-\frac12$.
The last two sentences are of course an abbreviation. In truth, we have to
replace $\eps$ with a variable $\delta$ and do everything exactly the
same way as it was done in the previous case.} 

For simplicity, suppose that ${\sgn}(u)=1$.
Putting all these calculations together, we get:
\begin{equation}\label{eq:W_s2}
\begin{aligned}
&\mu(W_{\frac{n+1}2, U_{\eps}})=\frac12\lef^2-\frac12;\\
&\mu(W_{\frac{n-1}2, U_{\eps}})=0;\\
&\mu(W_{\frac{n-1}2, U_{1}})=\lef^2(\lef-1);\\
&\mu(W_{\frac{n+1}2, U_{1}})=\frac12(\lef-1)^2.
\end{aligned}
\end{equation}
The case $n$ even, and the calculation of $\mu(W_{\lambda,U_0})^{(0)}(h)$ 
are similar to the case {$z=s_1$}.

\subsubsection{The ramified cases}

\begin{lemma}\label{lemma:ram.elements}
Let $x$ be a vertex in the Bruhat-Tits building for $\GKq$. Suppose $h\in
\Gamma_{z,n,\pm}$,  with { $z\in\{t_0,t_1,t_2,t_3\}$}. 
Then the sets $W_{U,\lambda}^{x}(h)$ are definable for all non-negative
integers $\lambda$,  their motivic volumes are independent of $h$
and the choice of $\eps$, and can be explicitly computed.
\end{lemma}

\begin{proof}
We prove this lemma only  for $x=(0)$. The proof in the case $x=(1)$
is very similar; the  results of these calculations become part of 
the expressions for $M_{z,n,\pm}^{(1),0}$, $M_{z,n,\pm}^{(1),1}$ and $M_{z,n,\pm}^{(1),\eps}$ appearing in 
Tables~\ref{table:m-char.1.s} and \ref{table:m-char.1.t}. 
So everywhere in this proof $x=(0)$, and we drop
the superscript $x$ from the notation $W_{U,\lambda}^x(h)$.  

The argument is very similar to the unramified case; the only
difference is that the actual calculation of the motivic volumes 
of the corresponding sets $W_{U,\lambda}(h)$ is simpler. Here we carry out the proof for the
case {$z=t_2$}. The other three ramified cases are almost identical to
it. First, as in the previous subsection, we write the elements $h$ of
the set $\Gamma_{{t_2},n}$ explicitly as:
\begin{equation}\label{ram.element}
h=(1+\ast)
\left[
\begin{matrix}
1+\ast & u\varpi^n\\
\eps u\varpi^{n+1} &1+\ast
\end{matrix}
\right], \text{where} \quad u\in\Uq,
\end{equation}
and $\ast$ denotes the terms of order at least $2n$, as before. 

As in the proof of the previous lemma, we let
$y=\left[\begin{smallmatrix}a&b\\c&d\end{smallmatrix}\right]$ be a variable
running over $\GKq_{(0)}$ (so that the symbols for its entries $a,b,c,d$ will
become the formal variables in the Pas's language formulas defining the sets 
$W_{U,\lambda}^{(0)}(h)$). 

{As before, we compute $y^{-1}hy$,  
which leads to the following conditions on $y, \lambda$ for the element 
$a_{\lambda}^{-1}y^{-1}hya_{\lambda}$ to be in $\GKq_{(0)}$:}
\begin{equation}\label{ram.conditions}
\begin{aligned}
-2\lambda+\ord(ud^2\varpi^n-b^2\eps u\varpi^{n+1})&\ge 0 \text {  and}\\ 
2\lambda+\ord(ua^2\eps\varpi^{n+1}-c^2u\varpi^n)&\ge 0.
\end{aligned}
\end{equation}
As before, it is convenient to consider the cases $n$ even and $n$ odd
separately. 

Suppose $n$ is even. Looking at the left-hand side of the inequalities
(\ref{ram.conditions}), we see that  
if $\lambda<n/2$, then the set $W_{U_0,\lambda}(h)$ coincides with the 
whole of $\GKq_{(0)}$, and the sets 
$W_{U_1,\lambda}(h)$, $W_{U_{\eps},\lambda}(h)$ are empty.

If $\lambda>n/2$, then, since $n$ is even and $\lambda$ is an integer, 
$2\lambda$ is at least $n+2$, which forces the element 
$a_{\lambda}^{-1}y^{-1}hya_{\lambda}$ outside $\GKq_{(0)}$, and
all three sets are empty.

Finally, in the case $\lambda=n/2$, the outcome depends on the entry
$d$:
if $\ord(d)= 0$, then the term $d^2u\varpi^n$ in the expression 
$d^2u\varpi^{n}-b^2\eps u\varpi^{n+1}$ dominates, and therefore 
$\rho_{(0),0}(a_{\lambda}^{-1}y^{-1}hya_{\lambda})\in U_1$ if
${\sgn}(u)=1$, and 
$\rho_{(0),0}(a_{\lambda}^{-1}y^{-1}hya_{\lambda})\in U_{\eps}$ if
${\sgn}(u)=-1$.
If $\ord(d)>0$, then 
$\rho_{(0),0}(a_{\lambda}^{-1}y^{-1}hya_{\lambda})\in U_{0}$. 

Hence, for $h\in\Gamma_{{t_2},n,\pm}$, the set 
$W_{U_0,n/2}(h)$ is defined by the formula 
'$\ord(d)>0$' in conjunction with the formulas defining $\GKq_{(0)}$.
Note that the volume of this set is the same as that of $\GKq_{(01)}$,
\ie, equals $\lef(\lef-1)$.
The sets $W_{U_1,n/2}(h)$, $W_{U_{\eps},n/2}(h)$ are defined by 
the formula '$\ord(d)=0$' in conjunction with the formulas defining
$\GKq_{(0)}$
for $h\in \Gamma_{t_2,n,+}$ and $h\in\Gamma_{t_2,n,-}$, respectively,
and are respectively empty for $h\in \Gamma_{t_2,n,-}$ 
and $h\in\Gamma_{t_2,n,+}$.

The case $n$ odd is very similar.
If  $\lambda=\frac{n+1}2$ (the most interesting case),
and $\ord(d)= 0$, the second one of the conditions
(\ref{ram.conditions}) is not satisfied, so all three sets 
$W_{U,\frac{n+1}2}(h)$ are
empty. If $\ord(d)>0$, then automatically $\ord(b)=0$, the leading
term in the expression $d^2u\varpi^{n}-b^2\eps u\varpi^{n+1}$ is  
$-b^2u\eps\varpi^{n+1}$, which has sign opposite to ${\sgn}(-1){\sgn}(u)$.

We get:  if $\lambda \leq (n-1)/2$ then 
 \[
 W_{U_0,\lambda}=\GKq_{(0)},\quad \text{ and} \quad
 W_{U_1,\lambda}(h)=W_{U_\eps,\lambda}(h)=\emptyset, \text{ for any } 
h\in \Gamma_{t_2,n};
 \]
 if $\lambda > (n+1)/2$ then all three sets are empty; if $\lambda =
 (n+1)/2$, we have
\begin{equation}
\begin{aligned}
W_{U_1,\frac{n+1}2}(h)=\GKq_{(0)}\cap \{\ord(d)>0\}, &\quad
W_{U_{\eps},\frac{n+1}2}(h)=\emptyset, 
&\quad h\in\Gamma_{t_2,n,{\sgn}(-1)};\\
W_{U_{\eps},\frac{n+1}2}(h)=\GKq_{(0)}\cap\{\ord(d)>0\}, &\quad
W_{U_1,\frac{n+1}2}(h)=\emptyset, &\quad
h\in\Gamma_{t_2,n,-{\sgn}(-1)}.\\
\end{aligned}
\end{equation}
It follows that $\mu(W_{U_1,\frac{n+1}2}(h))$ equals $0$ or equals 
$\mu(G_{(01)})=\lef(\lef-1)$ depending on whether $h\in \Gamma_{t_2,n,{\sgn}(-1)}$
or $h\in \Gamma_{t_2,n,-{\sgn}(-1)}$; the same is true for 
$\mu(W_{U_{\eps},\frac{n+1}2}(h))$.

As we see from this proof, a different choice of $\eps$ could not have
affected the motivic volumes of these sets; also, clearly no ``bad''
primes were acquired.
\end{proof}

\subsubsection{The case $x=(1)$}

In all the proofs in the present section we have been assuming that $x=(0)$. The only major difference
of the case $x=(1)$ is that the element 
$y=\left[\begin{smallmatrix}a&b\\c&d\end{smallmatrix}\right]$
now belongs to $\GKq_{(1)}$, not $\GKq_{(0)}$. The requirement is that
the element $a_{\lambda}^{-1}y^{-1}hya_{\lambda}$ belongs to
$G_{(1)}$, and and the reduction map
applied to this element  is $\rho_{(1),0}$ instead of $\rho_{(0),0}$.
Therefore all the cases look slightly different, but no new
varieties appear in the calculations of the motivic volumes of the new
sets.

\subsection{Motives corresponding to the distribution characters}

\begin{proposition}\label{prop:m-ch}
The Harish-Chandra character of each depth zero supercuspidal
representation is constant on each set ${}^G\Gamma_{z,n,\pm}$. 
Moreover, there exist virtual motives $M_{z,n,\nu}^{x,0}$, $M_{z,n,\nu}^{x,1}$ and
$M_{z,n,\nu}^{x,\eps}$ (where $z$ is any cocycle defined in 
Section~\ref{section:outline}, $n$ is a positive integer for 
$z\in\{s_0,s_1,s_2\}$ and a non-negative integer for
$z\in\{t_0,t_1,t_2,t_3\}$,
the sign  $\nu$ is $\pm$ and $x$ is a vertex $(0)$ or $(1)$) such that
\begin{equation}\label{motivic character}
\begin{aligned}
&\frac{1}{m({}^G\Gamma_{z,n,\nu})} \Theta_{\pi}(f_{z,n,\nu})\\
& =
\chi\vert_{U_0}\frob M_{z,n,\nu}^{x,0}
+\chi\vert_{U_1}\frob M_{z,n,\nu}^{x,1} +\chi\vert_{U_{\eps}}\frob M_{z,n,\nu}^{x,\eps}.
\end{aligned}
\end{equation} 
The virtual motives $M_{z,n,\nu}^{x,0}$, $M_{z,n,\nu}^{x,1}$ and 
$M_{z,n,\nu}^{x,\eps}$ are explicitly given in
Tables~\ref{table:m-char.0.s} and \ref{table:m-char.0.t} in the case $x=(0)$, and in
Tables~\ref{table:m-char.1.s} and \ref{table:m-char.1.t} in the case $x=(1)$.
\end{proposition}

\begin{proof}
We will show the details of the  proof of this proposition in the case
$x=(0)$ only. The results of the similar calculations in the case
$x=(1)$ are summarized in the tables.

Let us apply Proposition~\ref{proposition: N_U} to the test functions
$f_{z,n,\nu}$
(so that the set $H$ in that Proposition is
${}^G\Gamma_{z,n,\nu}$).
Note that in $\GKq$, for all $z$ but $s_0$, 
${\mathfrak h}_z$ is elliptic, \ie, it is  a compact Cartan
subalgebra. We consider the elliptic cases first. 
As we will explicitly see below, the sum over $\lambda$ that appears
in the expression (\ref{eqn: N_U}) has only finitely many terms
in these cases, so we can permute the two sums, and obtain
\begin{equation}
\Theta_{\pi}(f_{z,n,\nu})=
\sum_U\chi(U)\sum_{\lambda}\int_{{}^G\Gamma_{z,n,\nu}} N_{U,\lambda}^{x}(h)\,dh,
\end{equation}  
where
the index $U$ runs over $U_0$, $U_1$, and $U_{\eps}$.
That is, the character is already expressed as a linear combination of
the values $\chi\vert_{U_0}$, $\chi\vert_{U_1}$, and $\chi\vert_{U_{\eps}}$.
All we need to do, is ``evaluate'' the coefficients
$\sum_{\lambda}\int_{{}^G\Gamma_{z,n,\nu}}N_{U,\lambda}^x(h)$.  
We recall Lemma~\ref{lemma: W_U} which relates the numbers 
$N_{U,\lambda}^x(h)$ to the volumes of the sets
$W_{U,\lambda}^x(h)$. Then  we evaluate their {\it motivic} volumes 
for $h$ in each of the sets ${}^G\Gamma_{z,n,\nu}$ and
sum them over $\lambda$ with coefficients that come from 
Lemma~\ref{lemma: W_U}.

We start with $z=s_1$.
Let $h$ be an element of $\Gamma_{s_1,n,\pm}$, and let $\lambda$ be a
non-negative integer. 
By the comparison theorem, the  equality (\ref{equation: W_U.1}) of 
Lemma~\ref{lemma: W_U} can be written in a ``motivic'' form:
for all primes $p\neq 2$ (recall that $q$ is a power of $p$), we have 
$N_{U,\lambda}^x(h)=\frob_q M_{U,\lambda}^x$, where
$$
M_{U,\lambda}^x=[(q+1)q^{2\lambda-1}]\frac{\mu(W_{U,\lambda}^x)(h)}{\mu(\GKq_x)}=
(\lef+1)\lef^{2\lambda-1}\frac{\mu(W_{U,\lambda}^x)(h)}{\mu(\GKq_x)}.
$$
For now, let $x=(0)$.
Note that, {\it a priori}, the right hand side depends on the element
$h$. However, { by} Lemma~\ref{lemma:unram.elements} 
for all $h\in\Gamma_{s_1,n,\pm}$,  and for every $U$, 
the motivic volume of $W_{U,\lambda}^x(h)$ does not depend on $h$ and equals:
\begin{equation*}
\begin{aligned}
\frac12\lef(\lef^2-1),& \quad \text{if\ }\lambda=n/2,\  n\text{\
  even,\ } U=U_1
\text{\ or \ } U_{\eps},\\
\lef(\lef^2-1),&\quad \text{if\ }\lambda<n/2, \ U=U_0\\
0, & \quad \text{otherwise}.
\end{aligned}
\end{equation*}

Let us define the virtual Chow motive $M_{s_1,n,\pm}^{(0),0}$ that
corresponds to the conjugacy class $U_0$  by the formulas
\begin{equation}
\begin{aligned}
&M_{s_1,n,\pm}^{(0),0}&=&1+(\lef+1)\sum_{\lambda=1}^{n/2-1}\lef^{2{\lambda}-1}
=\frac{\lef^{n-1}-1}{\lef-1} &\quad\text{if } n \text{ is even}, \\
&M_{s_1,n,\pm}^{(0),0}&=&1+(\lef+1)\sum_{\lambda=1}^{(n-1)/2}\lef^{2{\lambda}-1}=
\frac{\lef^n-1}{\lef-1} &\quad\text{if } n \text{ is odd.}
\end{aligned}
\end{equation} 
Also let
\begin{equation} 
M_{s_1,n,\pm}^{(0),1} = M_{s_1,n,\pm}^{(0),\eps}=
\begin{cases}
	\frac12(\lef+1)\lef^{n-1}, &\text{ if $n$ is even}\\
	0, & \text{ if $n$ is odd}
\end{cases}
\end{equation}

Then, combining the equations above,  we obtain, for $p\neq 2$, 
\begin{equation}\label{mch}
\begin{aligned}
\hskip-20pt 
& \frac{\Theta_{\pi}(f_{s_1,n,\nu})}{m({}^G\Gamma_{s_1,n,\nu})} 
=\Theta_{\pi}(h)\\
& =  \chi\vert_{U_0}\frob\,M_{s_1,n,\nu}^{(0),0}+\chi\vert_{U_1}\frob
M_{s_1,n,\nu}^{(0),1}+\chi\vert_{U_\eps}\frob M_{s_1,n,\nu}^{(0),\eps}, 
\end{aligned}
\end{equation}
where we write $\Theta_{\pi}$ both for the distribution character 
and for the  locally integrable function on the regular set that
represents it.  
Note that there is no difference in the formulas for
$\Theta_{\pi}(f_{s_1,n,+})$ and $\Theta_{\pi}(f_{s_1,n,-})$, 
because $M_{s_1,n,\pm}^{(0),1}$ in any case 
coincides with $M_{s_1,n,\pm}^{(0),\eps}$.
The proposition in the case $z=s_1$ is proved.

In the case $z=s_2$, the calculation is very similar to the case
$z=s_1$, except that when $n$ is odd, the expressions for  
$M_{s_2,n,\nu}^{(0),1}$  and $M_{s_2,n,\nu}^{(0),{\eps}}$ 
contain one or two terms depending on $\nu$, yet the final answer
is the same in the both cases. 
Using equation~(\ref{eq:W_s2}) we get:
\begin{equation*}
\begin{aligned}
M_{s_2,n,+}^{(0),1}=M_{s_2,n,-}^{(0),{\eps}}
&= \frac{1}{2}\frac{(\lef-1)^2}{\lef(\lef^2-1)}\lef^{2\frac{n+1}2-1}(\lef+1)+\frac{\lef}{\lef+1}\lef^{2\frac{n-1}2-1}(\lef+1)\\
&=\frac12\lef^{n-1}(\lef+1);\\
M_{s_2,n,-}^{(0),1}=M_{s_2,n,+}^{(0),{\eps}}
&=\frac{1}{2}\frac{(\lef^2-1)}{\lef(\lef^2-1)}\lef^{2\frac{n+1}2-1}(\lef+1)\\
&=\frac12\lef^{n-1}(\lef+1).
\end{aligned}
\end{equation*}
The values of $M_{s_2,n,\nu}^{(0),0}$ are computed similarly to the
case $z=s_1$.  

Let us now prove the proposition for the ramified elements. All
ramified cases are very similar to each other. 
We show the details for the  case $z=t_2$. The argument is
exactly the same as in the case $z=s_1$, but we have to use
Lemma~\ref{lemma:ram.elements} instead of Lemma~\ref{lemma:unram.elements}.
In the case $n$ even, we get: 
\begin{equation}
\begin{aligned}
M_{t_2,n,\pm}^{(0),0}&=1+\sum_{\lambda=1}^{n/2-1}\lef^{2\lambda-1}
+\frac{\mu(\GKq_{(01)})}{\mu(\GKq_{(0)})}\lef^{2\frac n2-1}(\lef+1)\\
&=\frac{\lef^{n-1}-1}{\lef-1}
+\frac1{\lef+1}\lef^{2\frac n2-1}(\lef+1)=\frac{\lef^n-1}{\lef-1};\\
M_{t_2,n,+}^{(0),1}&=M_{t_2,n,-}^{(0),\eps}
=\frac{\mu(\GKq_{(0)})-\mu(\GKq_{(01)})}{\mu(\GKq_{(0)})}\lef^{2\frac n2-1}(\lef+1)\\
&=\frac{\lef}{\lef+1}\lef^{2\frac n2-1}(\lef+1)=\lef^n; \\
M_{t_2,n,-}^{(0),1}&=M_{t_2,n,+}^{(0),\eps}=0.
\end{aligned}
\end{equation}
If $n$ is odd, a similar calculation yields:
\begin{equation}
\begin{aligned}
M_{t_2,n,\pm}^{(0),0}& = 1+(\lef+1)\sum_{\lambda=1}^{\frac{n-1}2}\lef^{2\lambda-1}(\lef+1) =\frac{\lef^n-1}{\lef-1}.\\
M_{t_2,n,{\sgn}(-1)}^{(0),1}& = M_{t_2,n,-{\sgn}(-1)}^{(0),\eps} =
\frac{\lef}{\lef+1}\lef^{2\frac {n+1}2-1}(\lef+1) =\lef^n. \\
M_{t_2,n,-{\sgn}(-1)}^{(0),1}& = M_{t_2,n,{\sgn}(-1)}^{(0),\eps} =0.
\end{aligned}
\end{equation}
The Proposition for $z=t_2$ follows.

Finally, let us address the case $s=s_0$.
The elements of the set $\Gamma_{s_0,n}$ have the form
$h=\left[\begin{smallmatrix}1+u\varpi^n
    &0\\0&(1+u\varpi^n)^{-1}\end{smallmatrix}\right]$, where $u$ is a unit.
Following the pattern of the previous section, we take 
$y=\left[\begin{smallmatrix}a&b\\c&d\end{smallmatrix}\right]\in
\GKq_{(0)}$, and write down the
conditions ensuring that the element  
$a_{\lambda}^{-1}y^{-1}hya_{\lambda}$ belongs to $G_{(0)}$ and
projects to a given conjugacy class under the map $\rho_{(0),0}$. As
  before, we see that this depends on the entry in the right-hand
  corner of the matrix 
$a_{\lambda}^{-1}y^{-1}hya_{\lambda}$, which equals 
$\varpi^{-2\lambda}bd((1+u\varpi^n)-(1+u\varpi^n)^{-1})$.
We have (when $p\neq 2$):
\begin{equation}
\begin{aligned}
&\ord(\varpi^{-2\lambda}bd((1+u\varpi^n)-(1+u\varpi^n)^{-1}))&=&
-2\lambda+n+\ord(bd),\\
&\ac{\varpi^{-2\lambda}bd((1+u\varpi^n)-(1+u\varpi^n)^{-1})}&=& 2\ac{bd}\ac{u}.
\end{aligned}
\end{equation}
Here the situation is quite different from the elliptic cases, because
the valuation of $bd$ can be arbitrarily large, and therefore there
are {\it infinitely many} values of $\lambda$ such that 
$a_{\lambda}^{-1}y^{-1}hya_{\lambda}$  belongs to $G_{(0)}$. We know,
of course, that the sum in the equation (\ref{eqn: N_U}) 
has to be finite anyway. 
Here we will see explicitly that it happens because both for $\chi=Q_T$
and $\chi=Q_G$ (and therefore, for any linear combination of these two
functions also), for large values of $\lambda$ the sum of three terms 
$\sum_U\chi\vert_U N_{U,\lambda}^{(0)}(h)$ vanishes.
Indeed, if $\chi=Q_G$, it is easy to see that this sum is always zero,
because $Q_G$ vanishes on the class $U_0$ and takes opposite values
on $U_1$ and $U_{\eps}$.
Suppose $\chi=Q_T$. Then
for each positive integer $k$, 
we will need the volumes of the subsets of $G_{(0)}$ defined by the
formulas $\{\ord(b)\ge k\}$ and $\{\ord(d)\ge k\}$.
Note that these sets are disjoint
and have equal volumes. From the point of view of motivic
integration, it is easy to see that these  sets are stable at level
$k$ in the language of \cite{DL.ar}, and their motivic volumes equal
$(\lef-1)\lef^{-(k-2)}$. The only varieties appearing in this
calculation are affine spaces, so we acquire no bad primes. 
Now it is easy to see that when $2\lambda>n$, the value 
$\chi\vert_{U_0}$ appears with the coefficient
$(\lef-1)\lef^{-(k-2)}$, and the values  $\chi\vert_{U_1}$, 
$\chi\vert_{U_{\eps}}$ each appear with the coefficient 
$\frac12(\lef-1)^2\lef^{-(k-1)}$, which leads to the cancellation in
the case $\chi=Q_T$.
Finally, we are again in a situation similar to all the previous
cases, where we only need to sum over all $\lambda$ not exceeding
$n/2$. We omit the details of getting the answers that appear in the
first rows of  Tables~\ref{table:m-char.0.s} and ~\ref{table:m-char.1.s}.
This ends the proof of Proposition \ref{prop:m-ch}.
\end{proof}

\begin{table}[htbp]
\caption{Virtual motives for the characters of $\pi(0,\theta)$, $\pi(0,+)$ and $\pi(0,-)$ at $Y_{z,n}(u)$ for $z\in \{s_0,s_1,s_2\}$.}
  \centering
  \begin{tabular}{@{}| c | c | c | c | @{}}
\hline
&&&\\
	$z$ 
	& $M_{z,n,\nu}^{(0),0} $ 
	& $M_{z,n,\nu}^{(0),1}$ 
	& $M_{z,n,\nu}^{(0),\eps}$ \\
&&&\\
\hline
&&&\\
${s_0}$ & $\frac{\lef^{n-1}-1}{\lef-1}$ & $\lef^n$ &$\lef^n$\\
&&&\\
\hline
&&&\\
	${s_1}$ 
	& $\begin{matrix}\frac{\lef^{n}-1}{\lef-1}, & n \text{ odd}\\  \frac{\lef^{n-1}-1}{\lef-1}, &n \text{ even}
\end{matrix}$ 
 	& $\begin{matrix} 0, & n \text{ odd}\\ \frac{1}{2}(\lef+1)\lef^{n-1}, & n \text{ even} \end{matrix}$ 
 	& $\begin{matrix} 0, & n \text{ odd}\\ \frac{1}{2}(\lef+1)\lef^{n-1}, & n \text{ even} \end{matrix}$  \\
&&&\\
	${s_2}$ 
	&  $\begin{matrix} \frac{\lef^{n}-1}{\lef-1}, &n \text{ odd}\\ \frac{\lef^{n-1}-1}{\lef-1}, & n \text{ even}
\end{matrix}$  
	& $\begin{matrix} 0, & n\text{ odd} \\ \frac{1}{2}(\lef+1)\lef^{n-1}, & n \text{ even} \end{matrix}$ 
	& $\begin{matrix} 0, & n\text{ odd}\\ \frac{1}{2}(\lef+1)\lef^{n-1}, & n \text{ even} \end{matrix}$ \\ 
&&&\\
\hline
  \end{tabular}
\label{table:m-char.0.s}  
\end{table}

\begin{table}[htbp]
\caption{Virtual motives for the characters of $\pi(0,\theta)$, $\pi(0,+)$ and $\pi(0,-)$ at $Y_{z,n}(u)$ for $z\in \{t_0,t_1,t_2,t_3\}$. (Recall that $\zeta^2={\sgn}(-1)$.)}
  \centering
  \begin{tabular}{@{}| c | c | c | c | @{}}
\hline
&&&\\
	$z$ 
	& $M_{z,n,\nu}^{(0),0} $ 
	& $M_{z,n,\nu}^{(0),1}$ 
	& $M_{z,n,\nu}^{(0),\eps}$ \\
&&&\\
\hline
&&&\\
	${t_0}$ 
	& $\frac{\lef^{n}-1}{\lef-1}$
	& $\begin{matrix} \lef^n, &\zeta^{2n}= \nu\\ 0, & \text{otherwise}\end{matrix}$
	& $\begin{matrix} 0, & \zeta^{2n}=\nu \\ \lef^n, & \text{otherwise}\end{matrix}$ \\
&&&\\
	${t_1}$  
	& $\frac{\lef^{n}-1}{\lef-1}$ 
	& $\begin{matrix} \lef^n, &\zeta^{2n}=-\nu\\ 0, & \text{otherwise}\end{matrix}$
	& $\begin{matrix} 0, & \zeta^{2n}=-\nu\\ \lef^n, & \text{otherwise}\end{matrix}$ \\
&&&\\
	${t_2}$  
	&  $\frac{\lef^{n}-1}{\lef-1}$
	&  $\begin{matrix} \lef^n, & \zeta^{2n}=(-1)^n\nu\\ 0, & \text{otherwise}\end{matrix}$
	&  $\begin{matrix} 0, & \zeta^{2n}=(-1)^n\nu\\ \lef^n, & \text{otherwise}\end{matrix}$ \\
&&&\\
	${t_3}$  
	& $\frac{\lef^{n}-1}{\lef-1}$ 
	& $\begin{matrix} \lef^n, & \zeta^{2n}=(-1)^{n+1}\nu\\ 0, & \text{otherwise}\end{matrix}$
	& $\begin{matrix} 0, & \zeta^{2n}=(-1)^{n+1}\nu\\ \lef^n, & \text{otherwise}\end{matrix}$\\
&&&\\
\hline
  \end{tabular}
\label{table:m-char.0.t}  
\end{table}

\begin{table}[htbp]
\caption{Virtual motives for the characters of $\pi(1,\theta)$, $\pi(1,+)$ and $\pi(1,-)$ at $Y_{z,n}(u)$ for $z\in \{{ s_0,s_1,s_2}\}$.}
  \centering
  \begin{tabular}{@{}| c | c | c | c | @{}}
\hline
&&&\\
	$z$ 
	& $M_{z,n,\nu}^{(1),0} $ 
	& $M_{z,n,\nu}^{(1),1}$ 
	& $M_{z,n,\nu}^{(1),\eps}$ \\
&&&\\
\hline
&&&\\
${s_0}$ & $\frac{\lef^{n-1}-1}{\lef-1}$ & $\lef^n$ &$\lef^n$\\
&&&\\
\hline
&&&\\
	${s_1}$ 
	& $\begin{matrix}\frac{\lef^{n}-1}{\lef-1}, & n \text{ even}\\  \frac{\lef^{n-1}-1}{\lef-1}, &n \text{ odd}
\end{matrix}$ 
 	& $\begin{matrix} 0, & n \text{ even}\\ \frac{1}{2}(\lef+1)\lef^{n-1}, & n \text{ odd} \end{matrix}$ 
 	& $\begin{matrix} 0, & n \text{ even}\\ \frac{1}{2}(\lef+1)\lef^{n-1}, & n \text{ odd} \end{matrix}$  \\
&&&\\
	${s_2}$ 
	&  $\begin{matrix} \frac{\lef^{n}-1}{\lef-1}, &n \text{ even}\\ \frac{\lef^{n-1}-1}{\lef-1}, & n \text{ odd}
\end{matrix}$  
	& $\begin{matrix} 0, & n\text{ even} \\ \frac{1}{2}(\lef+1)\lef^{n-1}, & n \text{ odd} \end{matrix}$ 
	& $\begin{matrix} 0, & n\text{ even}\\ \frac{1}{2}(\lef+1)\lef^{n-1}, & n \text{ odd} \end{matrix}$ \\ 
&&&\\
\hline
  \end{tabular}
\label{table:m-char.1.s}  
\end{table}

\begin{table}[htbp]
\caption{Virtual motives for the characters of $\pi(1,\theta)$, $\pi(1,+)$ and $\pi(1,-)$ at $Y_{z,n}(u)$ for $z\in \{t_0,t_1,t_2,t_3\}$.}
  \centering
  \begin{tabular}{@{}| c | c | c | c | @{}}
\hline
&&&\\
	$z$ 
	& $M_{z,n,\nu}^{(1),0} $ 
	& $M_{z,n,\nu}^{(1),1}$ 
	& $M_{z,n,\nu}^{(1),\eps}$ \\
&&&\\
\hline
&&&\\
	${t_0}$ 
	& $\frac{\lef^{n}-1}{\lef-1}$
	& $\begin{matrix} \lef^n, &\zeta^{2n+2}= \nu\\ 0, & \text{otherwise}\end{matrix}$
	& $\begin{matrix} 0, & \zeta^{2n+2}=\nu \\ \lef^n, & \text{otherwise}\end{matrix}$ \\
&&&\\
	${t_1}$  
	& $\frac{\lef^{n}-1}{\lef-1}$ 
	& $\begin{matrix} \lef^n, & \zeta^{2n}=-\nu\\ 0, & \text{otherwise}\end{matrix}$
	& $\begin{matrix} 0, & \zeta^{2n}=-\nu\\ \lef^n, & \text{otherwise}\end{matrix}$ \\
&&&\\
	${t_2}$  
	&  $\frac{\lef^{n}-1}{\lef-1}$
	&  $\begin{matrix} \lef^n, & \zeta^{2n+2}=(-1)^{n+1}\nu\\ 0, & \text{otherwise}\end{matrix}$
	&  $\begin{matrix} 0, & \zeta^{2n+2}=(-1)^{n+1}\nu\\ \lef^n, & \text{otherwise}\end{matrix}$ \\
&&&\\
	${t_3}$  
	& $\frac{\lef^{n}-1}{\lef-1}$ 
	& $\begin{matrix} \lef^n, & \zeta^{2n+2}=(-1)^{n}\nu\\ 0, & \text{otherwise}\end{matrix}$
	& $\begin{matrix} 0, & \zeta^{2n+2}=(-1)^{n}\nu\\ \lef^n, & \text{otherwise}\end{matrix}$\\
&&&\\
\hline
  \end{tabular}
\label{table:m-char.1.t}  
\end{table}

\section{Motives for the Fourier transforms of our orbital integrals}\label{section:frobenius orbital}

In this section we use the notation of \cite{CH}.

\subsection{The Fourier transform of good orbital integrals}\label{subsection:good}

For any rational number $s$, let $\gKq_s$ denote the union of the 
{Moy-Prasad lattices} 
$\gKq_{x,s}$ as $x$ ranges over all points in the extended Bruhat-Tits building $I(\Ksch{G},\Kq)$ for $\GKq$ (see \cite{MP} for the definition of $\gKq_{x,s}$). Let $\hecke{\gKq}$ denote the Hecke algebra of locally constant, compactly supported functions $f: \gKq \to \EE$. { As in \cite[\S 1.3]{CH}, f}or any pair of rational numbers $s \leq r$ { we write} $\hecke{\gKq}^s_r$ for the $\EE$-vector space of elements of $\hecke{\gKq}$ such that $f$ is supported by $\gKq_s$ and $\hat{f}$ is supported by $\gKq_{-r}$. (Recall that the Fourier transform is taken with respect to an additive character of $\Kq$ with conductor $\Rq$.) Then the Fourier transform defines an isomorphism of $\EE$-vector spaces
\begin{eqnarray*}
	\hecke{\gKq}^s_r &\to& \hecke{\gKq}^{-r}_{-s}\\
	f &\mapsto& \hat{f}.
\end{eqnarray*}
{ Again following \cite[\S 1.3]{CH}, } we write $\hecke{\gKq}^s$ for the union of the spaces $\hecke{\gKq}^s_r$ with $s\leq r$, and $\hecke{\gKq}_r$ for the union of the spaces $\hecke{\gKq}^s_r$ with $s\leq r$.

For any $\varphi: \qgkq_{x,r}\to \EE$, we write $\varphi_{x,r}$ for the 
element of $\hecke{\gKq}^r_r$ such that $\varphi_{x,r}(Y) = (\varphi\circ\rho_{x,r})(Y)$ if $Y\in \gKq_{x,r}$ and $\varphi_{x,r}(Y) =0$ otherwise. As explained in \cite[\S 1]{CH}, on the level of reductive quotients we have another Fourier transform taking functions on $\qgkq_{x,r}$ to functions on $\qgkq_{x,-r}$. With respect to these definitions we have
\begin{equation}
	\widehat{\varphi_{x,r}} = \vol(\gKq_{x,r}) \hat\varphi_{x,-r},
\end{equation}
where $\vol$ refers to the measure on $\gKq$. For elaboration and proofs, the reader is referred to \cite[\S 1]{CH}.

{ Before stating the next proposition, we remind the reader that if $X$ is regular elliptic, then there is unique point $x$ in the Bruhat-Tits building for $\GKq$ corresponding to the centraliser $X$ in $\GKq$, since $\GKq$ has compact centre. Moreover, in this case, the depth of $X$ in $\gKq$ is the unique real number $r$ (rational number, actually) such that $X\in \gKq_{x,r}$ and $X{\not\in}\gKq_{x,r^+}$.}

\begin{proposition}\label{thm:good}
Suppose $X$ is a regular, elliptic, good element of $\gKq$. Let $x$ be the point { in the Bruhat-Tits building for $\GKq$} corresponding to the centraliser of $X$ in $\GKq$ and let $r$ be the depth of $X$.  Let $\bar{X}$ denote the image of $X$ under $\rho_{x,r}: \gKq_{x,r} \to \qgkq_{x,r}$ and let $\varphi: \qgkq_{x,r} \to \Q$ denote the characteristic function of the $\qGkq_x$-orbit of $\bar{X} \in \qgkq_{x,r}$ divided by the cardinality of that orbit. If $f\in \hecke{\gKq}^{-r}$ then
\begin{equation}\label{F.t.}
	\foi_{X}(f) 	
	=  
	\int\limits_{\GKq} \int\limits_{\gKq} f(\Ad(g)Y) \hat\varphi_{x,-r}(Y)\, dY\, dg.
\end{equation}
\end{proposition}

\begin{proof}
Suppose $f\in \hecke{\gKq}^{-r}$. Then $f\in \hecke{\gKq}_{-r}^{-s}$ for some some $-r\leq -s$. Thus, $\hat{f} \in \hecke{\gKq}^{s}_{r}$, so $\hat{f}\in \hecke{\gKq}_{r}$. Now, by \cite[Prop 1.22]{CH} and elementary properties of the Fourier transform,
\begin{equation}\label{oi}
\begin{aligned}
	&\foi_{X}(f) \\
	&= \oi_{X}(\hat{f})\\
	&= {\vol(\gKq_{x,r})}^{-1} 
\int\limits_{\GKq} \int\limits_{\gKq} \hat{f}(\Ad(g)Y) \varphi_{x,r}(Y)\, dY\, dg\\
	&= {\vol(\gKq_{x,r})}^{-1}
\int\limits_{\GKq} \int\limits_{\gKq} f(\Ad(g)Y) \widehat{\varphi_{x,r}}(Y)\, dY\, dg.
\end{aligned}
\end{equation}
By \cite[Prop 1.13]{CH}, 
\[
\widehat{\varphi_{x,r}}(Y) = \vol(\gKq_{x,r}) \hat\varphi_{x,-r}(Y),
\] 
so
\begin{equation}
	\foi_{X}(f) 	
	=  
	\int\limits_{\GKq} \int\limits_{\gKq} f(\Ad(g)Y) \hat\varphi_{x,-r}(Y)\, dY\, dg,
\end{equation}
as claimed.
\end{proof}

\begin{remark}\label{remark:good}
We will sometimes write $\varphi_{\bar X}$ (resp. $\hat{\varphi}_{\bar X}$) for the function $\varphi$ (resp. $\hat{\varphi}$) appearing in Proposition~\ref{thm:good} above; in that case, $\varphi_{x,r}$ becomes $(\varphi_{\bar X})_{x,r}$ and $\hat\varphi_{x,-r}$ becomes $(\hat\varphi_{\bar X})_{x,-r}$.
\end{remark}

\subsection{Application to our orbital integrals}\label{subsection: applications}

In order to apply \cite[Prop 1.22]{CH} to the Lie algebra $\gKq$ and the orbital integrals appearing in Theorem~\ref{thm:C}, we must find the function $(\hat\varphi_{\bar{X}_z})_{x_z,-r_z}$ for each $X_z$ appearing in Theorem~\ref{thm:C}, where $x_z$ is the point in the Bruhat-Tits building for $X_z$ and $r_z$ is the depth of $X_z$ in $\gKq$.

\subsubsection{Case: $z=s_1$}\label{subsub:orbital-s1}
Recall (from Section~\ref{subsection:Xz}) that
\begin{equation*}
X_{s_1}(v) \ceq \begin{bmatrix} 0 & v  \\  \eps v  & 0 \end{bmatrix}.
\end{equation*}  
{Here we will assume that $v$ is a unit.}
The point $x_{s_1}$ is the standard vertex of the Bruhat-Tits building for $\GKq$ and the depth $r_{s_1}$ of $X_{s_1}(v)$ is $0$; in other words, $x_{s_1} = (0)$ and $r_{s_1} = 0$. Thus,
\begin{equation*}
\gKq_{x_{s_1},r_{s_1}} = \gKq_{(0),0} = \left\{ \begin{bmatrix} z  & x \\ y  & -z   \end{bmatrix} \tq x,y,z \in \Rq \right\}.
\end{equation*}
The reduction map $\rho_{x_{s_1},r_{s_1}}$ is given by
\begin{equation*}
\begin{bmatrix} z  & x  \\ y  & - z  \end{bmatrix} \mapsto \begin{bmatrix} \bar{z} & \bar{x} \\ \bar{y} & -\bar{z}\end{bmatrix},
\end{equation*}
where $\bar{x}$, $\bar{y}$ and $\bar{z}$ denote the image of $x$, $y$ and $z$ respectively under $\Rq \to \kq$. Let $\bar{X}_{s_1}(v)$ denote the image of $X_{s_1}(v)$ under $\rho_{x_{s_1},r_{s_1}}$.
The $\quo{G}_{x_{s_1}}(\kq)$-orbit of $\bar{X}_{s_1}(v)$ in $\quo{\g}_{x_{s_1},r_{s_1}}$ is  
\begin{equation*}
\left\{  \begin{bmatrix} z & x \\ y & -z \end{bmatrix} \in { \fsl(2,}\kq) \Big\vert xy+z^2 = \bar\eps \bar{v}^2 \right\},
\end{equation*}
which has cardinality $q(q-1)$. Thus, $\varphi_{\bar{Y}_{s_1}(v)} : \quo{\g}_{x_{s_1},r_{s_1}} \to \EE$ is given by
\begin{equation}\label{eqn:s1.0}
\varphi_{\bar{X}_{s_1}(v)}\left(\begin{bmatrix} z & x \\ y & -z \end{bmatrix}\right) = 
	\begin{cases} 
	\frac{1}{q(q-1)} & xy+z^2 = \bar\eps \bar{v}^2 \\
	0 & \text{\ otherwise}.
	\end{cases}
\end{equation}

In order to find the (relative) Fourier transform of this function (in the sense of \cite{CH}) we observe that the Killing form $\langle X,Y\rangle \ceq \trace(XY)$ gives a pairing between lattices
\begin{eqnarray*}
\gKq_{x_{s_1},r_{s_1}} \times \gKq_{x_{s_1},-r_{s_1}} &\to& \Rq\\
\left( \begin{bmatrix} z & x \\ y  & -z  \end{bmatrix} , \begin{bmatrix} c  & a  \\ b  & -c  \end{bmatrix} \right) &\mapsto& xb+ya+2zc,
\end{eqnarray*}
which in turn gives a bilinear form $\qgkq_{x_{s_1},r_{s_1}} \times \qgkq_{x_{s_1},-r_{s_1}} \to \kq$. The relative Fourier transform is taken with respect to this form.
The image of the set of topologically nilpotent elements in $\gKq_{x_{s_1},-r_{s_1}}$ under $\rho_{x_{s_1},-r_{s_1}}$ is the cone
\begin{equation*}
\left\{  \begin{bmatrix} c & a \\ b & -c \end{bmatrix} \in { \fsl(2,}\kq) \Big\vert ab+c^2 = 0 \right\}.
\end{equation*}
In Section~\ref{subsection: finite} we will see that, if $ab+c^2=0$ then
\begin{equation}\label{eqn:s1.1}
Q_{\ksch{T}}\left(\begin{bmatrix} c & a \\ b & -c \end{bmatrix}\right)
=
(1-q)\hat\varphi_{\bar{X}_{s_1}(v)}\left(\begin{bmatrix} c & a \\ b & -c \end{bmatrix}\right),
\end{equation}
where $Q_\ksch{T}$ is given in equation (\ref{eqn:ff7.T}). This completes our description of the relevant properties of $\hat\varphi_{\bar{X}_{s_1}(v)}$. 

\subsubsection{Case: $z=s_2$}\label{subsub:orbital-s2}
Since this case is very similar to the case above, we only summarize the results here.
Recall (from Section~\ref{subsection:Xz}) that
\begin{equation*}
X_{s_2}(v) \ceq \begin{bmatrix} 0 & v \varpi \\  \eps v \varpi^{-1} & 0 \end{bmatrix},
\end{equation*}  
{ where $v$, as above, is assumed to be a unit.}
The point $x_{s_2}$ is the vertex $(1)$ of the Bruhat-Tits building for $\GKq$ (see Section~\ref{subsection: cuspidal}) and the depth $r_{s_2}$ of $X_{s_2}(v)$ is $0$. Thus,
\begin{equation*}
\gKq_{x_{s_2},r_{s_2}} = \gKq_{(1),0} = \left\{ \begin{bmatrix} z  & x \varpi \\ y \varpi^{-1} & -z   \end{bmatrix} \tq x,y,z \in \Rq \right\}.
\end{equation*}
The reduction map $\rho_{x_{s_2},r_{s_2}}$ is given by
\begin{equation*}
\begin{bmatrix} z   & x \varpi \\ y \varpi^{-1} & - z  \end{bmatrix} \mapsto \begin{bmatrix} \bar{z} & \bar{x} \\ \bar{y} & -\bar{z}\end{bmatrix}.
\end{equation*}
The function $\varphi_{\bar{X}_{s_2}(v)} : \quo{\g}_{x_{s_2},r_{s_2}} \to \EE$ is exactly as in the preceding case, so
\begin{equation}\label{eqn:s2.1}
Q_{\ksch{T}}\left(\begin{bmatrix} c & a \\ b & -c \end{bmatrix}\right)
=
(1-q)\hat\varphi_{\bar{X}_{s_2}(v)}\left(\begin{bmatrix} c & a \\ b & -c \end{bmatrix}\right),
\end{equation}
as above. This completes our description of the relevant properties of $\hat\varphi_{\bar{X}_{s_2}(v)}$. 

\subsubsection{Case: $z\in \{t_0,t_1,t_2,t_3\}$.} Then the point $x_z$ is $(01)$ and the depth $r_z$ is $\frac{1}{2}$. Thus, 
\begin{equation*}
\gKq_{x_z,r_z} = \gKq_{(01),\frac{1}{2}} = \left\{ \begin{bmatrix} z \varpi & x \\ y \varpi & -z \varpi  \end{bmatrix} \tq x,y,z \in \Rq \right\}.
\end{equation*}
Thus, $\quo{\g}_{x_z,r_z} = \A^2(\kq)$ and the reduction map $\rho_{x_z,r_z}: \gKq_{x_z,r_z} \to \A^2(\kq)$ is given by
\begin{eqnarray*}
\begin{bmatrix} z \varpi  & x  \\ y \varpi & - z \varpi \end{bmatrix} &\mapsto& (\bar{x},\bar{y}).
\end{eqnarray*}
The reduction map on $\GKq_{x_z} \to {\GL}(1,\kq)$ is given by
\begin{eqnarray*}
\begin{bmatrix} a & b \\ \varpi c & d \end{bmatrix} &\mapsto& \bar{a}.
\end{eqnarray*}
Thus, the action of $\quo{G}_{x_z}$ on $\quo{\g}_{x_z, r_z}$ corresponds to the action of ${\GL}(1,\kq)$ on $\A^2(\kq)$ given by $t\cdot(x,y) \ceq (t^2x,t^{-2}y)$. It follows immediately from the definitions above that
\begin{eqnarray*}
\varphi_{\bar{X}_{t_0}(v)}(x,y) &=& 
	\begin{cases} 
	\frac{2}{q-1} & xy = \bar{v}^2 \AND {\sgn}(x) = {\sgn}(v) \\
	0 & \text{\ otherwise}
	\end{cases}\\
\varphi_{\bar{X}_{t_1}(v)}(x,y) &=& 
	\begin{cases} 
	\frac{2}{q-1} & xy = \bar{v}^2 \AND {\sgn}(x) = {\sgn}(\eps v) \\
	0 & \text{\ otherwise}
	\end{cases}\\
\varphi_{\bar{X}_{t_2}(v)}(x,y) &=& 
	\begin{cases} 
	\frac{2}{q-1} & xy = \bar{\eps} \bar{v}^2 \AND {\sgn}(x) = {\sgn}(v) \\
	0 & \text{\ otherwise}
	\end{cases}\\
\varphi_{\bar{X}_{t_3}(v)}(x,y) &=& 
	\begin{cases} 
	\frac{2}{q-1} & xy = \bar{\eps} \bar{v}^2 \AND {\sgn}(x) = {\sgn}(\eps v) \\
	0 & \text{\ otherwise}.
	\end{cases}
\end{eqnarray*}

\subsection{Two functions on $\A^2(\kq)$}\label{subsection: two functions}

In this section we introduce two functions that play a crucial role in our proof of Theorem~\ref{thm:C}. Using notation from the preceding subsection, define
\begin{equation}\label{eqn:varphi01}
\begin{aligned}
\varphi_{(0)} &\ceq (\varphi_{\bar{X}_{t_0}} - \varphi_{\bar{X}_{t_1}}) + (\varphi_{\bar{X}_{t_2}}  - \varphi_{\bar{X}_{t_3}}) \\
\varphi_{(1)} &\ceq (\varphi_{\bar{X}_{t_0}} - \varphi_{\bar{X}_{t_1}}) - ( \varphi_{\bar{X}_{t_2}}  - \varphi_{\bar{X}_{t_3}}).
\end{aligned}
\end{equation}
To find the Fourier transform of these functions (in the sense of \cite{CH}), observe that the Killing form $\langle X,Y\rangle \ceq \trace(XY)$ gives a pairing between lattices
\begin{eqnarray*}
\gKq_{x_{z},r_{z}} \times \gKq_{x_{z},-r_{z}} &\to& \Rq\\
\left( \begin{bmatrix} z \varpi & x  \\ y \varpi & -z \varpi  \end{bmatrix} , \begin{bmatrix} c  & a \varpi^{-1} \\ b  & -c   \end{bmatrix} \right) &\mapsto& xb+ya + 2 \varpi zc,
\end{eqnarray*}
which in turn gives the bilinear form 
\begin{eqnarray*}
\qgkq_{x_{z},r_{z}} \times \qgkq_{x_{z},-r_{z}} &\to& \kq\\
\left( (x,y), (a,b) \right) &\mapsto& xb+ya.
\end{eqnarray*}
The relative Fourier transform is taken with respect to this form.
The image of the set of topologically nilpotent elements in $\gKq_{x_{z},-r_{z}}$ under $\rho_{x_{z},-r_{z}}$ is the normal crossing
\begin{equation*}
\left\{  (a,b) \in \A^2(\kq) \tq  ab = 0 \right\}.
\end{equation*}
In Section~\ref{subsection: finite} we will see that
\begin{equation}\label{eqn:hatvarphi01}
\begin{aligned}
\hat\varphi_{(0)}(a,0) &= 0\\
\hat\varphi_{(0)}(0,b) &= \frac{2^2}{q-1}\sqrt{q}\zeta^{3} {\sgn}(b) \\
\hat\varphi_{(1)}(a,0) &= \frac{2^2}{q-1}\sqrt{q}\zeta {\sgn}(a) \\
\hat\varphi_{(1)}(0,b) &= 0.
\end{aligned}
\end{equation}

\subsection{Some finite field calculations}\label{subsection: finite}

In this section we defend equations (\ref{eqn:s1.1}), (\ref{eqn:s2.1}) and (\ref{eqn:hatvarphi01}).

Recall the definition of the function $\gauss_\pm : \A^2(\kq) \to \EE$ from equation~(\ref{eqn: gauss}). Recall also that we equip $\A^2(\kq)$ with the bilinear form $\langle (x,y),(a,b)\rangle = xb + ya$, as explained in Section~\ref{subsection: applications}. Finally, recall the definition of $\varphi_{\bar{X}_{z}(v)}$ for $z\in \{ t_0,t_1,t_2,t_3\}$.
Then
\begin{eqnarray*}
\hat\varphi_{\bar{X}_{t_0}(v)}(a,b) 
	&=& \sum_{(x,y)\in \A^2(\kq)} \bar\psi \langle(x,y),(a,b)\rangle\ \varphi_{\bar{X}_{t_0}(v)}(x,y)\\
	&=& \frac{2}{q-1}  \sum_{xy= \bar{v}^2 \atop {\sgn}(x)={\sgn}(v) } \bar\psi(xb)\bar\psi(ya)\\
	&=& \frac{2}{q-1} \sum_{{\sgn}(x)={\sgn}(v) } \bar\psi(xb) \bar\psi(\bar{v}^2 x^{-1}a).
\end{eqnarray*}
Thus,
\begin{eqnarray*}
\hat\varphi_{\bar{X}_{t_0}(v)}(a,0) 
	&=& \frac{2}{q-1} \sum_{{\sgn}(x)={\sgn}(v) }\bar\psi(0) \bar\psi(\bar{v}^2 x^{-1}a)\\
	&=& \frac{2}{q-1} \sum_{{\sgn}(t)={\sgn}(v) } \bar\psi(ta)\\
	&=& \frac{2}{q-1}\gauss_{{\sgn}(v)}(a),
\end{eqnarray*}
and
\begin{eqnarray*}
\hat\varphi_{\bar{X}_{t_0}(v)}(0,b) 
	&=& \frac{2}{q-1} \sum_{{\sgn}(x)={\sgn}(v) }\bar\psi(xb) \bar\psi(0)\\
	&=& \frac{2}{q-1}\gauss_{{\sgn}(v)}(b).
\end{eqnarray*}
Similar arguments show that
\begin{eqnarray*}
\hat\varphi_{\bar{X}_{t_1}(v)}(a,0) &=& \frac{2}{q-1}\gauss_{{\sgn}(\eps v)}(a)\\
\hat\varphi_{\bar{X}_{t_1}(v)}(0,b) &=& \frac{2}{q-1}\gauss_{{\sgn}(\eps v)}(b)
\end{eqnarray*}
and
\begin{eqnarray*} 
\hat\varphi_{\bar{X}_{t_2}(v)}(a,0) &=& \frac{2}{q-1}\gauss_{{\sgn}(\eps v)}(a)\\
\hat\varphi_{\bar{X}_{t_2}(v)}(0,b) &=& \frac{2}{q-1}\gauss_{{\sgn}(v)}(b)
\end{eqnarray*}
and
\begin{eqnarray*}
\hat\varphi_{\bar{X}_{t_3}(v)}(a,0) &=& \frac{2}{q-1}\gauss_{{\sgn}(v)}(a)\\
\hat\varphi_{\bar{X}_{t_3}(v)}(0,b) &=& \frac{2}{q-1}\gauss_{{\sgn}(\eps v)}(b).
\end{eqnarray*}
Thus,
\begin{eqnarray*}
	&& \hskip-20pt \left(\hat\varphi_{\bar{X}_{t_0}(v)}-\hat\varphi_{\bar{X}_{t_1}(v)}+\hat\varphi_{\bar{X}_{t_2}(v)}-\hat\varphi_{\bar{X}_{t_3}(v)}\right)(a,0)\\
	&=& \frac{2}{q-1}\left( \gauss_{{\sgn}(v)} - \gauss_{{\sgn}(\eps v)} + \gauss_{{\sgn}(\eps v)} - \gauss_{{\sgn}(v)}\right)(a)\\
	&=& 0,
\end{eqnarray*}
and
\begin{eqnarray*}
	&& \hskip-20pt \left(\hat\varphi_{\bar{X}_{t_0}(v)}-\hat\varphi_{\bar{X}_{t_1}(v)}+\hat\varphi_{\bar{X}_{t_2}(v)}-\hat\varphi_{\bar{X}_{t_3}(v)}\right)(0,b)\\
	&=& \frac{2}{q-1}\left( \gauss_{{\sgn}(v)} - \gauss_{{\sgn}(\eps v)} + \gauss_{{\sgn}(v)} - \gauss_{{\sgn}(\eps v)}\right)(b)\\
	&=& \frac{2^2}{q-1} (\gauss_{{\sgn}(v)} - \gauss_{{\sgn}(\eps v)})(b)\\
	&=& \frac{2^2}{q-1} {\sgn}(v) \widehat{\sgn}(b)\\
	&=& \frac{2^2}{q-1} {\sgn}(v) \sqrt{q} \zeta^{3} {\sgn}(b).
\end{eqnarray*}
Letting $v=1$ we recover the first two parts of equation (\ref{eqn:hatvarphi01}).
Likewise,
\begin{eqnarray*}
	&& \hskip-20pt \left(\hat\varphi_{\bar{X}_{t_0}(v)}-\hat\varphi_{\bar{X}_{t_1}(v)}-\hat\varphi_{\bar{X}_{t_2}(v)}+\hat\varphi_{\bar{X}_{t_3}(v)}\right)(a,0)\\
	&=& \frac{2}{q-1}\left( \gauss_{{\sgn}(v)} - \gauss_{{\sgn}(\eps v)} - \gauss_{{\sgn}(\eps v)} + \gauss_{{\sgn}(v)}\right)(a)\\
	&=& \frac{2^2}{q-1} (\gauss_{{\sgn}(v)} - \gauss_{{\sgn}(\eps v)})(a)\\
	&=& \frac{2^2}{q-1} {\sgn}(v) \widehat{\sgn}(a)\\
	&=& \frac{2^2}{q-1} {\sgn}(v) \sqrt{q} \zeta^{1} {\sgn}(a),
\end{eqnarray*}
and
\begin{eqnarray*}
	&& \hskip-20pt \left(\hat\varphi_{\bar{X}_{t_0}(v)}-\hat\varphi_{\bar{X}_{t_1}(v)}-\hat\varphi_{\bar{X}_{t_2}(v)}+\hat\varphi_{\bar{X}_{t_3}(v)}\right)(0,b)\\
	&=& \frac{2}{q-1}\left( \gauss_{{\sgn}(v)} - \gauss_{{\sgn}(\eps v)} - \gauss_{{\sgn}(v)} + \gauss_{{\sgn}(\eps v)}\right)(b)\\
	&=& 0.
\end{eqnarray*}
Letting $v=1$ we recover the last two parts of equation (\ref{eqn:hatvarphi01}).

Equation~(\ref{eqn:s1.1}) (and therefore equation (\ref{eqn:s2.1}) also) is exactly Springer's Hypothesis in our context. In fact, it is quite easy to verify this equality by direct calculation of the relevant Fourier transforms, but since this paper is already long, we omit the details.

\subsection{Motives for the Fourier transforms of orbital integrals}\label{subsection:m-oi}

As we see from Section~\ref{subsection: applications}, for our purposes, there are exactly three important points in the Bruhat-Tits building: $(0)$, $(01)$ and $(1)$.  For the remainder of this section, the symbol $x$ is reserved for the vertices $(0)$ and $(1)$ while $y=(01)$. 

Let $f$ be a test function supported on the topologically nilpotent elements in $\g(\Kq)$. All we need to know in order to find $\foi_{X_z}(f)$ is the values of the corresponding function $\hat\varphi_{X_z}$ on the corresponding reductive quotient, and the information about the projections of the elements of the form $h^{-1}Yh$ for the elements $Y$ in the support of $f$.

Recall the notation  $\h_{z,n}=\{Y_{z,n}(u)\mid u\in\Uq\}\subset\gKq$,
 and $\tilde f_{z,n,\pm}=\cay^{\ast}f_{z,n,\pm}$ (the characteristic function 
of the set ${}^G\h_{z,n,\pm}$,  where $\GKq$ acts by adjoint action).

\begin{proposition}\label{proposition: frobenius orbits}
Suppose {$z\in \{ s_0, s_1, s_2, t_0, t_1, t_2, t_3\}$}.
Let $\tilde f_{z,n,\pm}$ be the characteristic function of the set
${}^G\h_{z,n,\pm}$ and let $\mathfrak{u}_0$, $\mathfrak{u}_1$ and $\mathfrak{u}_\eps$ denote the nilpotent orbits in ${ \fsl(2,}\kq)$ corresponding to the unipotent conjugacy classes $U_0$, $U_1$ and $U_\eps$ respectively.
Let $\kappa$ be an arbitrary linear combination of the functions $Q_T$ and $Q_G$ viewed as functions on ${ \fsl(2,}\kq)$ 
(see Section~\ref{subsection: finite} for their definition). 
Then
\begin{equation}
\begin{aligned}
&\int_G\int_{\gKq}\tilde f_{z,n,\pm}(gYg^{-1})\kappa_{x,0}(Y)dYdg\\
&= \kappa\vert_{\mathfrak{u}_0}\frob M_{z,n,\pm}^{x,0} +\kappa\vert_{\mathfrak{u}_1}\frob M_{z,n,\pm}^{x,1} +\kappa\vert_{\mathfrak{u}_{\eps}}\frob M_{z,n,\pm}^{x,\eps},
\end{aligned}
\end{equation}
where $M_{z,n,\pm}^{x,0}$, $M_{z,n,\pm}^{x,1}$ and $M_{z,n,\pm}^{x,\eps}$ are defined in Tables~\ref{table:m-char.0.s} and \ref{table:m-char.0.t} in the case $x=(0)$ and Tables~\ref{table:m-char.1.s} and \ref{table:m-char.1.t} in the case $x=(1)$.
\end{proposition}

\begin{proof}
This Proposition essentially follows from
Proposition~\ref{prop:m-ch}. We have
\begin{equation*}
\begin{aligned}
&\int_G\int_{\gKq}\tilde f_{z,n,\pm}(gYg^{-1})\kappa_{x,0}(Y)dYdg\\
&=\int_G\int_{\gKq}\tilde f_{z,n,\pm}(Y)\kappa_{x,0}(g^{-1}Yg)dYdg\\
&=\int_G\int_{{}^G\Gamma_{z,n,\pm}}(\kappa_{x,0}\circ\cay^{\ast})(g^{-1}\gamma g)d\gamma dg.
\end{aligned}
\end{equation*} 
Note that in the above integral, the function $\kappa_{x,0}\circ \cay^{\ast}$
is evaluated only at topologically unipotent elements.
The { modified Cayley transform $\cay$} is measure-preserving on this set, and that's why the integral
can be rewritten as a double integral over the group $\GKq$.
This is exactly the expression that appears in the equation
(\ref{eq:frobenius_char}), and therefore by Proposition~\ref{prop:m-ch} 
it has the required form.  
\end{proof}
 
\begin{corollary}\label{cor: unramified orbital}
Suppose $z\in \{s_1,s_2\}$ and 
{$z'\in \{s_0, s_1, s_2, t_0,t_1,t_2,t_3\}$}.
For each $n\in \N$, the Fourier transform $\foi_{X_{z}}$ of the
orbital integral at $X_{z}$ is constant on the set
${}^G\h_{z',n,\pm}$, and 
\begin{eqnarray*}\label{eq:unramified-orbital}
&& \hskip -20pt \foi_{X_z}(\tilde f_{z',n,\pm}) \\
&=& 
\hat\varphi_{\bar{X}_{z}}\vert_{\mathfrak{u}_0}\frob M_{z',n,\pm}^{x,0} 
+\hat\varphi_{\bar{X}_{z}}\vert_{\mathfrak{u}_1}\frob M_{z',n,\pm}^{x,1} 
+\hat\varphi_{\bar{X}_{z}}\vert_{\mathfrak{u}_{\eps}}\frob M_{z',n,\pm}^{x,\eps},
\end{eqnarray*}
where $x=(0)$ if $z=s_1$ and $x=(1)$ if $z=s_2$. 
See Remark~\ref{remark:good} for the definition of  $\varphi_{\bar{X}_{z}}$.
\end{corollary}

\begin{proof}
From the Section~\ref{subsection: applications} we see that our elements
$X_{z}$ correspond to vertices $x$ in the building, and their depths
are all $r=0$.  
First, recall that by Proposition~\ref{thm:good}, 
we have 
$$
\Phi(X_{z},\tilde f_{z',n,\pm})
=\int_G\int_{\gKq}\tilde f_{z',n,\pm}(gYg^{-1})(\hat\varphi_{\bar{X}_{z}})_{x,0}(Y).
$$
Now we can plug in $\kappa=\hat\varphi_{\bar{X}_{z}}$ in the Proposition~\ref{proposition: frobenius orbits}.
By Springer's Hypothesis (see Section \ref{subsection: finite}), the function 
$\hat\varphi_{\bar{X}_{z}}$ restricted to the set of nilpotent
elements is a constant multiple of the Green's polynomial $Q_T$
(thought of as a function on the Lie algebra), so
the assumptions of the Proposition are satisfied. 
\end{proof}

Let us now consider the elements $X_z$ of the ramified elliptic tori, that is, $z=t_0$, $t_1$, $t_2$ or $t_3$. Then (see Section~\ref{subsection: applications}), $y=(01)$ and $r=\frac{1}{2}$, and $\bar\g_{y,-r}(\kq)=\A^2(\kq)$; in the rest of this section we write $y$ for $(01)$ and $r$ for $\frac{1}{2}$. The image of the set of topologically nilpotent elements in $\gKq_{y,-1/2}$ under $\rho_{y,-r}$ is contained in $\{ (x,y)\in \A^2(\kq) \tq xy=0\}$. This set is the union of the following five orbits of the action of ${\GL}(1,\kq)$ on $\A^2(\kq)$: 
\begin{equation}\label{eqn:V}
\begin{aligned}
& V^0 \ceq \{(0,0)\} &\\
&V^{1,+} \ceq \{(0, x)\mid {\sgn}(x)=1\} \qquad
 & V^{1,-} \ceq \{(0, x)\mid {\sgn}(x)=-1\}\\
& V^{2,+} \ceq \{(x, 0)\mid {\sgn}(x)=1\} \qquad
 &V^{2,-} \ceq \{(x, 0)\mid {\sgn}(x)=-1\}
 \end{aligned}
 \end{equation}

\begin{proposition}\label{orb.int}
Let $\kappa$ be a ${\GL}(1)$-invariant function defined on $\A^2(\kq)$
with respect to the action defined in Section~\ref{subsection:
  finite}. Then, for each 
{ $z\in\{s_1,s_2,t_0,t_1,t_2,t_3\}$}, each $\nu=\pm$, and each
non-negative integer $n$ in the case {
$z\in\{t_0,t_1,t_2,t_3\}$}, and each
  positive integer $n$ in the case {$z\in\{s_1,s_2\}$}, there
  exist  virtual Chow motives ${\mathcal N}_{z,n,\nu}^0$,
${\mathcal N}_{z,n,\nu}^{1,\pm}$ and ${\mathcal N}_{z,n,\nu}^{2,\pm}$, 
such that
\begin{equation*} 
\begin{aligned}
& \int_{G(\Kq)}\int_{\g(\Kq)}\tilde{f}_{z,n,\nu}(\text{Ad}(g)Y)\kappa_{y,-r}(Y)\, dY\, dg\\
&= \kappa\vert_{V_0}\frob {\mathcal N}_{z,n,\nu}^{0}+\sum_{\alpha=\pm}
\kappa\vert_{V^{1,\alpha}}\frob {\mathcal N}_{z,n,\nu}^{1,\alpha} + \sum_{\alpha=\pm}\kappa\vert_{V^{2,\alpha}}\frob {\mathcal N}_{z,n,\nu}^{2,\alpha},
\end{aligned}
\end{equation*}
Moreover, ${\mathcal N}_{z,n,\nu}^0$ and ${\mathcal N}_{z,n,\nu}^{1,\pm}$ and
${\mathcal N}_{z,n,\nu}^{2,\pm}$ are rational functions of $\lef$. 
The virtual Chow motives ${\mathcal N}_{z,n,\nu}^{1,\pm}$ and
${\mathcal N}_{z,n,\nu}^{2,\pm}$ are given in Tables~\ref{table:m-orb.s} and \ref{table:m-orb.t}.
\end{proposition}

\begin{proof}
We rewrite the left-hand side 
as the sum over the nilpotent orbits in the reductive
quotient, as we did with the character in
Section~\ref{sub:tricks}. In order to do this, we make the following definition.
Suppose $V$ is one of the ${\GL}(1,\kq)$-orbits in $\A^2(\kq)$ appearing in equation~(\ref{eqn:V}).
Let $n$ be an  integer and recall that
$y$ denotes the point $(01)$ in the Bruhat-Tits building for $\GKq$. For
any regular topologically nilpotent element $Y$ of $\gKq$, let
$N_{V,\lambda}^{y}(Y)$ denote the number of $\GKq_x$-cosets $g\GKq_y$ 
inside the double coset space $\GKq_x \backslash \GKq / \GKq_y$ that satisfy the following condition:
\[
g^{-1}Yg\in \gKq_{y,-r}
\quad \AND \quad 
\rho_{y,-r}(g^{-1}Yg)\in V.
\]

Let $\tilde{A}$ be the set of diagonal matrices of the form 
$\text{diag}(\varpi^{\lambda},\varpi^{-\lambda})$ where $\lambda$ is
an arbitrary integer. 
Then $\GKq$ has the decomposition $\GKq=\GKq_x\tilde{A} \GKq_y$ (note the difference with Cartan
decomposition where $\lambda$ is non-negative). Using the $\GKq_y$-invariance of the function
$\kappa_{y,-r}$, we obtain, for an arbitrary test function $f$:
\begin{equation}\label{int1}
\begin{aligned}
&\int_{\GKq} \int_{\gKq}f(Y)\kappa_{y,-r}(g^{-1}Yg)\,dY\,dg\\
&=\int_{\GKq/\GKq_y}\int_{\GKq_y}\int_{\gKq}f(Y)\kappa_{y,-r}(y^{-1}h^{-1}Yhy)\,dY\,dy\,dh\\
&=\mes(\GKq_y)\sum_{a\in \GKq_x\backslash \GKq /\GKq_y}
\sum_{h\in \GKq_x a\GKq_y/\GKq_y}\int_{\gKq}f(Y)\kappa_{y,-r}(h^{-1}Yh)\,dY.
\end{aligned}
\end{equation}
Note that the summation index in the outside sum in fact runs over $\Z$.  
As with the case of the character, the sum 
in fact contains only finitely many nonzero terms, since at the moment
we are 
considering only the elliptic elements $Y$.

The rest of the argument follows the pattern of the proof of Proposition~\ref{prop:m-ch}, taking equation~(\ref{int1}) as the starting point. It also proceeds case-by-case. Here we carry out the proof for the
test functions $\tilde{f}_{z,n,\pm}$ with $z=s_1$ and $z=t_2$.
 The other cases are very similar; the results of these calculations are recorded in Tables~\ref{table:m-orb.s} and \ref{table:m-orb.t}.

As in Section~\ref{sub:tricks}, we can continue the chain of
equalities (\ref{int1}) by writing the integral inside the sum as a
sum over ${\GL}(1,\F_q)$-orbits:
\begin{equation*}
\begin{aligned}
&\mes(\GKq_y)\sum_{a\in \GKq_x\backslash \GKq /\GKq_y}
\sum_{h\in \GKq_x a\GKq_y/\GKq_y}\int_{\gKq}\tilde{f}_{z,n,\pm}(Y)\kappa_{y,-r}(h^{-1}Yh)\,dY\\
&=\sum_{a\in \GKq_x\backslash \GKq /\GKq_y}
\sum_{V}\kappa\vert_{V}
\int_{{}^G\mathfrak{h}_{z,n,\pm}} N_{V,\lambda}(Y)dY \\
&=\sum_{V}\kappa\vert_{V}\sum_{\lambda=-\infty}^{\infty}\int_{{}^G\mathfrak{h}_{z,n,\pm}} N_{V,\lambda}(Y)dY . 
\end{aligned}
\end{equation*}

Note that since there are, in fact, only finitely many non-zero terms,
the permutation of the two sums is valid.

Now it remains to ``calculate'' the  numbers  $N_{V,\lambda}(Y)$, \ie, to express them in terms of  motivic volumes of some
definable sets. This is done by brute force, in a manner similar to the calculation of the character.  
Our calculation will make it transparent that these numbers are constant on each of the sets ${}^G\h_{z,n,\pm}$.

We will need the formula for the number of $\GKq_y$-cosets inside
each double coset $\GKq_x a_{\lambda} \GKq_y$ ($\lambda\in\Z$):
the cardinality $\# \GKq_x a_{\lambda} \GKq_y/\GKq_y$ equals $q^{2\lambda-1}$ 
if $\lambda>0$, and $q^{2|\lambda|}$ if $\lambda\leq 0$, as can be shown using the affine Bruhat decomposition for $\GKq$ (see \cite{BT}, for example).
Recall that with the notation of Section \ref{subsubsection: ch.m.}, we can write  
$[\GKq_x a_{\lambda}\GKq_y/\GKq_y]=\lef^{2\lambda-1}$ 
when $\lambda$ is a positive integer, 
$[\GKq_x a_{\lambda}\GKq_y/\GKq_y]=\lef^{-2\lambda}$ when $\lambda\leq 0$.
 
Let $ha_{\lambda}$ be a representative of a coset $\GKq_x a_{\lambda}
\GKq_y/\GKq_y$; write
$h=\left[\begin{smallmatrix}a&b\\c&d\end{smallmatrix}\right]\in
\GKq_x$. 

A statement completely analogous to Lemma \ref{lemma: W_U} relates
 the numbers $N_{V,\lambda}$ to 
 the motivic volumes of the sets 
$\{h\mid
\rho_{y,-r}(a_{\lambda}^{-1}h^{-1}Y_{z,n}(u)ha_{\lambda})=(x,0)
\}$, and $\{h\mid
\rho_{y,-r}(a_{\lambda}^{-1}h^{-1}Y_{z,n}(u)ha_{\lambda})=(0,x)
\}$, with $x$ square or non-square, respectively (of course, we also need to show
that these four sets are definable).
 Then to obtain the virtual motives 
${\mathcal N}_{z,n,\nu}^{1,\pm}$, and ${\mathcal N}_{z,n,\nu}^{2,\pm}$, 
we need to sum these motivic volumes over all values
of $\lambda$.

The following is the list of possibilities for
the element $\rho_{y,-r}(a_{\lambda}^{-1}h^{-1}Y_{z,n}(u)ha_{\lambda})$ of $\A^2(\kq)$ in the two cases {$z=t_2$ and $z=s_1$}:

\noindent{\bf Case: $z=t_2$.}
The conditions on $a,b,c,d$ and $\lambda$  
for $a_{\lambda}^{-1}h^{-1}Y_{z,n}(u)h a_{\lambda}$ to be in $\GKq_{y,-r}$ are:
\begin{equation}\label{conditions.orb.t2}
\begin{aligned}
-2\lambda+\ord(d^2u\varpi^n-b^2\eps u\varpi^{n+1})&\ge -1, \text{and}\\
2\lambda+\ord(-c^2u\varpi^n +a^2\eps u \varpi^{n+1})&\ge 0.
\end{aligned}
\end{equation}
Applying the reduction map $\rho_{y,-r}$ to 
$a_{\lambda}^{-1}h^{-1}Y_{z,n}(u)h a_{\lambda}$ ,
we see that there are a few possibilities:

\begin{enumerate}
\item[1.]$n$ is even: 
Recall the notation $\nu={\sgn}(u)$.
		\begin{enumerate}
		\item[(a)]
If $-n/2<\lambda\le n/2$, then $\rho_{y,-r}(a_{\lambda}^{-1}h^{-1}Y_{z,n}(u)h a_{\lambda})=(0,0)$.
		\item[(b)]
Suppose $\lambda=n/2+1$. If $\ord(d)>0$ then \[\rho_{y,-r}(a_{\lambda}^{-1}h^{-1}Y_{z,n}(u)h a_{\lambda})=(x,0)\] 
with ${\sgn}(x)={\sgn}(-\eps u)=-{\sgn}(-u)=-\zeta^2{\sgn}(u)$; on the other hand, if $\ord(d)=0$, we have $a_{\lambda}^{-1}h^{-1}Y_{z,n}(u)h a_{\lambda}\notin \GKq_{y,-r}$.
		\item[(c)] 
Suppose $\lambda=-n/2$. If $\ord(c)>0$  then \[\rho_{y,-r}(a_{\lambda}^{-1}h^{-1}Y_{z,n}(u)h a_{\lambda})=(0,0);\] on the other hand, if $\ord(c)=0$ then \[\rho_{y, -r}(a_{\lambda}^{-1}h^{-1}Y_{z,n}(u)h a_{\lambda})=(0,y),\] with ${\sgn}(y)={\sgn}(-u)=\zeta^2\nu$. 
		\end{enumerate}
 
The case (b) allows us to calculate the virtual motives responsible
for the two orbits of the form $(x,0)$: 
the virtual Chow motive {${\mathcal N}_{t_2,n,\nu}^{1,+}$} 
corresponding to
the orbit of $(1,0)$ equals $0$ if $\zeta^2\nu=1$, and equals 
$\frac{\mu(\GKq_y)}{\mu(\GKq_x)}[\GKq_x a_{\frac{n}2+1} \GKq_y/\GKq_y]
=\frac{1}{\lef+1}\lef^{n+1}$ if $\zeta^2\nu=-1$. 
For the orbit of $(\epsilon, 0)$, the answer is the reverse:  
{${\mathcal N}_{t_2,n,\nu}^{1,-}$} 
equals $0$ if  $\zeta^2\nu=-1$, and
equals $\frac{1}{\lef+1}\lef^{n+1}$ 
if $\zeta^2\nu=1$.

From the case (c), we get the motives corresponding to the other two
orbits: ${\mathcal
  N}_{t_2,n,\nu}^{2,+}=\frac{\mu(\GKq_x)-\mu(\GKq_y)}{\mu(\GKq_x)}[\GKq_x
  a_{-\frac{n}2} \GKq_y/\GKq_y] =\frac{\lef}{\lef+1}\lef^{n}$ if
$\zeta^2\nu=1$,  and $0$ if  $\zeta^2\nu=-1$, and for
${\mathcal N}_{t_2,n,\nu}^{2,-}$ these answers are reversed.

\item[2.] $n$ is odd:
		\begin{enumerate}
		\item[(a)]
If $-(n+1)/2<\lambda<(n+1)/2$, then  $\rho_{y,-r}(a_{\lambda}^{-1}h^{-1}Y_{z,n}(u)h a_{\lambda})=(0,0)$.
		\item[(b)] 
If $\lambda=(n+1)/2$, there are two possibilities: if $\ord(d)=0$ then we have  \[\rho_{y,-r}(a_{\lambda}^{-1}h^{-1}Y_{z,n}(u)h a_{\lambda})=(x,0)\] with ${\sgn}(x)={\sgn}(u)$; on the other hand, if $\ord(d)>0$ then  \[\rho_{y,-r}(a_{\lambda}^{-1}h^{-1}Y_{z,n}(u)h a_{\lambda})=(0,0).\]
		\item[(c)] 
Suppose $\lambda=-(n+1)/2$.  If $\ord(c)=0$, we have  $a_{\lambda}^{-1}h^{-1}Y_{z,n}(u)h a_{\lambda}\notin \GKq_{y,-r}$; if $\ord(c)>0$ then \[\rho_{y,-r}(a_{\lambda}^{-1}h^{-1}Y_{z,n}(u)h a_{\lambda})=(0,y),\] with ${\sgn}(y)=-{\sgn}(u)$.
		\end{enumerate}
As before, from the calculations in the sub-case (b), we see that the virtual Chow motive ${\mathcal N}_{t_2,n,\nu}^{1,+}$
that corresponds to the orbit of $(1,0)$ equals $0$ if $\nu=-1$, and equals $\frac{\mu(\GKq_x)-\mu(\GKq_y)}{\mu(\GKq_x)}[\GKq_x a_{\frac{n+1}2} \GKq_y/\GKq_y] =\frac{\lef}{\lef+1}\lef^{n+1-1}=\frac1{\lef+1}\lef^{n}$ if $\nu=1$. For the orbit of $(\epsilon, 0)$, the answer is the reverse. The subcase (c) yields  ${\mathcal N}_{t_2,n,\nu}^{2,+}=\frac{\mu(\GKq_y)}{\mu(\GKq_x)}[\GKq_x a_{-\frac{n+1}2} \GKq_y/\GKq_y]
=\frac{1}{\lef+1}\lef^{n+1}$ if $\nu=-1$,  and $0$ if $\nu=1$. 
\end{enumerate} 

\noindent{\bf Case: $z=s_1$.} In this case, the conditions
(\ref{conditions.orb.t2}) are replaced with: 
\begin{equation}\label{conditions.orb.s1}
\begin{aligned}
-2\lambda+\ord(d^2\varpi^n-b^2\eps u\varpi^{n})&\ge -1, \text{and}\\
2\lambda+\ord(-c^2u\varpi^n +a^2\eps u \varpi^{n})&\ge 0.
\end{aligned}
\end{equation}

\begin{enumerate}
\item[1.] Suppose $n$ is even.
Then we obtain:
		\begin{enumerate}
		\item[(a)]
If $-n/2< \lambda \le n/2$ then  \[\rho_{y,-r}(a_{\lambda}^{-1}h^{-1}Y_{z,n}(u)h a_{\lambda})=(0,0).\]
		\item[(b)] 
If $\lambda=-n/2$ then  \[\rho_{y,-r}(a_{\lambda}^{-1}h^{-1}Y_{z,n}(u)h a_{\lambda})=(0,y)\]  with  ${\sgn}(y)={\sgn}(u){\sgn}(\ac{-c^2+a^2\eps})$.
		\end{enumerate}	
We see that ${\mathcal N}_{s_1,n,\nu}^{1,\pm}$ are both zero when $n$
is even, for any value of $\nu$. 
In order to find ${\mathcal N}_{s_1,n,\nu}^{2,\pm}$, 
we need to calculate the ratio
of the volume of the subset of $\GKq_{x}$ defined by the formula
'$\nexists \eta\neq 0, \ac{-c^2+a^2\eps}=\eta^2$' to the total volume
of $\GKq_{x}$. Note that $-c^2+a^2\eps=-(c^2-a^2\eps)$. (This
formula should be understood as an abbreviation. We should first
consider the formula with an extra free variable $\delta$: 
'$\nexists \eta\neq 0, \ac{-c^2+a^2\delta}=\eta^2$', do the motivic
calculation,
then plug in our value of $\eps$, and the calculation is very similar to the 
one 
in Lemma~\ref{lemma:unram.elements}, carried out in \cite{JG: overview}.)
Note that the expression $c^2-a^2\eps$ is a square of a nonzero
element for exactly half of the elements of $\GKq_{x}$. Then this ratio is $1/2$.
Therefore, in the both cases $\nu=1$ and $\nu=-1$, the answer
is the same: 
$$
{\mathcal N}_{s_1,n,\nu}^{2,\pm}=\frac12[\GKq_xa_{n/2}\GKq_y/\GKq_y]=\frac12\lef^{2|-n/2|}=\frac12\lef^n.
$$
 
 \item[2.] Suppose $n$ is odd.
In this case, we have
		\begin{enumerate}
		\item[(a)]
If $-(n+1)/2<\lambda<(n+1)/2$ then  \[\rho_{y,-1/2}(a_{\lambda}^{-1}h^{-1}Y_{z,n}(u)h a_{\lambda})=(0,0).\]
		\item[(b)] 
If $\lambda=-(n+1)/2$ then the element
$a_{\lambda}^{-1}h^{-1}Y_{z,n}(u)h a_{\lambda}$ is not in
$\GKq_x$ for any $h$.  		
\item[(c)] 
If $\lambda=(n+1)/2$ then \[\rho_{y,-r}(a_{\lambda}^{-1}h^{-1}Y_{z,n}(u)h a_{\lambda})=(x,0),\] with ${\sgn}(x)={\sgn}(u){\sgn}(\ac{d^2-\eps b^2})$.
		\end{enumerate}
Similarly to the previous case, we get:
${\mathcal N}_{s_1,n,\nu}^{2,\pm}=0$; ${\mathcal N}_{s_1,n,\nu}^{1,\pm}=\frac12\lef^n$ (again using
lemma~\ref{lemma:unram.elements} to show that the volume of the subset of $\GKq_x$ defined by
the formula '$\nexists \eta\neq 0, \ac{d^2-b^2\eps}=\eta^2$' equals
half of the volume of $\GKq_x$. 
	\end{enumerate} 

The calculation in the case $z=s_2$ follows the pattern of the case
$z=s_1$, except there are the same additional complications as we saw
in the case of the character. Not surprisingly, the answer is still
the same as in the case of $s_1$ except the roles of the cases $n$
even and $n$ odd are switched. 
The results of similar calculations for the remaining  ramified 
cases are summarized in Tables~\ref{table:m-orb.s} and \ref{table:m-orb.t}.   

The virtual Chow motives ${\mathcal N}_{z,n,\nu}^0$ can be computed following arguments as above. 
However, since they do not appear in any of our further calculations (all the functions $\kappa$ we {are}
interested in vanish at the origin), we omit this calculation.  
\end{proof}

\begin{corollary}\label{cor: ramified orbital}
Let $\h_{z}$ be ramified elliptic (so $z\in \{t_0,t_1,t_2,t_3\}$).
For each $n\in \N$, the Fourier transform $\foi_{X_{z}}$ of the
orbital integral at $X_{z}$ is constant on the set
${}^G\h_{z',n,\nu}$ for $\nu=\pm$, 
{ $z'\in \{s_1,s_2,t_0,t_1,t_2,t_3\}$} and a
positive (resp., non-negative, if 
{$z'\in\{t_0{,}t_1,t_2,t_3\}$)}
integer $n$, and we have  
\begin{equation*}
\begin{aligned}
&\foi_{X_z}(\tilde f_{z',n,\nu})\\
&= \hat{\varphi}_{X_{z}}\vert_{V_0}\frob {\mathcal N}_{z',n,\nu}^{0}
+\sum_{\alpha=\pm}
\hat{\varphi}_{X_{z}}\vert_{V^{1,\nu}}\frob {\mathcal
  N}_{z',n,\nu}^{1,\alpha}
+ \sum_{\alpha=\pm}\hat{\varphi}_{X_{z}}\vert_{V^{2,\nu}}\frob {\mathcal N}_{z',n,\nu}^{2,\alpha},
\end{aligned}
\end{equation*}
\end{corollary}

\begin{proof}
From Section~\ref{subsection: applications}, we see that our elements
$X_{z}$ correspond to the point  $y=(01)$ in the building, and their depth
is $r=1/2$.  
First, recall that by Proposition~\ref{thm:good}, 
we have 
$$
\oi_{X_{z}}(\tilde f_{z',n,\pm})
=\int_G\int_{\gKq}\tilde f_{z',n,\pm}(gYg^{-1})(\hat\varphi_{\bar{X}_{z}})_{y,-r}(Y).
$$
Note that the support of the functions $\tilde f_{z,n,\pm}$ for all
$z$ and all $n>0$ is contained in $\mathfrak{h}_{y,-r}$, so the
Proposition~\ref{thm:good} is applicable. 
Now it remains to plug in the function
$\kappa=\hat\varphi_{\bar{X}_{z}}$, which is now a
${\GL}(1,\kq)$-invariant function on $\A^2(\kq)$  
in the Proposition~\ref{orb.int} to complete the proof.
\end{proof}

\begin{table}[htbp]
\caption{The virtual motives for the Fourier transform of orbital
  integrals at elements $Y_{z,n}(u)$, for {$z\in \{ s_1,s_2\}$}. 
Here $\nu={\sgn}(u)$.} 
\centering
  \begin{tabular}{@{}|  c  |c  |c |  @{}}
\hline
&&\\
	$z$
& ${\mathcal N}_{z,n,\nu}^{1,\pm}$ & ${\mathcal N}_{z,n,\nu}^{2,\pm}$ \\
&&\\
\hline
&&\\
$s_1$
 & $\begin{aligned} 
{\mathcal N}_{z,n,\nu}^{1, +}&=&{\mathcal N}_{z,n,\nu}^{1,-}&=\frac12 \lef^n &\quad n\text{ odd}\\ 
 {\mathcal N}_{z,n,\nu}^{1, +}&=&{\mathcal N}_{z,n,\nu}^{1,-}&=0 &\quad n\text{ even}
\end{aligned}$
& $\begin{aligned} {\mathcal N}_{z,n,\nu}^{2,+}&=&{\mathcal N}^{2,-}&=0 &\quad n\text{ odd}\\ 
 {\mathcal N}_{z,n,\nu}^{2,+}&=&{\mathcal N}_{z,n,\nu}^{2, -}&=\frac12\lef^{n}
 &\quad n\text{ even}
\end{aligned}$ \\
&&\\
\hline
&&\\
$s_2$&
$\begin{aligned} {\mathcal N}_{z,n,\nu}^{1,+}&=&{\mathcal N}^{1,-}&=0 &\quad n\text{ odd}\\ 
 {\mathcal N}_{z,n,\nu}^{1,+}&=&{\mathcal N}_{z,n,\nu}^{1,-}&=\frac12\lef^{n}
 &\quad n\text{ even}
\end{aligned}$
& $\begin{aligned} 
{\mathcal N}_{z,n,\nu}^{2, +}&=&{\mathcal N}_{z,n,\nu}^{2,-}&=\frac12 \lef^n &\quad n\text{ odd}\\ 
 {\mathcal N}_{z,n,\nu}^{2, +}&=&{\mathcal N}_{z,n,\nu}^{2,-}&=0 &\quad n\text{ even}
\end{aligned}$
\\ 
&&\\
\hline
  \end{tabular}
\label{table:m-orb.s}  
\end{table}

\begin{table}[htbp]
\caption{The virtual motives for the Fourier transform of orbital
  integrals at elements $Y_{z,n}(u)$, for {$z\in \{
  t_0,t_1,t_2,t_3\}$}. (Recall that $\zeta^2={\sgn}(-1)$ and $\nu={\sgn}(u)$.)} 
\centering
  \begin{tabular}{@{}|  c  |c  |c |  @{}}
\hline
&&\\
	$z$
& ${\mathcal N}_{z,n,\nu}^{1,\pm}$ & ${\mathcal N}_{z,n,\nu}^{2,\pm}$ \\
&&\\
\hline
&&\\
$t_0$
&$\begin{aligned}
 &{\mathcal
    N}_{z,n,\nu}^{1,+}=\frac{\lef^{n+1}}{\lef+1}\quad\text{and}\quad
  {\mathcal N}_{z,n,\nu}^{1,-}=0, \\ 
& \text{if } \zeta^{2n+2} =\nu \\ 
& {\mathcal N}_{z,n,\nu}^{1,+}= 0\quad\text{and}\quad
{\mathcal N}_{z,n,\nu}^{1,-}=\frac{\lef^{n+1}}{\lef+1}\\
&\text{otherwise }
\end{aligned}$ 
&{$\begin{aligned}
&{\mathcal N}_{z,n,\nu}^{2,+}=\frac{\lef^{n+1}}{\lef+1}
    \quad\text{and}\quad  {\mathcal N}_{z,n,\nu}^{2,-}=0 \\ 
&\text{if } \zeta^{2n+2} =\nu\\ 
&{\mathcal N}_{z,n,\nu}^{2,+}= 0\quad\text{and}\quad {\mathcal
      N}_{z,n,\nu}^{2,-}=\frac{\lef^{n+1}}{\lef+1} \\
 & \text{otherwise}\end{aligned}$ }\\
&&\\
\hline
&&\\
$t_1$
& $\begin{aligned}
& {\mathcal
    N}_{z,n,\nu}^{1,+}=\frac{\lef^{n+1}}{\lef+1}\quad\text{and}\quad
  {\mathcal N}_{z,n,\nu}^{1,-}=0, \\  
&\text{if } \zeta^{2n+2} =- \nu \\ 
 &{\mathcal N}_{z,n,\nu}^{1,+}= 0\quad\text{and} \quad 
{\mathcal N}_{z,n,\nu}^{1,-}=\frac{\lef^{n+1}}{\lef+1}\\
&\text{otherwise }
\end{aligned}$ 
&{$\begin{aligned}
&{\mathcal N}_{z,n,\nu}^{2,+}=\frac{\lef^{n+1}}{\lef+1}\quad\text{and} \quad {\mathcal
  N}_{z,n,\nu}^{2,-}=0, \\ 
&\text{if } \zeta^{2n+2} =-\nu\\ 
&{\mathcal N}_{z,n,\nu}^{2,+}= 0\quad\text{and}\quad {\mathcal
      N}_{z,n,\nu}^{2,-}=\frac{\lef^{n+1}}{\lef+1}\\
&\text{otherwise}\end{aligned}$ }\\
&&\\
\hline
&&\\
$t_2$
&$\begin{aligned}
 &{\mathcal N}_{z,n,\nu}^{1,+}=\frac{\lef^{n+1}}{\lef+1}
  \quad\text{and}\quad 
{\mathcal N}_{z,n,\nu}^{1,-}=0, \\ 
& \text{if } \zeta^{2n+2}
=(-1)^{n+1} \nu\\ 
 & {\mathcal N}_{z,n,\nu}^{1,+}= 0 \quad\text{and}\quad
{\mathcal N}_{z,n,\nu}^{1,-}=\frac{\lef^{n+1}}{\lef+1}\\
&\text{otherwise}
\end{aligned}$ 
&{$\begin{aligned}
&{\mathcal
      N}_{z,n,\nu}^{2,+}=\frac{\lef^{n+1}}{\lef+1}\quad\text{and}\quad{\mathcal N}_{z,n,\nu}^{2,-}=0, \\ 
& \text{if } \zeta^{2n+2}
=(-1)^{n} \nu\\ 
& {\mathcal N}_{z,n,\nu}^{2,+}= 0\quad\text{and}\quad {\mathcal
  N}_{z,n,\nu}^{2,-}=\frac{\lef^{n+1}}{\lef+1}\\
&\text{otherwise}\end{aligned}$ }\\
&&\\
\hline
&&\\
$t_3$
&$\begin{aligned}
 &{\mathcal N}_{z,n,\nu}^{1,+}=\frac{\lef^{n+1}}{\lef+1}\quad \text{and} \quad {\mathcal
    N}_{z,n,\nu}^{1,-}=0, \\  
&\text{if } \zeta^{2n+2}
=(-1)^{n} \nu\\ 
& {\mathcal N}_{z,n,\nu}^{1,+}= 0 \quad \text{and} \quad
{\mathcal N}_{z,n,\nu}^{1,-}=\frac{\lef^{n+1}}{\lef+1}\\
&\text{otherwise }
\end{aligned}$ 
&{$\begin{aligned}
&{\mathcal
      N}_{z,n,\nu}^{2,+}=\frac{\lef^{n+1}}{\lef+1}\quad\text{and}\quad
    {\mathcal N}_{z,n,\nu}^{2,-}=0, \\
& \text{if } \zeta^{2n+2} 
=(-1)^{n+1}\nu\\ 
&{\mathcal N}_{z,n,\nu}^{2,+}= 0\quad\text{and} \quad {\mathcal
  N}_{z,n,\nu}^{2,-}=\frac{\lef^{n+1}}{\lef+1} \\
&\text{otherwise}\end{aligned}$ }
\\
&&\\
\hline
  \end{tabular}
\label{table:m-orb.t} 
\end{table}

\section{Motivic proof of the Character formula}\label{section: motivic_proof}

Now we are ready to prove Theorem~\ref{thm:C}, and to calculate the
 coefficients that appear in Table~\ref{table:C} in the process.
  Recall that semi-simple character expansion is an equality of two
 distributions on the topologically nilpotent regular set in the Lie
 algebra. These distributions are represented by locally
 integrable functions, which are constant on the sets 
${}^G\mathfrak{h}_{z,n,\pm}$, see Proposition \ref{prop:m-ch},
 Corollary~\ref{cor: unramified orbital} and
 Corollary~\ref{cor: ramified orbital}.  
We prove the semi-simple character expansion by checking the equality on
 each of these sets.
 
Let $\pi$ be a representation as in
 Section~\ref{subsection:representations}.
We start with the case when $\pi$ is obtained from a
Deligne-Lusztig representation, where the proof is straightforward,
and requires no consideration of separate cases. 
If  $\pi$ is of the type $\pi(x,\theta)$ with $x=(0)$ or $x=(1)$, then
by Proposition~\ref{prop:m-ch}, we have, for each 
$z'\in{\{s_0,s_1,s_2,t_0,t_1,t_2,t_3\}}$
\begin{equation*}\label{ch_as_sum:D-L}
 \frac 1{m({}^G\Gamma_{z',n,\nu})}\Theta_{\pi(x,\theta)}(f_{z',n,\nu})
 = -\sum_{U}Q_T\vert_U \frob M_{z',n,\nu}^{x,U}.
\end{equation*} 
On the other hand, for $z=s_1$ or $s_2$, by
Corollary~\ref{cor: unramified orbital}, 
$$\foi_{X_{z}}(\tilde f_{z',n,\pm})=
\sum_{U} \varphi_{\bar{X}_{z}}\frob M_{z',n,\pm}^{x,U},
$$
where $x=(0)$ if  $z=s_1$, and $x=(1)$ if $z=s_2$.
As shown in Sections \ref{subsub:orbital-s1} and \ref{subsub:orbital-s2}, 
$$Q_T=(1-q)\hat\varphi_{\bar X_{s_1}}=(1-q)\hat\varphi_{\bar X_{s_2}}.$$
It follows immediately that 
\begin{equation}\label{eq:Q_T orbital}
\frac 1{m({}^G\mathfrak {h}_{z',n,\nu})}\Theta_{\pi}(\tilde f_{z',n,\nu}\circ \cay)=
\begin{cases}
(1-q)\foi_{X_{s_1}}(\tilde f_{z',n,\nu}) &\quad \text{if } \pi
=\pi(0,\theta)\\
(1-q)\foi_{X_{s_2}}(\tilde f_{z',n,\nu}) &\quad \text{if }
\pi=\pi(1,\theta)\end{cases}
\end{equation}
This proves the theorem in the case when $\pi$ comes from a
Deligne-Lusztig representation.

Let us now turn to the non-Deligne-Lusztig case. 
Let { $\pi=\pi(x,+)$, where $x=(0)$ or $x=(1)$} 
(the other two cases can be obtained by changing the
sign in front of $Q_G$ everywhere below). Then
\begin{equation} 
\begin{aligned}
&\frac 1{m({}^G\mathfrak {h}_{z',n,\nu})}\Theta_{\pi}(\tilde
f_{z',n,\nu}\circ \cay)\\
&= -\frac12\sum_{U}Q_T(U) \frob M_{z',n,\nu}^{x,U}-\frac12\sum_U Q_G(U)\frob M_{z',n,\nu}^{x,U}.
\end{aligned}
\end{equation}
As we have seen in the equation (\ref{eq:Q_T orbital}), 
the first term (i.e. the part of the character
that comes from inflating $Q_T$) is a multiple of  
$\foi_{X_{s_1}}(\tilde f_{z',n,\nu})$ in the case $\pi=\pi(0,+)$, and of 
$\foi_{X_{s_2}}(\tilde f_{z',n,\nu})$ in the case $\pi=\pi(1,+)$.
The coefficient is $\frac{q-1}{2}$ since on the left $Q_T$ appears with the
coefficient $-\frac{1}{2}$.

It remains to express the second term,  
 $-\frac12\sum_U Q_G(U)\frob M_{z',n,\nu}^{x,U}$, as a linear
 combination of the Fourier transforms of orbital integrals.
In order to do that, recall the functions $\varphi_{(0)}$ and 
$\varphi_{(1)}$ defined Section~\ref{subsection: two functions}.   
By Corollary~\ref{cor: ramified orbital}, if $z'$ is elliptic, we have 
\begin{equation}\label{ecf.orbital side.1}
\begin{aligned}
&\sum_{l=1,2\atop{\alpha=\pm}}\hat\varphi_{(0)}\vert_{V^{l,\alpha}}
\frob {\mathcal N}_{z',n,\nu}^{l,\alpha}\\
&= \foi_{X_{t_0}}(f_{z',n,\nu})-\foi_{X_{t_1}}(f_{z',n,\nu})
+\foi_{X_{t_2}}(f_{z',n,\nu})-\foi_{X_{t_3}}(f_{z',n,\nu})
\end{aligned}
\end{equation}
and
\begin{equation}\label{ecf.orbital side.2}
\begin{aligned}
&\sum_{l=1,2\atop{\alpha=\pm}}\hat\varphi_{(1)}\vert_{V^{l,\alpha}}
\frob {\mathcal N}_{z',n,\nu}^{l,\alpha} \\
&= \foi_{X_{t_0}}(f_{z',n,\nu})-\foi_{X_{t_1}}(f_{z',n,\nu})-\foi_{X_{t_2}}(f_{z',n,\nu})+\foi_{X_{t_3}}(f_{z',n,\nu}).
\end{aligned}
\end{equation}
Also, if $z'=s_0$, then, since the both functions $\hat\varphi_{(0)}$
and $\hat\varphi_{(1)}$ vanish at the origin and take the opposite
values on the orbits $V^{1,+}, V^{1,-}$, and  $V^{2,+}, V^{2,-}$, 
it is easy to see in a way that mimics the proof of 
Proposition ~\ref{prop:m-ch} in the case {$z=s_0$}, that the right-hand side of
equation~(\ref{ecf.orbital side.1}) and (\ref{ecf.orbital side.2}) vanishes on ${\h_{s_0,n}}$ for $n>0$. 

We claim that for any $z$ and any positive (resp., non-negative) $n$,   
\begin{equation}\label{eq:ecf.mot'}
-\frac12\sum_U Q_G(U)\frob M_{z',n,\nu}^{x,U} =c\sum_{l=1,2\atop{\alpha=\pm}}\hat\varphi_{x}\vert_{V^{l,\alpha}}
\frob {\mathcal N}_{z,n,\nu}^{l,\alpha},
\end{equation}
with some constant $c$ that we will calculate below (we will see that 
$c=-\frac{q^2-1}{2^3q}$).
Note that the left-hand side depends on the choice of the vertex
$x$. On the right, it is the constant $c$ and the function $\varphi_x$
that depend on $x$.  
This  equation could be called the motivic version of our character
formula.
It is the core of the proof  -- here we are comparing the inflation of two
functions that live on different reductive quotients.
Note that once we prove this claim, Theorem~\ref{thm:C} will follow
immediately.
On the other hand, the proof of the claim is automatic: all we need to
do is plug in the values of the functions $\varphi_{x}$ from
Section~\ref{subsection: two functions}, the values of $Q_G$ from
Section~\ref{subsection: cuspidal},  and the motivic coefficients from Tables~\ref{table:m-char.0.s} and \ref{table:m-char.0.t} if $x=(0)$ and Tables~\ref{table:m-char.1.s} and \ref{table:m-char.1.t} if $x=(1)$ on the left,
and from Tables~\ref{table:m-orb.s} and \ref{table:m-orb.t} on the right.   

The equality  (\ref{eq:ecf.mot'}) has to be
checked on each of the sets $\mathfrak {h}_{{z'},n,\nu}$
(recall that $\nu=\pm$).

Observe that the function $\hat\varphi_{(0)}$ vanishes at $V^{1,\pm}$,
and the function  $\hat\varphi_{(1)}$ vanishes at $V^{2,\pm}$, so that
the right-hand side in any case has only two nonzero terms.
On the left, since $Q_G$ vanishes at the identity, there are also only
two nonzero terms.

For ${z'}=s_0$ the equality we want to prove is trivial, since the both
sides vanish, as discussed above. 
For ${z'=s_1}$ and ${z=s_2}$, 
the equality also turns out to be  trivial.
On the right, the two nonzero opposite values of the function
$\hat\varphi_{(0)}$ or $\hat\varphi_{(1)}$ appear with the same
coefficient, so the right-hand side equals zero. 
 On the left,  the function $Q_G$ takes opposite values  on
  $U_1$ and $U_{\eps}$, and they also appear on the
    left-hand side with equal coefficients, so the left-hand side of 
equation~(\ref{eq:ecf.mot'}) is also zero. This, however, gives no
information of the constant $c$. 

For {$z'\in \{t_0,t_1,t_2, t_3\}$}, 
we see from Tables~\ref{table:m-orb.s} and \ref{table:m-orb.t} that 
only one of the nonzero values of $\hat\varphi_x$ appears with a nonzero
coefficient, and we see from Tables~\ref{table:m-char.0.s}, \ref{table:m-char.0.t}, \ref{table:m-char.1.s} and
\ref{table:m-char.1.t} that also only one
nonzero term appears on the left.
The left-hand side of equation~(\ref{eq:ecf.mot'})
equals $-\frac12\sqrt{q}\zeta^3q^n\text{SIGN}_1$,
where $\text{SIGN}_1=1$ if the coefficient 
$M_{{z'},n,\nu}^{x,U_1}$ is nonzero, and $\text{SIGN}_1=-1$
if  $M_{{z'},n,\nu}^{x,U_{\eps}}$ is nonzero.

The right-hand side equals:
$$\begin{cases}
c\frac{2^2}{q-1}\sqrt{q}\zeta^{3} \frac{q^{n+1}}{q+1}\text{SIGN}_2&\quad
\text{ if } x=(0),\\ 
c\frac{2^2}{q-1}\sqrt{q}\zeta \frac{q^{n+1}}{q+1}\text{SIGN}_2
&\quad \text{ if } x=(1),
\end{cases}
$$ 
where 
$\text{SIGN}_2$ is the sign that depends on which one of the
virtual Chow motives ${\mathcal N}^{1,2,\pm}_{z',n,\nu}$ 
is nonzero.
Let  
$c=-\frac{q^2-1}{2^3q}{\zeta^2}$, which is the constant that appears 
in Table~\ref{table:C}.
Then, to finish the proof of the theorem it remains to show that 
on every set $\mathfrak {h}_{{z'},n,\nu}$ with 
{$z'\in \{t_0,\dots, t_3\}$},
we have the identity
$\text{SIGN}_1=\zeta^2\text{SIGN}_2$ in the case $x=(1)$, and
$\text{SIGN}_1=\text{SIGN}_2$ in the case $x=(0)$.
Here is the comparison of the two SIGNs in the case $x=(0)$:
\begin{equation*} 
\begin{aligned}
&z=t_0:&\quad \text{SIGN}_1&=\nu\zeta^{2n}\quad 
&\text{SIGN}_2&=\zeta^{2(n+1)}.\\
&z=t_1:&\quad \text{SIGN}_1&=-\nu\zeta^{2n}\quad 
&\text{SIGN}_2&=-\nu\zeta^{2(n+1)}.\\
&z=t_2:&\quad \text{SIGN}_1&=\nu\zeta^{2n}(-1)^n\quad 
&\text{SIGN}_2&=\nu\zeta^{2(n+1)}(-1)^n.\\
&z=t_3:&\quad \text{SIGN}_1&=\nu\zeta^{2n}(-1)^{n+1}\quad 
&\text{SIGN}_2&=\nu\zeta^{2(n+1)}(-1)^{n+1}.
\end{aligned}
\end{equation*}
We see that in all cases $\text{SIGN}_1=\zeta^2\text{SIGN}_2$, which
completes the proof in the case $x=(0)$. 
 
The case $x=(1)$ is identical, except we need to use Tables~\ref{table:m-char.1.s} and \ref{table:m-char.1.t}
instead of Tables~\ref{table:m-char.0.s} and \ref{table:m-char.0.t} to calculate $\text{SIGN}_1$, and
the first column of Tables~\ref{table:m-orb.s} and \ref{table:m-orb.t} instead of the second
column, to calculate $\text{SIGN}_2$.

\section{Final comments}\label{subsection: comments}

\subsection{Theorem~\ref{thm:C} and our choices.}\label{section:choices}

Let us begin by reviewing all the choices made in this paper before the proof of Theorem~\ref{thm:C}.
 
First, in the pre-amble to Section~\ref{section: basic notions}, we began with an odd prime $p$, a $p$-adic field $\Kq$, a prime $\ell$ different from $p$ (\eg  $\ell=2$), and an algebraic closure $\EE$ of $\Q_\ell$. 

Next, in Section~\ref{subsection: gauss} we fixed an additive character $\bar\psi : \kq \to \EE$ and a square root $\sqrt{q}$ of $q$ in $\EE$. These choices (via Gauss sums) determined a fourth-root of unity $\zeta\in \EE$ such that $\zeta^2={\sgn}(-1)$ (see Remark~\ref{remark: zeta}). This, in turn, determined how we labelled the two representations $\sigma_+$ and $\sigma_-$ in the Lusztig series for $(T,\theta_0)$ and therefore our definition of $Q_\ksch{G}$ (see Section~\ref{subsection: cuspidal}), and therefore our definition of $\pi(0,+)$, $\pi(0,-)$, $\pi(1,+)$ and $\pi(1,-)$ (see Remark~\ref{remark: signs}).

Independently, we fixed cocycles $\{s_1,s_2, t_0,t_1,t_2,t_3\}$ in $Z^1(\Kq,N)$ such that their cohomology classes lay in the kernel of the map $H^1(\Kq,N)\to H^1(\Kq,G)$ induced by inclusion $N\to G$. This choice determined { a uniformizer $\varpi$ for $\Kq$ and} a non-square unit $\eps$ in $\Rq$ (see Remark~\ref{remark: eps}). We remind the reader that, if ${\sgn}(-1)=-1$ then the cohomology class for $t_0$ equal the cohomology class for $t_1$ and the cohomology class for $t_2$ equal the cohomology class for $t_3$.

Finally, for each  $z\in \{ s_1,s_2,t_0,t_1,t_2,t_3\}$ we chose an element $X_z\in \gKq$ with minimal non-negative depth in its Cartan subalgebra; they are listed in Table~\ref{table:orbits}. This last step amounted to the choice of a unit $v$ in $\Rq$ (see Section~\ref{section:outline}).

Also, the Fourier transform of the orbital integral $\oi_{X_z}$ is taken with respect to an additive character $\psi:\Kq\to \EE$ with conductor $\Rq$ such that the induced additive character of $\kq$ is $\bar\psi$. Also, we must use compatible Killing forms, and the correct measures everywhere, as we did.

With all these choices we made, the values of the coefficients $c_z(\pi)$ in our semi-simple character expansion are presented in Table~\ref{table:C}. { From that table we make three observations. 
\begin{enumerate}
\item
If $\pi$ is induced from a Deligne-Lusztig representation then, for each cocycle $z$, $c_z(\pi)$ is a rational function of $q$ with integer coefficients.
\item
If the cocycle $z$ is unramified (by which we mean the Cartan $T_z$ is unramified) then $c_z(\pi)$ is a rational function in $q$ with integer coefficients, for every depth zero supercuspidal irreducible representation $\pi$.
\item
If $\pi$ is depth zero supercuspidal irreducible representation but $\pi$ is not induced from a Deligne-Lusztig representation and if the cocycle $z$ is ramified (by which we mean the Cartan $T_z$ is ramified), then $c_z(\pi)$ is a rational number for every $p$, but it is \emph{not} a rational function of $q$. Instead, in this case, $c_z(\pi)$ is a rational function of $q$ multiplied by $\sgn_q(-1)$. Observe that $\sgn_q(-1)$ cannot be expressed as a rational function in $q$.
\end{enumerate}
}
It is reasonable to ask if these properties continue to hold had we made different choices above. { The answer is affirmative. To see why,} we say a few words about how the semi-simple character expansion would change if the definition of the $X_z$s were modified. Since we have gathered complete information about the Fourier transforms of all regular elliptic orbital integrals $\foi_{X}$ (when $X$ has minimal non-negative depth in its Cartan subalgebra) evaluated at topologically nilpotent elements $Y$, we can explore the dependence of the coefficients in the semi-simple character expansion of the orbits we chose. (Of course, each $X_z$ can be replaced by any element in the same adjoint orbit as $X_z$, since the distributions $\foi_{X_z}$ and $\Theta_\pi$ are invariant under the group action.)

To begin, consider the local constancy of $X\mapsto \foi_{X}(f)$ for a fixed Schwartz function $f$ supported by topologically nilpotent elements. Using equations (\ref{eqn:s1.1}) and (\ref{eq:frobenius_char}) we see that  
\begin{equation}\label{eqn:s1.2}
\int_{\GKq} \int_{\gKq} f(\Ad(g)Y)\, (Q_\ksch{T})_{(0),0}(Y)\, dY\, dg = (1-q) \foi_{X_{s_1}(v)}(f).
\end{equation}
Moreover, since $Q_\ksch{T}$ (see Section~\ref{subsection: cuspidal}) does not depend on $v\in \Uq$ it follows that $\foi_{X_{s_1}(v)}(f)$ is independent of $v$; thus, for any $f$ as above, 
\begin{equation}\label{eqn:s1.3}
\int_{\Uq} \foi_{X_{s_1}(v)}(f)\, dv = (q-1)\foi_{X_{s_1}}(f),
\end{equation}
where the Haar measure on $\Uq$ is chosen so that $\Uq$ has measure $q-1$.

From our point of view, it would have been more natural to re-write equation~(\ref{eqn:s1.2}) (and therefore the first line in Table~\ref{table:C}) in the form
\begin{equation}\label{eqn:s1.4}
\Theta_{\pi(0,\theta)}(f) = \int_{\Uq} \foi_{X_{s_1}(v)}(f)\, dv,
\end{equation}
from which we see that our choice for $v$ was completely unimportant here. (See Sections~\ref{subsection:representations} and \ref{subsection: cuspidal} for the definition of $\pi(0,\theta)$.) Similarly, we can re-write the last line of Table~\ref{table:C} in the form
\begin{equation}\label{eqn:s2.4}
\Theta_{\pi(1,\theta)}(f) = \int_{\Uq} \foi_{X_{s_2}(v)}(f)\, dv.
\end{equation}
These expressions depend on our choice for the cocycles $s_1$ and $s_2$ only. Therefore, once the cocycle $s_1$ is chosen, the element $X_{s_1}$ may be replaced by \emph{any} element in the Cartan uniquely determined by the cocycle $s_1$ as long as that element has minimal non-negative depth in that Cartan; likewise, for the cocycle $s_2$.

Regarding the totally ramified Cartan subgroups, things are a bit more subtle. From Table~\ref{table:C} we see that
if $f$ is a Schwartz function supported by topologically nilpotent elements, then
\begin{equation*}
\begin{aligned}
& \hskip-20pt \int_{\GKq} \int_{\gKq} f(\Ad(g)Y)\, (Q_\ksch{G})_{(01),\frac{1}{2}}(Y)\, dY\, dg \\
&= \frac{q^2-1}{2^2q} \left( \foi_{X_{t_0}}(f) - \foi_{X_{t_1}}(f)
+ \foi_{X_{t_2}}(f) - \foi_{X_{t_3}}(f) \right).
\end{aligned}
\end{equation*}
Invoking the local constancy of $X\mapsto\foi_{X}(f)$ we may re-write the equation above in the form
\begin{equation*}
\begin{aligned}
& \hskip-20pt \frac{2q}{q+1} \int_{\GKq} \int_{\gKq} f(\Ad(g)Y)\, (Q_\ksch{G})_{(0),0}(Y)\, dY\, dg \\
&= \int\limits_{\{ {\sgn}(v) = +1\}} \foi_{X_{t_0}(v)}(f)\, dv 
-  \int\limits_{\{ {\sgn}(v) = +1\}} \foi_{X_{t_1}(v)}(f)\, dv \\
& + \int\limits_{\{ {\sgn}(v) = +1\}} \foi_{X_{t_2}(v)}(f)\, dv 
- \int\limits_{\{ {\sgn}(v) = +1\}} \foi_{X_{t_3}(v)}(f)\, dv.
\end{aligned}
\end{equation*}
Similar observations apply to representations induced from $\GKq_{(1)}$. In this way we see that only ${\sgn}(v)$ was important in our definition of the ramified orbits appearing in Theorem~\ref{thm:C}. More precisely, once the cocycle $t_0$ is chosen, the element $X_{t_0}$ may be replaced by  $X_{t_0}(v)$ as long as ${\sgn}(v)=1$; likewise for the cocycles $t_1$, $t_2$ and $t_3$. 
{ However, if ${\sgn}(v)=-1$ and $z\in \{ t_0,t_1,t_2,t_3\}$, then $\foi_{X_{z}}$ is \emph{not} equal to $\foi_{X_{t_1}(v)}(Y)$, even when restricted to functions supported by topologically nilpotent elements of $\gKq$.}

Having dealt with $v$ we return to our choice of non-square unit $\eps$. From the local constancy of the Fourier transform of the orbital integrals we  have the following immediate consequence. For any Schwartz function $f$ supported by topologically nilpotent elements of $\gKq$,
\begin{equation}
\begin{aligned}
\int\limits_{\{ {\sgn}(\delta)= -1\}} \foi_{\lb\begin{smallmatrix} 0 & 1 \\ \delta \varpi & 0 \end{smallmatrix}\rb}(f)\, d\delta 
& = \frac{q-1}{2} \foi_{X_{t_2}}(f).
\end{aligned}
\end{equation}
Similar observations hold for all our orbits. { Motivated by the motivic integration view of things,} we might have replaced each of our orbital integrals with an integral over $v$ and $\delta$ ($\eps$ is a particular value of $\delta$). In fact, that is exactly what we did when we used motivic integration with parameters in Section~\ref{subsection: anything} and allowed the parameter to vary over the set of all non-squares in the residue field. The result is easily related back to orbital integrals appearing in representation theory because of relations of the following form.
\begin{equation}
\begin{aligned}
\int\limits_{\{ {\sgn}(\delta)= -1\}} \int\limits_{\{ {\sgn}(v)=+1\}} \foi_{\lb\begin{smallmatrix} 0 & v \\ \delta \varpi v & 0 \end{smallmatrix}\rb}(f)\, dv\, d\delta 
& = \frac{(q-1)^2}{2^2} \foi_{X_{t_2}}(f).
\end{aligned}
\end{equation} 
Similar observations hold for all our orbital integrals. 
Consequently, replacing any of our orbital integrals with these `smeared' orbital integrals would have the effect of changing the coefficients in the semi-simple character expansion in very simple ways.

{ In summary, we have shown in this section 
that the three observations made above 
concerning the rationality of $c_z(\pi)$ 
are completely independent of all our choices. 

Also, since we have been switching freely
between polynomials in $q$ and virtual Chow motives, we notice
that the coefficients $c_z(\pi)$ can always be
interpreted as elements of the ring $\mot$, as we do in Table~\ref{table:Cmot}. Indeed, if we replace the
polynomials in $q$ with the corresponding elements of the ring $\mot$
as we have been doing so far, we can see that the denominators of
$c_z(\pi)$ are invertible in the ring
$\mot$. Moreover, a motivic expression for
$\zeta^2$ is found by considering the $0$-dimensional
variety defined by the equation $x^2=-1$ (this variety appears for the
same purpose in \cite{Tom.whatis}): if $\mathbb{S}$ denotes the class of $x^2=-1$
in the ring $\mot$, then $\zeta^2=-\frob(1-\mathbb{S})$. Note, moreover, that this motive appears exactly when we study depth zero supercuspidal representations $\pi$ which are not induced from Deligne-Lusztig representations and when we consider cocycles $z$ for which the corresponding Cartan $T_z$ is ramified. We will have more to say about this phenomenon in Section~\ref{section:endoscopy}.

{
\begin{table}[htbp]  \caption{Motives for the coefficients $c_z(\pi)$ appearing in Theorem~\ref{thm:C}. Here we write $\mathbb{M}$ for $(\mathbb{L}^2-1)(1-\mathbb{S})$ in order to save space.}
  \centering
  \begin{tabular}{@{} |c || c  c |  cc | c c | @{}}
    \hline
  &&&&&&\\
$c_z(\pi)$	& $z=s_1$ &  $z=s_2$ & $z=t_0$ & $z=t_1$ & $z=t_2$ & $z=t_3$ \\ 
 &&&&&&\\
    \hline\hline
 &&&&&&\\
${ \pi=\pi(0,\theta)}$ & $\mathbb{L}-1$ &  $0$ & $0$ & $0$ & $0$ & $0$  \\
 &&&&&&\\
${ \pi= \pi(1,\theta)}$ & $0$  & $\mathbb{L}-1$ &  $0$ & $0$ & $0$ & $0$ \\ 
 &&&&&&\\
 \hline\hline
 &&&&&&\\
${ \pi= \pi(0,+)}$ & $\frac{\mathbb{L}-1}{2}$ &  $0$ &  $+\frac{\mathbb{M}}{2^3\mathbb{L}}$ & $-\frac{\mathbb{M}}{2^3\mathbb{L}}$ & $+\frac{\mathbb{M}}{2^3\mathbb{L}}$ & $-\frac{\mathbb{M}}{2^3\mathbb{L}}$ \\
 &&&&&&\\
${ \pi=\pi(0,-)}$ & $\frac{\mathbb{L}-1}{2}$  & $0$ & $-\frac{\mathbb{M}}{2^3\mathbb{L}}$ & $+\frac{\mathbb{M}}{2^3\mathbb{L}}$ & $-\frac{\mathbb{M}}{2^3\mathbb{L}}$ & $+\frac{\mathbb{M}}{2^3\mathbb{L}}$ \\
 &&&&&&\\
${ \pi=\pi(1,+)}$ & $0$  &  $\frac{\mathbb{L}-1}{2}$ &  $+\frac{\mathbb{M}}{2^3\mathbb{L}}$ & $-\frac{\mathbb{M}}{2^3\mathbb{L}}$ & $-\frac{\mathbb{M}}{2^3\mathbb{L}}$ & $+\frac{\mathbb{M}}{2^3\mathbb{L}}$ \\
 &&&&&&\\
${\pi=\pi(1,-)}$ & $0$  &  $\frac{\mathbb{L}-1}{2}$ & $-\frac{\mathbb{M}}{2^3\mathbb{L}}$ & $+\frac{\mathbb{M}}{2^3\mathbb{L}}$ & $+\frac{\mathbb{M}}{2^3\mathbb{L}}$ & $-\frac{\mathbb{M}}{2^3\mathbb{L}}$ \\
 &&&&&&\\
\hline
  \end{tabular}
  \label{table:Cmot}
\end{table}
}

\subsection{Theorem~\ref{thm:C} and endoscopy}\label{section:endoscopy}

This paper has focused on the motivic nature of the values of characters of depth-zero supercuspidal representations of $p$-adic $\SL(2)$, on the motivic nature of the Fourier transform of some associated orbital integrals, and on the relations between the associated motives in the Chow ring. In particular, the techniques used in the paper make it clear that, once these motives are determined, the character formula of Theorem~\ref{thm:C} admits a proof which could easily be automated. While it is seems promising to illustrate a strategy showing that certain results from local harmonic analysis can be proved algorithmically, this method of proof does not explain some of the striking patterns in Table~\ref{table:C}. In this section we explain these patterns.

First, one must note that, although Table~\ref{table:C} is a $6\times 6$ matrix, the rank of this matrix is $4$. It is easy to understand why the rank of the matrix is at most $5$: it follows immediately from the nilpotent characters expansion that characters of depth zero supercuspidal representations of $\SL(2,\Kq)$ span a space of dimension at most $5$, since that is the number of nilpotent orbits in $\fsl(2,\Kq)$ when the residual characteristic of $\Kq$ is odd. 

But why is the rank of Table~\ref{table:C} exactly $4$?  Going back to Section~\ref{subsection: cuspidal}, observe that characters of cuspidal representations of $\SL(2,\kq)$ are linearly dependent when restricted to unipotent elements. In fact, the set $\{ \trace\sigma_\theta, \trace\sigma_+, \trace\sigma_-\}$ admits exactly one linear relation on unipotent elements: if $g\in \fsl(2,\kq)$ is unipotent, then
\begin{equation}
\trace\sigma_+(g) + \trace\sigma_-(g) = \trace\sigma_\theta(g).
\end{equation}
(This relation is best understood through Lusztig's work, but this would take us too far afield.)
Two linear relations involving characters of depth zero supercuspidal representations on topologically nilpotent elements of $\fsl(2,\Kq)$ follow directly from this observation: if $f$ is supported by topologically nilpotent elements, then
\begin{equation}
\begin{aligned}
\Theta_{\pi(0,+)}(\cay^*f) + \Theta_{\pi(0,-)}(\cay^*f) &= \Theta_{\pi(0,\theta)}(\cay^*f)\\
\end{aligned}
\end{equation}
and
\begin{equation}
\begin{aligned}
\Theta_{\pi(1,+)}(\cay^*f) + \Theta_{\pi(1,-)}(\cay^*f) &= \Theta_{\pi(1,\theta)}(\cay^*f).\\ 
\end{aligned}
\end{equation}
Together with the fact that the Fourier transform of our orbital integrals are linearly independent on functions supported by topologically nilpotent elements (which can be seen using the techniques of \cite[\S1]{CH}), this explains why the rank of Table~\ref{table:C} is exactly $4$.

But a deeper understanding of Table~\ref{table:C} begins with the following observation: if $f$ is supported by topologically nilpotent elements, then
\begin{equation}
\begin{aligned}
\Theta_{\pi(0,\theta)}(\cay^*(f)) + \Theta_{\pi(1,\theta)}(\cay^*f) =& 
(q-1) \left( \foi_{X_{s_1}}(f) + \foi_{X_{s_2}}(f) \right);\\
\end{aligned}
\end{equation}
moreover, the right-hand side is the Fourier transform of the stable distribution
\begin{equation}
\oi_{X_{s_1}}^\st \ceq \oi_{X_{s_1}} + \oi_{X_{s_2}},
\end{equation}
and so it follows from the work of Waldspurger that $\foi_{X_{s_1}}^\st$ is a stable distribution. Thus, the sum of characters on the left-hand side is a stable distribution on the set of topologically nilpotent elements. In fact,
\[
\{ \pi(0,\theta),\pi(1,\theta)\}
\]
is an L-packet (see \cite[\S12]{Lind}). Likewise, from Table~\ref{table:C} we see that if $f$ is supported by topologically nilpotent elements, then
\begin{equation}
\begin{aligned}
&\Theta_{\pi(0,+)}(\cay^*f) + \Theta_{\pi(1,+)}(\cay^*f) +\Theta_{\pi(1,+)}(\cay^*f) + \Theta_{\pi(1,-)}(\cay^*f) \\
&= (q-1) \foi_{X_{s_1}}^\st(f),
\end{aligned}
\end{equation}
which is the same stable distribution appearing above. In fact, 
\[
\{ \pi(0,+),\pi(0,-),\pi(1,+),\pi(1,-)\}
\]
is an L-packet (see \cite[\S12]{Lind}). Thus, if $\pi$ is a depth zero supercuspidal irreducible representation of $\SL(2,\Kq)$ then the L-packet containing $\pi$ has cardinality two if $\pi$ is induced from a Deligne-Lusztig representation; otherwise, the L-packet containing $\pi$ has cardinality four.

Now we are ready to understand the pattern seen in Table~\ref{table:C} through the theory of endoscopy. There are exactly five endoscopic groups $H$ for $\GKq$; four elliptic endoscopic groups, and one non-elliptic endoscopic group. The elliptic endoscopic groups are $\SL(2)$ itself, and three copies of $\U(1)$, one for each quadratic extension of $\Kq$. In Table~\ref{table:Y} we write $\U_\eps(1)$ for the special unitary group splitting over $\Kq(\sqrt{\eps})$, $\U_\varpi(1)$ for the special unitary group splitting over $\Kq(\sqrt{\varpi})$,  and $\U_{\eps\varpi}(1)$  for the special unitary group splitting over $\Kq(\sqrt{\eps\varpi})$. The non-elliptic endoscopic group for $\GKq$ is $\GL(1)$, which plays the role of a Levi subgroup of $\SL(2)$; had we considered non-supercuspidal depth zero representations in this paper, it would have played an important part.

\begin{table}[htbp]  \caption{Elliptic endoscopic groups $H$ for $\SL(2)$ over $\Kq$; one $\SL(2)$-regular element $Y_H$ from each $\hKq = \Lie H$; an image $X_H\in \fsl(2,\Kq)$ under the Langlands-Shelstad map; and the $\kappa$-orbital integral determined by $Y_H$.}
  \centering
  \begin{tabular}{@{} | c || c | c | c | c |  @{}}
    \hline
&&&&\\
$H$ & $\SL(2)$ & $\U_{\varepsilon}(1)$ & $\U_{\varpi}(1)$ & $\U_{\varepsilon\varpi}(1)$ \\ 
&&&&\\
\hline\hline
&&&&\\
$Y_H$ & $\left[\begin{smallmatrix} 0 & 1 \\ \eps & 0 \end{smallmatrix}\right]$ &  $\sqrt{\eps}$ & $\sqrt{\varpi}$ & $\sqrt{\eps\varpi}$ \\
&&&&\\
\hline
&&&&\\
$X_H$ & $X_{s_1}$ &  $X_{s_1}$ & $X_{t_0}$ & $ X_{t_2}$ \\
&&&&\\
\hline
&&&&\\
$\oi^{G,H}_{Y_H}$ & $\oi_{X_{s_1}}^\st$ &  $\oi_{X_{s_1}}^\sgn$ & $\oi_{X_{t_0}}^\sgn$ & $\oi_{X_{t_2}}^\sgn$ \\
&&&&\\
\hline
\end{tabular}
\label{table:Y}
\end{table}%

In Table~\ref{table:Y} we pick one good elliptic element $Y_H$ from the Lie algebra $\hKq$ of each elliptic endoscopic group for $\GKq$.  Each $Y_H\in \hKq$ is $\gKq$-regular; for each $Y_H$ we choose an image $X_H = X_z$ in $\gKq$ from the list of elements appearing in Theorem~\ref{thm:C}, such that $\Delta_{\gKq,\hKq}(X_H,Y_H) =1$, where $\Delta_{\gKq,\hKq}$ is the Langlands-Shelstad transfer factor for the pair $(\gKq,\hKq)$. Each $Y_H \in \hKq$ thus determines a $\kappa$-orbital integral on $\gKq$ according to the formula
\begin{equation}
\oi^{G,H}_{Y_H} = \sum_{X'} \Delta_{\gKq,\hKq}(X',Y_H) \oi_{X'},
\end{equation}
where the sum is taken over adjoint orbits in $\gKq$. In our cases the results are:
\begin{equation*}
\begin{aligned}
\oi^{G,\SL(2)}_{Y_{\SL(2)}} &= \oi_{X_{s_1}} + \oi_{X_{s_2}} = \oi_{X_{s_1}}^\st\\
\oi^{G,\U_\eps(1)}_{Y_{\U_\eps(1)}} &= \oi_{X_{s_1}} - \oi_{X_{s_2}} = \oi_{X_{s_1}}^\sgn\\
\oi^{G,\U_\varpi(1)}_{Y_{\U_\varpi(1)}} &= \oi_{X_{t_0}} - \oi_{X_{t_1}} = \oi_{X_{s_1}}^\sgn\\
\oi^{G,\U_{\eps\varpi}(1)}_{Y_{\U_{\eps\varpi}(1)}} &= \oi_{X_{t_2}} - \oi_{X_{t_3}} = \oi_{X_{t_2}}^\sgn,
\end{aligned}
\end{equation*}
as recorded in Table~\ref{table:Y}. It is now clear from Table~\ref{table:C} that, when restricted to Schwartz functions supported by topologically nilpotent elements, the Fourier transform of these four distributions span exactly the same space spanned by the characters of depth zero supercuspidal representations of $\GKq$ (on the pull-back by the Cayley transform of the same space of functions). In fact, the distributions 
\[
\{ \foi_{X_{s_1}}^\st, \foi_{X_{s_1}}^\sgn, \foi_{X_{s_1}}^\sgn, \foi_{X_{t_2}}^\sgn\}
\]
are linearly independent on the set of Schwartz functions supported by topologically nilpotent elements, and provide a \emph{natural} basis for the 
characters of depth zero supercuspidal representations of $\GKq$ (on the pull-back by the Cayley transform of the same space of functions). The result is Theorem~\ref{thm:Cendo}, which enjoys one advantage over Theorem~\ref{thm:C}: the coefficients are \emph{unique} with the choices made above. Moreover, regarding these choices, the techniques of Section~\ref{section:choices} apply here too, as do all other techniques from this paper.

\begin{theorem}\label{thm:Cendo}
Let $\Kq$ be a $p$-adic field with $p\ne 2$. For each
depth zero supercuspidal representation $\pi$ of $\GKq$ and for
each elliptic endoscopic group $H$ for $\GKq$ there is a good elliptic $Y_H\in \Lie H$ with minimal non-negative depth in $\Lie H$ and a unique rational number $c_H(\pi)$ such that
\begin{equation}\label{eqn:Cendo}
 	\Theta_\pi(\cay^*f)
 	= 
	\sum_{H} c_H(\pi) \foi^{G,H}_{Y_H}(f)
\end{equation}
for all Schwartz functions $f$ supported by topologically nilpotent elements in $\gKq$. Moreover, the coefficients are motivic, in the sense explained in this paper. Motives for the coefficients $c_z(\pi)$ are given in Table~\ref{table:Cendo}.
\end{theorem}

\begin{table}[htbp]  \caption{Motives for the \emph{unique} coefficients $c_H(\pi)$ appearing in Theorem~\ref{thm:Cendo}. }
  \centering
  \begin{tabular}{@{} |c || c | c | c | c |  @{}}
    \hline
&&&&\\
 & $H=SL(2)$ & $H=\U_{\varepsilon}(1)$ & $H=\U_{\varpi}(1)$ & $H=\U_{\varepsilon\varpi}(1)$ \\ 
 &&&&\\
    \hline\hline
 &&&&\\
${ \pi=\pi(0,\theta)}$ & $+\frac{\mathbb{L}-1}{2}$ & $+\frac{\mathbb{L}-1}{2}$ & $0$ & $0$\\
 &&&&\\
 ${\pi= \pi(1,\theta)}$ & $+\frac{\mathbb{L}-1}{2}$ & $-\frac{\mathbb{L}-1}{2}$  & $0$ & $0$\\ 
  &&&&\\
 \hline\hline
 &&&&\\
${ \pi=\pi(0,+)}$ &  $+\frac{\mathbb{L}-1}{2^2}$ & $+\frac{\mathbb{L}-1}{2^2}$ & $+\frac{\mathbb{L}^2-1}{2^3\mathbb{L}}(1-\mathbb{S}) $ & $+\frac{\mathbb{L}^2-1}{2^3\mathbb{L}}(1-\mathbb{S}) $  \\
 &&&&\\
${ \pi=\pi(0,-)}$ &  $+\frac{\mathbb{L}-1}{2^2}$ & $+\frac{\mathbb{L}-1}{2^2}$  & $-\frac{\mathbb{L}^2-1}{2^3\mathbb{L}}(1-\mathbb{S})$ & $-\frac{\mathbb{L}^2-1}{2^3\mathbb{L}}(1-\mathbb{S}) $  \\
 &&&&\\
${ \pi=\pi(1,+)}$ & $+\frac{\mathbb{L}-1}{2^2}$ & $-\frac{\mathbb{L}-1}{2^2}$  & $+\frac{\mathbb{L}^2-1}{2^3\mathbb{L}}(1-\mathbb{S})$ & $-\frac{\mathbb{L}^2-1}{2^3\mathbb{L}}(1-\mathbb{S})$ \\
 &&&&\\
${\pi=\pi(1,-)}$ &  $+\frac{\mathbb{L}-1}{2^2}$ & $-\frac{\mathbb{L}-1}{2^2}$  & $-\frac{\mathbb{L}^2-1}{2^3\mathbb{L}}(1-\mathbb{S})$ & $+\frac{\mathbb{L}^2-1}{2^3\mathbb{L}}(1-\mathbb{S})$ \\
 &&&&\\
 \hline
  \end{tabular}
  \label{table:Cendo}
\end{table}

}

\bibliographystyle{amsplain}

\end{document}